\documentclass[12pt]{article}
\pdfoutput=1

\usepackage{geometry,amsthm,amsmath,amssymb,nicefrac,epsfig,mathrsfs,amscd,verbatim,bbm,enumerate}
\usepackage[fleqn,tbtags]{mathtools}

\usepackage{etoolbox}
\newtoggle{arXiv:v3}

\toggletrue{arXiv:v3}

\allowdisplaybreaks

\theoremstyle{plain}
\newtheorem{theorem}{Theorem}[section]
\newtheorem{corollary}[theorem]{Corollary}
\newtheorem{lemma}[theorem]{Lemma}
\newtheorem{proposition}[theorem]{Proposition}

\newtheorem{definition}[theorem]{Definition}
\theoremstyle{remark}

\numberwithin{equation}{section}

\newcommand{\C}{\mathbb{C}}
\newcommand{\R}{\mathbb{R}}

\newcommand{\N}{\mathbb{N}}

\newcommand{\E}{{\mathbb E}}
\renewcommand{\P}{{\mathbb P}}

\newcommand{\calC}{\mathscr{C}}
\newcommand{\calF}{\mathscr{F}}

\newcommand{\calP}{\mathscr{P}}

\newcommand{\eps}{\varepsilon}

\DeclarePairedDelimiter{\abs}{\lvert}{\rvert}
\DeclarePairedDelimiter{\norm}{\lVert}{\rVert}
\newcommand{\dd}{\mathrm{d}}
\usepackage{bm}
\newcommand{\ud}{\,\mathrm{d}}	
\newcommand{\ds}{\ud s}
\newcommand{\dWs}{\ud W_s}
\newcommand\yesnumber{\refstepcounter{equation}\tag{\theequation}} 
\makeatletter
\newcommand{\pushright}[1]{\ifmeasuring@#1\else\omit\hfill$\displaystyle#1$\fi\ignorespaces}
\makeatother

\usepackage[backend=bibtex,style=numeric,firstinits=true,sorting=nyt,maxbibnames=99,useprefix=true]{biblatex}
\DeclareNameAlias{default}{last-first}
\renewbibmacro{in:}{}
\begin{document}

\title{Convergence in H\"older norms with applications\\ to Monte
Carlo methods in infinite dimensions}

\author{Sonja Cox$^1$,
Martin Hutzenthaler$^2$,
Arnulf Jentzen$^3$, \\
Jan van Neerven$^4$,
and Timo Welti$^3$
\bigskip
\\
\footnotesize{$^1$Korteweg-de Vries Institute for Mathematics,
University of Amsterdam,}
\\
\footnotesize{1090 GE Amsterdam, The Netherlands,
e-mail: S.G.Cox$\,$(at)$\,$uva.nl}
\smallskip
\\
\footnotesize{$^2$Faculty of Mathematics,
University of Duisburg-Essen,}
\\
\footnotesize{45127 Essen, Germany,
e-mail: martin.hutzenthaler$\,$(at)$\,$uni-due.de}
\smallskip
\\
\footnotesize{$^3$Seminar for Applied Mathematics,
ETH Z\"{u}rich,}
\\
\footnotesize{8092 Z\"{u}rich, Switzerland,
e-mail: arnulf.jentzen/timo.welti$\,$(at)$\,$sam.math.ethz.ch}
\smallskip
\\
\footnotesize{$^4$Delft Institute of Applied Mathematics,
Delft University of Technology,}
\\
\footnotesize{2628 CD Delft, The Netherlands,
e-mail: J.M.A.M.vanNeerven$\,$(at)$\,$tudelft.nl}
}

\maketitle

\begin{abstract}
 We show
 that if 
 a sequence of piecewise affine linear processes converges 
 in the strong sense with a
 positive rate to a stochastic process
 which is strongly
 H\"older continuous in time,
 then this sequence converges in the strong sense
 even
 with respect to much stronger H\"older norms
 and the convergence rate is essentially reduced by the H\"older exponent.
 Our first application hereof establishes pathwise convergence rates
 for spectral Galerkin approximations of stochastic partial differential equations.
 Our second application derives strong convergence rates
 of multilevel Monte Carlo approximations of expectations of Banach space valued stochastic processes.
\end{abstract}
\tableofcontents

\section{Introduction}
\label{sec:introduction}

In this article we study convergence rates for general stochastic processes in H\"older norms.
In particular, in the main results of this work
(see Corollary~\ref{cor:hoelder1} and Corollary~\ref{cor:hoelder2} in Subsection~\ref{subsec:holder2} below)
we reveal estimates for uniform H\"older errors of general stochastic processes.
In this introductory section we now sketch these results and thereafter outline several applications of the general estimates, which can be found in subsequent sections of this article
(see Corollary~\ref{cor:hoelder3} in Subsection~\ref{subsec:holder2}, Corollary~\ref{cor:convLocLip} in Subsection~\ref{subsec:convLocLip}, and Corollary~\ref{c:mlmc.conv} in Subsection~\ref{sec:mlmc} below).
To illustrate the key results of this work, we consider the following framework throughout this section.
Let
$ 
  T \in (0,\infty)
$
be a real number,
let
$
  \left( \Omega, \mathscr{F},
  \P \right) 
$
be a probability space,
let
$
  \left( 
    E, \left\| \cdot \right\|_E
  \right)
$
be an $ \R $-Banach space,
and
for every function $ f \colon [ 0, T ] \to E $ and every natural number $ N \in \N = \{ 1, 2, 3, \ldots \} $
let
$ [ f ]_N \colon [0, T] \to E $
be the function which satisfies
for all $ n \in \{ 0, 1, \dots, N - 1 \} $, $ t \in \bigl[ \frac{nT}{N}, \frac{(n+1) T}{N} \bigr] $ that
\begin{equation}
[ f ]_N ( t ) = \bigl( n + 1 - \tfrac{tN}{T} \bigr) \cdot f \bigl( \tfrac{nT}{N} \bigr)  +  \bigl( \tfrac{tN}{T} - n \bigr) \cdot f \bigl( \tfrac{(n+1)T}{N} \bigr)
\end{equation}
(the piecewise affine linear interpolation of
$ f |_{ \{ 0, \nicefrac{ T }{ N } , \nicefrac{ 2 T }{ N } , \ldots, \nicefrac{ ( N - 1 ) T }{ N } , T \} } $,
cf.\ \eqref{eq:lin_intpolation} below).

\begin{theorem}
\label{thm:intorduction}
Assume the above setting.
Then
for all 
$ p \in (1,\infty) $,
$ \eps \in ( \nicefrac{ 1 }{ p }, 1 ] $,
$ \alpha \in [ 0, \eps - \nicefrac{ 1 }{ p } ) $
there exists
$ C \in ( 0, \infty ) $
such that it holds for all
$ \beta \in [ \eps, 1 ] $,
$ N \in \N $
and all
$ ( \mathscr{F}, \left\| \cdot \right\|_E ) $-strongly measurable
stochastic processes
$ X, Y \colon [0,T] \times \Omega \rightarrow E $
with continuous sample paths
that
\begin{equation} \label{eq:thmintro}
\begin{split}
& \Bigl( \E \Bigl[
    \| X - [ Y ]_N \|_{ \calC^{ \alpha }( [0,T], \left\| \cdot \right\|_E )  }^p
\Bigr] \Bigr)^{ \nicefrac{1}{p} }
\\ &
\leq
C  N^\eps  \biggl(
\sup_{ n \in \{ 0, 1, \ldots, N \}}
    \bigl\| X_{ \frac{nT}{N} } - Y_{ \frac{nT}{N} } \bigr\|_{ \mathscr{L}^p( \P; \left\| \cdot \right\|_E ) }
+ N^{ -\beta }
\| X \|_{
\calC^{ \beta }( [0,T], \left\| \cdot \right\|_{  \mathscr{L}^p( \P; \left\| \cdot \right\|_E )  } )
} \biggr).
\end{split}
\end{equation}
\end{theorem}

The H\"older and $ \mathscr{L}^p $-norms in \eqref{eq:thmintro} are to be understood in the usual sense (see Subsection~\ref{notation} below for details).
Theorem~\ref{thm:intorduction} is a direct consequence of the more general result in Corollary~\ref{cor:hoelder2} in Subsection~\ref{subsec:holder2} below, which establishes an estimate similar to \eqref{eq:thmintro} also for the case of non-equidistant time grids. Moreover, Corollary~\ref{cor:hoelder1} in Subsection~\ref{subsec:holder2} provides an estimate similar to \eqref{eq:thmintro} but with
$ \bigl( \E \bigl[
    \| X - Y \|_{ \calC^{ \alpha }( [0,T], \left\| \cdot \right\|_E )  }^p
\bigr] \bigr)^{ \nicefrac{1}{p} } $
instead of
$ \bigl( \E \bigl[
\| X - [ Y ]_N \|_{ \calC^{ \alpha }( [0,T], \left\| \cdot \right\|_E )  }^p
\bigr] \bigr)^{ \nicefrac{1}{p} } $
on the left hand side and
with an appropriate H\"older norm of $ Y $ occurring on the right hand side.
Theorem~\ref{thm:intorduction} has a number of applications in the numerical approximation of stochastic processes,
as the next corollary, Corollary~\ref{cor:hoelder4}, clarifies. Corollary~\ref{cor:hoelder4} follows immediately from Theorem~\ref{thm:intorduction}.

\begin{corollary}
\label{cor:hoelder4}
Assume the above setting,
let 
$ \beta \in ( 0, 1 ] $,
and
let
$ X \colon [0,T] \times \Omega \rightarrow E $
and
$ Y^N \colon [0,T] \times \Omega \rightarrow E $,
$ N \in \N $,
be
$ ( \mathscr{F}, \left\| \cdot \right\|_E ) $-strongly measurable
stochastic processes
with continuous sample paths
which satisfy
for all $ p \in ( 1, \infty ) $
that
$ \forall \, N \in \N \colon Y^N = [ Y^N ]_N $
and
\begin{equation}
\label{eq:assumptions}
\| X \|_{
\calC^{ \beta }( [0,T] , \left\| \cdot \right\|_{  \mathscr{L}^p( \P; \left\| \cdot \right\|_E )  } )
}
+
\sup_{ N \in \N }
    \biggl[ N^\beta \sup_{ n \in \{ 0, 1, \ldots, N \}}
    \bigl\| X_{ \frac{nT}{N} } - Y_{ \frac{nT}{N} }^N \bigr\|_{ \mathscr{L}^p( \P; \left\| \cdot \right\|_E ) }
    \biggr] < \infty.
\end{equation}
Then it holds for all
$ p, \eps \in ( 0, \infty ) $
that
\begin{equation}
\sup_{ N \in \N } 
\Biggl[
N^{ \beta - \eps }
\biggl( \E \biggl[
\sup_{ t \in [ 0, T ] }
\| X_t - Y_t^N \|_E^p
\biggr] \biggr)^{ \nicefrac{1}{p} }
\Biggr] < \infty.
\end{equation}
\end{corollary}

It is assumed in \eqref{eq:assumptions} that a sequence of affine linearly interpolated
$ ( \mathscr{F}, \left\| \cdot \right\|_E ) $-strongly measurable
stochastic processes
$ ( Y^N )_{ N \in \N } $
converges
for every $ p \in ( 1, \infty ) $ in $ \mathscr{L}^p( \P; \left\| \cdot \right\|_E ) $
to an
$ ( \mathscr{F}, \left\| \cdot \right\|_E ) $-strongly measurable
stochastic process $ X $ with a positive rate uniformly on all grid points
and that this process $ X $ admits corresponding temporal H\"older regularity.
Corollary~\ref{cor:hoelder4} then shows that these assumptions are sufficient to obtain convergence for every
$ p \in ( 1, \infty ) $
in the uniform
$ L^p( \P; \left\| \cdot \right\|_{ C( [ 0, T ], \left\| \cdot \right\|_E ) } ) $-norm
with essentially the same rate.
Corollary~\ref{cor:hoelder3} in Subsection~\ref{subsec:holder2} below implies this result as a special case and includes the case of non-equidistant time grids.
Moreover, Corollary~\ref{cor:hoelder3} proves an analogous conclusion for convergence in uniform H\"older norms, where the obtained convergence rate is reduced by the considered H\"older exponent.
Corollary~\ref{cor:eulermethod} below demonstrates how this principle can be applied 
to Euler-Maryuama approximations for stochastic differential equations (SDEs) with globally Lipschitz coefficients.
Arguments related to Corollary~\ref{cor:hoelder3} can be found in Lemma~A1 in Bally, Millet, \& Sanz-Sol\'e~\cite{BallyMilletSanzSole1995} and in the second display on page 325 in~\cite{CoxVanNeerven2013}.

Corollary~\ref{cor:hoelder4} is particularly useful for the study of stochastic partial differential equations (SPDEs).
In general, a solution of an SPDE fails to be a semi-martingale.
As a consequence, Doob's maximal inequality cannot be applied to obtain estimates with respect to the
$ L^2( \P; \left\| \cdot \right\|_{ C( [ 0, T ], \left\| \cdot \right\|_E ) } )$-norm.
However, convergence rates with respect to the
$ C( [ 0, T ], \left\| \cdot \right\|_{ L^2( \P; \left\| \cdot \right\|_E ) } ) $-norm
are often feasible and Corollary~\ref{cor:hoelder4} can then be applied to obtain
convergence rates with respect to the
$ L^2( \P; \left\| \cdot \right\|_{ C( [ 0, T ], \left\| \cdot \right\|_E ) } )$-norm.
Estimates with respect to the
$ L^2( \P; \left\| \cdot \right\|_{ C( [ 0, T ], \left\| \cdot \right\|_E ) } )$-norm
are useful for using standard localisation arguments in order to extend results for SPDEs with globally Lipschitz continuous nonlinearities to results for SPDEs with nonlinearities that are only Lipschitz continuous on bounded sets.
We demonstrate this in Corollary~\ref{cor:convLocLip} in Subsection~\ref{subsec:convLocLip} below in the case of pathwise convergence rates for Galerkin approximations.
To be more specific, Corollary~\ref{cor:convLocLip} proves essentially sharp pathwise convergence rates for spatial Galerkin and noise approximations for a large class of SPDEs with non-globally Lipschitz continuous nonlinearities.
For example, Corollary~\ref{cor:convLocLip} applies to stochastic Burgers, stochastic Ginzburg-Landau, stochastic Kuramoto-Sivashinsky, and Cahn-Hilliard-Cook equations.

Another prominent application of Corollary~\ref{cor:hoelder4} are multilevel Monte Carlo methods in Banach spaces.
For a random variable 
$ X \in \mathscr{L}^2( \P; \left\| \cdot \right\|_E ) $
convergence in $ \mathscr{L}^2( \P; \left\| \cdot \right\|_E ) $
of Monte Carlo approximations of the expectation $ \E[X] \in E $ has only been established if $ E $ has so-called (Rademacher) type $ p $ for some $ p \in ( 1, 2 ] $ and in this case the convergence rate is given by $ 1 - \nicefrac{1}{p} $ (see, e.g., Heinrich~\cite{h01} or Corollary~\ref{cor:monteCarloBanach} in Subsection~\ref{sec:montecarlo} below).
However, the space $ C( [ 0, T ], E ) $ fails to have type $ p $ for any $ p \in ( 1, 2 ] $.
If $ X $ has more sample path regularity, this problem can nevertheless be bypassed.
More precisely, if it holds for some $ \alpha \in ( 0, 1 ] $, $ p \in ( \nicefrac{1}{\alpha}, \infty) $ that
$ X \in \mathscr{L}^2( \P; \left\| \cdot \right\|_{ \mathcal{W}^{ \alpha, p }( [ 0, T ], E ) } ) $,
then Monte Carlo approximations of
$ \E[X] \in \mathcal{W}^{ \alpha, p }( [ 0, T ], E ) $
have been shown to converge in
$ \mathscr{L}^2( \P; \left\| \cdot \right\|_{ \mathcal{W}^{ \alpha, p }( [ 0, T ], E ) } ) $
with rate $ 1 - \nicefrac{1}{ \min\{ 2, p \} } $ and, by the Sobolev embedding theorem,
also converge in
$ \mathscr{L}^2( \P; \left\| \cdot \right\|_{ C( [ 0, T ], \left\| \cdot \right\|_E ) } )$
with the same rate.
Here for any real numbers $ \alpha \in ( 0, 1 ] $, $ p \in ( \nicefrac{1}{\alpha}, \infty) $
we denote by
$ \mathcal{W}^{ \alpha, p } ( [ 0, T ], E ) $
the Sobolev space
with regularity parameter $ \alpha $ and integrability parameter $ p $
of continuous functions from $ [ 0, T ] $ to $ E $.
Informally speaking, in order to gain control over the variances appearing in multilevel Monte Carlo approximations it is therefore sufficient for the approximations to converge with respect to the
$ L^2( \P; \left\| \cdot \right\|_{ \calC^\alpha( [ 0, T ], \left\| \cdot \right\|_E ) } )$-norm
for some $ \alpha \in ( 0, 1] $.
For more details, we refer the reader to Section~\ref{sec:cubature} and, in particular, to Corollary~\ref{c:mlmc.conv},
which formalises this approach for the case of multilevel Monte Carlo approximations of expectations of Banach space valued stochastic processes.

Finally, we mention a few results in the literature
which employ some findings from this article.
In particular, Corollary~\ref{cor:hoelder23} in this article is applied
in the proof of Corollary~6.3 in Jentzen \& Pu\v{s}nik~\cite{JentzenPusnik2015}
to prove
uniform convergence in probability for spatial spectral Galerkin approximations of stochastic evolution equations (SEEs) with semi-globally Lipschitz continuous coefficients (see Proposition~6.4 in Jentzen \& Pu\v{s}nik~\cite{JentzenPusnik2015}).
Moreover,
Corollary~\ref{cor:convGlobLip} in this article is employed
in Subsection~5.2 and Subsection~5.3 in~\cite{CoxHutzenthalerJentzen2013}
for transferring
initial value regularity results
for finite-dimensional stochastic differential equations
to the case of infinite-dimensional SPDEs
using the examples of
the stochastic Burgers equation
and the Cahn-Hilliard-Cook equation.
Furthermore,
Corollary~\ref{cor:hoelder3} in this article is used in the proof of
Corollary~5.2 in Hutzenthaler, Jentzen, \& Salimova~\cite{HutzenthalerJentzenSalimova2018}
to establish
essentially sharp uniform strong convergence rates
for spatial spectral Galerkin approximations
of linear stochastic heat equations.

\subsection{Notation}
\label{notation}
In this subsection we introduce some of the notation which we use throughout this article.
For two sets $ A $ and $ B $ we denote by $ \mathbb{M}( A , B ) $
the set of all mappings from $ A $ to $ B $.
For measurable spaces
$ ( \Omega_1, \calF_1 ) $ and $ ( \Omega_2, \calF_2 ) $
we denote by
$ \mathscr{M}( \calF_1, \calF_2 ) $
the set of all $ \calF_1 / \calF_2 $-measurable mappings from $ \Omega_1 $ to $ \Omega_2 $.
For topological spaces 
$ ( E, \mathscr{E} ) $
and 
$ ( F, \mathscr{F} ) $
we denote by
$ \mathscr{B}( E ) $
the Borel $\sigma$-algebra on $ ( E, \mathscr{E} ) $
and we denote by
$ C( E, F ) $
the set of all continuous functions from $ E $ to $ F $.
We denote by
$ \left| \cdot \right| \colon \R \to [ 0, \infty ) $
the absolute value function on $ \R $.
We denote by
$ \Gamma \colon ( 0 , \infty ) \to ( 0 , \infty ) $
the gamma function, that is,
we denote by
$ \Gamma \colon ( 0 , \infty ) \to ( 0 , \infty ) $
the function which satisfies
for all $ x \in ( 0 , \infty  ) $
that $ \Gamma ( x ) = \int_0^{ \infty } t^{ ( x - 1 ) } \, e^{ - t }  \ud t $.
We denote by
$ \mathscr{E}_r \colon [ 0 , \infty ) \to [ 0 , \infty )$, $ r \in ( 0 , \infty ) $,
the mappings which satisfy
for all $ r \in ( 0 , \infty ) $, $ x \in [ 0 , \infty ) $ that
\begin{equation}
\mathscr{E}_{r} [x] = \Biggl( \sum_{ n = 0 }^{ \infty } \frac{ x^{2n} ( \Gamma (r) )^n }{ \Gamma ( nr + 1 ) } \Biggr)^{ \nicefrac{1}{2} }
=
\biggl(
    1
    + \frac{ x^2 \, \Gamma (r) }{ \Gamma ( r + 1 ) }
    + \frac{ x^4 ( \Gamma (r) )^2 }{ \Gamma ( 2r + 1 ) }
    + \ \ldots \
\biggr)^{ \nicefrac{1}{2} }
\end{equation}
(cf.\ Chapter~7 in Henry~\cite{Henry1981}).
As a notational device
to condense
the statements and proofs of many results in this article
in a mathematically rigorous way,
we next introduce the notion 
of an extendedly semi-normed vector space,
which, roughly speaking,
corresponds to a vector space
with a semi-norm-type function
that is allowed to attain infinity.
For a field $ \mathbb{K} \in \{ \R, \C \} $,
a $ \mathbb{K} $-vector space $ V $,
and a mapping $ \left\| \cdot \right\| \colon V \to [0,\infty] $
which satisfies for all 
$ v, w \in \left\{ u \in V \colon  \| u \| < \infty \right\} $,
$ \lambda \in \mathbb{K} $
that
$ \left\| \lambda v \right\|  =  \sqrt{ [ \mathrm{Re}( \lambda ) ]^2 + [ \mathrm{Im}( \lambda ) ]^2 } \left\| v \right\| $
and 
$ \left\| v + w \right\|  \leq \left\| v \right\| + \left\| w \right\| $
we call $ \left\| \cdot \right\| $ an extended semi-norm on $ V $
and we call $ ( V , \left\| \cdot \right\| ) $ an extendedly semi-normed vector space.
For a measure space $ ( \Omega, \mathscr{F}, \mu ) $,
a measurable space $ ( S, \mathscr{S} ) $,
a set $ R \subseteq S $,
and a function $ f \colon \Omega \to R $
we denote by
$ [ f ]_{ \mu, \mathscr{S} } $
the set given by
\begin{equation}
[ f ]_{ \mu, \mathscr{S} }
= \bigl\{  g \in \mathscr{M}( \mathscr{F}, \mathscr{S} ) \colon ( \exists A \in \mathscr{F} \colon \mu(A) = 0
                \text{ and }  \{ \omega \in \Omega \colon f(\omega) \neq g(\omega) \} \subseteq A )  \bigr\}.
\end{equation}
For a measure space $ ( \Omega, \mathscr{F}, \mu ) $,
a normed vector space $ ( V , \left\| \cdot \right\|_V ) $,
and real numbers $ p \in [ 0, \infty ) $, $ q \in ( 0, \infty ) $
we denote by $ \mathscr{L}^0( \mu; \left\| \cdot \right\|_V ) $
the set given by
\begin{equation}
\mathscr{L}^0( \mu; \left\| \cdot \right\|_V )
= \bigl\{  f \in \mathbb{M}( \Omega, V ) \colon
                f \text{ is } ( \mathscr{F}, \left\| \cdot \right\|_V ) \text{-strongly measurable} \bigr\},
\end{equation}
we denote by
$ \left\| \cdot \right\|_{ \mathscr{L}^q( \mu; \left\| \cdot \right\|_V ) } \colon
    \mathscr{L}^0( \mu; \left\| \cdot \right\|_V ) \to [0,\infty] $
the mapping which satisfies
for all $ f \in \mathscr{L}^0( \mu; \left\| \cdot \right\|_V ) $ that
\begin{equation}
\left\| f \right\|_{ \mathscr{L}^q( \mu; \left\| \cdot \right\|_V ) }
= \left[  \int_{ \Omega }
                \left\| f( \omega ) \right\|_V^q   \mu( \dd \omega )
   \right]^{ \nicefrac{ 1 }{ q } }
   \in [0,\infty],
\end{equation}
we denote by
$ \mathscr{L}^q( \mu; \left\| \cdot \right\|_V ) $
the set given by
\begin{equation}
\mathscr{L}^q( \mu; \left\| \cdot \right\|_V ) 
= \Bigl\{  f \in \mathscr{L}^0( \mu; \left\| \cdot \right\|_V ) \colon
    \left\| f \right\|_{ \mathscr{L}^q( \mu; \left\| \cdot \right\|_V ) } 
    < \infty \Bigr\},
\end{equation}
we denote by
$ L^p ( \mu; \left\| \cdot \right\|_V ) $
the set given by
\begin{equation}
L^p ( \mu; \left\| \cdot \right\|_V )
= \bigl\{  \{ g \in \mathscr{L}^0( \mu; \left\| \cdot \right\|_V ) \colon
                    \mu( f \neq g ) = 0 \} \subseteq \mathscr{L}^0( \mu; \left\| \cdot \right\|_V ) \colon
    f \in \mathscr{L}^p( \mu; \left\| \cdot \right\|_V )  \bigr\},
\end{equation}
and we denote by
$ \left\| \cdot \right\|_{ L^q( \mu; \left\| \cdot \right\|_V ) } \colon
    L^0( \mu; \left\| \cdot \right\|_V ) \to [0,\infty] $
the function which satisfies
for all $ f \in \mathscr{L}^0( \mu; \left\| \cdot \right\|_V ) $ that
\begin{equation}
\bigl\| \{ g \in \mathscr{L}^0( \mu; \left\| \cdot \right\|_V ) \colon
                    \mu( f \neq g ) = 0 \} \bigr\|_{ L^q( \mu; \left\| \cdot \right\|_V ) }
= \left\| f \right\|_{ \mathscr{L}^q( \mu; \left\| \cdot \right\|_V ) } \in [0,\infty].
\end{equation}
Note that
for every $ p \in [ 1, \infty ) $,
every measure space $ ( \Omega, \mathscr{F}, \mu ) $,
and every normed vector space $ ( V , \left\| \cdot \right\|_V ) $
it holds that
$ \bigl( 
    \mathscr{L}^0( \mu; \left\| \cdot \right\|_V ),
    \left\| \cdot \right\|_{ \mathscr{L}^p( \mu; \left\| \cdot \right\|_V ) } 
\bigr) $
and
$ \bigl( 
    L^0( \mu; \left\| \cdot \right\|_V ),
    \left\| \cdot \right\|_{ L^p( \mu; \left\| \cdot \right\|_V ) }
\bigr) $
are extendedly semi-normed vector spaces.
For a real number $ T \in [0,\infty) $,
a measurable space $ ( S, \mathscr{S} ) $,
a normed vector space $ ( V, \left\| \cdot \right\|_V ) $,
and a mapping $ X \colon [ 0, T ] \times S \to V $
which satisfies
for all $ t \in [ 0, T ] $ that
$ X_t \colon S \to V $ is an $ ( \mathscr{S}, \left\| \cdot \right\|_V ) $-strongly measurable mapping
we call $ X $ an
$ ( \mathscr{S}, \left\| \cdot \right\|_V ) $-strongly measurable stochastic process.
For a metric space
$ ( M, d ) $,
an extendedly semi-normed vector space
$ ( E, \left\| \cdot \right\| ) $,
a real number
$ r \in [0,1] $,
and a set
$ A \subseteq (0,\infty) $ 
we denote by
$
  \left| 
    \cdot  
  \right|_{
    \calC^{ r, A }( M, \left\| \cdot \right\| )
  },
  \left| \cdot \right|_{
    \calC^{ r }( M, \left\| \cdot \right\| )
  },
  \left\| \cdot \right\|_{
    C( M, \left\| \cdot \right\| )
  },
  \left\| \cdot \right\|_{
    \calC^{ r }( M, \left\| \cdot \right\| )
  } \colon
  \mathbb{M}( M, E ) \to [0, \infty] $
the mappings which satisfy for all $ f \in \mathbb{M}( M, E ) $ that
\begin{align}
  \left| 
    f  
  \right|_{
    \calC^{ r, A }( M, \left\| \cdot \right\| )
  }
&  
  =
  \sup\!\left(
    \left\{
      \tfrac{
        \left\| f(e_1) - f(e_2) \right\|
      }{
        \left| d(e_1, e_2) \right|^{ r }
      }
      \colon
      e_1, e_2 \in M,
      d(e_1,e_2) \in A
    \right\} \cup \left\{ 0 \right\}
  \right)
  \in [0,\infty],
\\
  \left| f \right|_{
    \calC^{ r }( M, \left\| \cdot \right\| )
  }
&  
  =
  \left| f \right|_{
    \calC^{ r, (0,\infty) }( M, \left\| \cdot \right\| )
  }
 \in [0,\infty],
\\
  \left\| f \right\|_{
    C( M, \left\| \cdot \right\| )
  }
  & 
  =
  \sup\!\left(
    \left\{
      \left\| f(e) \right\|
      \colon 
      e\in M
    \right\}
    \cup
    \{ 0 \}
  \right)
  \in [0,\infty] ,
\\
  \left\| f \right\|_{
    \calC^r( M, \left\| \cdot \right\| )
  }
  & 
  =
  \left\| f \right\|_{ C( M, \left\| \cdot \right\| ) }
  +
  \left| f \right|_{ \calC^r( M, \left\| \cdot \right\| ) }
  \in [0,\infty]
\end{align}
and we denote by 
$ 
  \calC^r( M, \left\| \cdot \right\| ) 
$ 
the set given by
\begin{equation}
\calC^r( M, \left\| \cdot \right\| ) 
=
  \Bigl\{ 
    f \in 
    C(M,E)
    \colon
    \left\| f \right\|_{ \calC^r( M, \left\| \cdot \right\| ) } < \infty
  \Bigr\}.
\end{equation}
For Hilbert spaces
$ ( H_i, \langle \cdot, \cdot \rangle_{H_i}, \left\| \cdot \right\|_{H_i} ) $, $ i \in \{ 1, 2 \} $,
we denote by
$ ( \mathrm{HS}( H_1, H_2 ), \langle \cdot, \cdot \rangle_{\mathrm{HS}( H_1, H_2 )}, $\linebreak$\left\| \cdot \right\|_{\mathrm{HS}( H_1, H_2 )} ) $
the Hilbert space of Hilbert-Schmidt operators from $ H_1 $ to $ H_2 $.
For a real number $ T \in (0,\infty) $
we denote by
$ \calP_T $
the set given by
\begin{equation}
\calP_T
= \bigl\{ \theta \subseteq [0,T] \colon
        \{ 0, T \} \subseteq \theta
        \text{ and } \#( \theta ) < \infty \bigr\}.
\end{equation}
We denote by
$ d_{ \max }, d_{ \min } \colon
    \cup_{ T \in (0,\infty) } \calP_T \to \R $
the functions which satisfy
for all $ \theta = \{ \theta_0 , \theta_1, \ldots, \theta_{ \#( \theta ) - 1 } \} \in \cup_{ T \in (0,\infty) } \calP_T $
with $ \theta_0 < \theta_1 < \ldots < \theta_{ \#( \theta ) - 1 } $
that
\begin{equation}
d_{\max}( \theta )
  = \max_{ j \in \{ 1, 2, \ldots, \#( \theta ) - 1 \} } 
      | \theta_j - \theta_{ j - 1 } |
\qquad \text{and} \qquad
d_{\min}( \theta )
  = \min_{ j \in \{1, 2, \ldots, \#( \theta ) - 1 \} }
      | \theta_j - \theta_{ j - 1 } |.
\end{equation}
For a normed vector space $ ( E , \left\| \cdot \right\|_E ) $,
an element
$ \theta = \{ \theta_0 , \theta_1, \ldots, \theta_{ \#( \theta ) - 1 } \} \in \cup_{ T \in (0,\infty) } \calP_T $
with $ \theta_0 < \theta_1 < \ldots < \theta_{ \#( \theta ) - 1 } $,
and a function
$ f \colon [ 0, \theta_{ \#( \theta ) - 1 } ] \to E $
we denote by
$ [ f ]_{ \theta } \colon $\linebreak
$ [ 0, \theta_{ \#( \theta ) - 1 } ] \to E $ 
the piecewise affine linear interpolation of 
$ f |_{ \{ \theta_0 , \theta_1, \ldots, \theta_{ \#( \theta ) - 1 } \} } $,
that is, 
we denote by
$ [ f ]_{ \theta } \colon
    [ 0, \theta_{ \#( \theta ) - 1 } ] \to E $ 
the function which satisfies
for all
$ j \in \{ 1, 2, \ldots, \theta_{ \#( \theta ) - 1 } \} $, 
$ s \in [ \theta_{ j - 1 } , \theta_j ] $ that
\begin{equation}
\label{eq:lin_intpolation}
[ f ]_{ \theta }( s )
= \frac{ ( \theta_j - s ) f( \theta_{ j - 1 } ) }%
          { ( \theta_j - \theta_{ j - 1 } ) }
  + \frac{ ( s - \theta_{ j - 1 } ) f( \theta_j ) }%
            { ( \theta_j - \theta_{j-1} ) } .
\end{equation}
\section{Convergence  in H\"ol\-der norms for Banach space valued stochastic processes}
\label{sec:holder1}

\subsection{Error bounds for the H\"older norm}
\label{subsec:holder1}

\begin{lemma}[An interpolation-type inequality]
\label{lem:hoelderconv}
Consider the notation in Subsection~\ref{notation},
let $ ( E, \left\| \cdot \right\|_E ) $ 
be a normed vector space,
let $ (M,d) $ be a metric space,
let
$ f \colon M \to E $
be a function,
and let
$ c \in (0,\infty) $,
$ \alpha, \beta, \gamma \in [0,1] $
satisfy $ \alpha \leq \beta \leq \gamma $.
Then 
\begin{equation}
\label{eq:toshow0}
  \left|
    f
  \right|_{
    \calC^{ \beta }(
      M, \left\| \cdot \right\|_E
    )
  }
\leq
  \max\!\left\{
    c^{  \alpha - \beta }
    \left|
      f
    \right|_{
      \calC^{\alpha, (c,\infty)
      }( M, \left\| \cdot \right\|_E )
    }
    ,
    c^{  \gamma - \beta  }
    \left| f \right|_{
      \calC^{ \gamma, (0,c]
      }( M, \left\| \cdot \right\|_E )
    }
  \right\}
\end{equation}
and
\begin{equation}
  \left|
    f
  \right|_{
    \calC^{ \beta }(
      M, \left\| \cdot \right\|_E
    )
  }
\leq
  \max\!\left\{
    c^{  \alpha - \beta  }
    \left|
      f
    \right|_{
      \calC^{\alpha, [c,\infty)
      }( M, \left\| \cdot \right\|_E )
    }
    ,
    c^{  \gamma - \beta  }
    \left| f \right|_{
      \calC^{ \gamma, (0,c)
      }( M, \left\| \cdot \right\|_E )
    }
  \right\} 
  .
\label{eq:toshow1}
\end{equation}
\end{lemma}
\begin{proof}[Proof
of Lemma~\ref{lem:hoelderconv}]
First of all, note that it holds for all $ e_1, e_2 \in M $ with 
$ d(e_1,e_2) \in (c,\infty) $ that
\begin{equation}
\label{eq:use1}
\begin{aligned}
&
  \frac{
    \left\|
        f( e_1 )
      -
        f( e_2 )
    \right\|_{
      E
    }
  }{
    \left| d( e_1, e_2 ) \right|^{ \beta }
  }
\leq
    \left| d( e_1, e_2 ) \right|^{
       \alpha - \beta
    }
    \left|
      f
    \right|_{
      \calC^{\alpha, (c,\infty)
      }( M, \left\| \cdot \right\|_E )
    }
\leq
    c^{ \alpha - \beta }
    \left|
      f
    \right|_{
      \calC^{\alpha, (c,\infty)
      }( M, \left\| \cdot \right\|_E )
    }.
\end{aligned}
\end{equation}
In addition, observe that
it holds for all 
$ e_1, e_2 \in M $
with
$
  d( e_1, e_2 ) \in (0, c]
$
that
\begin{equation}
\label{eq:use2}
\begin{aligned}
  \frac{
    \left\|
      f( e_1 )
      -
      f( e_2 )
    \right\|_{
      E
    }
  }{
    \left| d( e_1, e_2 ) \right|^{ \beta }
  }
& \leq
  \left| d( e_1, e_2 ) \right|^{
     \gamma - \beta
  }
    \left| f \right|_{
      \calC^{ \gamma, (0,c] }( M, \left\| \cdot \right\|_E
      )
    }
\leq
    c^{  \gamma - \beta  }
    \left| f \right|_{
      \calC^{ \gamma, (0,c] }( M, \left\| \cdot \right\|_E
      )
    }.
\end{aligned}
\end{equation}
Combining \eqref{eq:use1}
and \eqref{eq:use2} shows
\eqref{eq:toshow0}.
The proof of~\eqref{eq:toshow1} is analogous.
This finishes the proof of Lemma~\ref{lem:hoelderconv}.
\end{proof}

\begin{lemma}[Approximation error for affine linear interpolation]
\label{lem:estNaffine}
Consider the notation in Subsection~\ref{notation},
let $ T \in (0,\infty) $, $ \theta \in \calP_T $, 
$ \alpha \in [0,1] $,
let $ ( E, \left\| \cdot \right\|_E ) $ 
be a normed vector space, 
and let  
$ 
  f 
  \colon 
  [0,T]
  \rightarrow E
$
be a function.
Then 
\begin{equation}
\label{eq:appPiecewise}
  \left\|
    f - [ f ]_{ \theta } 
  \right\|_{
    C([0,T],\left\| \cdot \right\|_E)
  }
  \leq
  \big|
    \tfrac{d_{\max}(\theta)}{2 }
  \big|^{\alpha}
  \left|
    f
  \right|_{
    \calC^{ \alpha }( [0,T], \left\| \cdot \right\|_E )
  }
  .
\end{equation}
\end{lemma}

\begin{proof}[Proof of Lemma~\ref{lem:estNaffine}]
Throughout this proof let 
$ N \in \N $,
$ \theta_0, \theta_1, \dots, \theta_N \in [0,T] $
be the real numbers which satisfy
$
  0 = \theta_0 < \theta_1 < \ldots < \theta_N = T
$
and 
$
  \theta = \{ \theta_0, \theta_1, \dots, \theta_N \}
$,
let
$ s \in [0,T] \setminus \theta $,
let $ j \in \{ 1, 2, \ldots, N \} $ 
be the natural number
such that
$ s \in ( \theta_{ j - 1 } , \theta_j ) $,
and let
$ 
  g \colon [0,1] \rightarrow \R
$
be the function which satisfies
for all $ u \in [0,1] $ that
$
  g( u ) = (1 - u) \, u^\alpha + u \, ( 1 - u )^\alpha
$.
Observe that
the concavity of
the function
$ [ 0, \infty ) \ni x \mapsto x^\alpha \in \R $
shows
for all $ u \in [ 0, 1 ] $
that
\begin{equation}
\begin{split}
2^\alpha g(u)
& =
( 1 - u ) \, ( 2 u )^\alpha
+ u \, ( 2 ( 1 - u ) )^\alpha
\leq
( (1-u) \, 2 u + u \, 2 ( 1 - u ) )^\alpha
\\ &
= ( 4 u ( 1 - u ) )^\alpha
= ( 1 - ( 2 u - 1 )^2 )^\alpha
\leq 1.
\end{split}
\end{equation}
Note that this proves that
\begin{equation}
\label{eq:f_minus_ftheta}
\begin{aligned}
&
  \left\| f(s) - [ f ]_{ \theta }( s ) \right\|_E
\leq 
  \tfrac{ ( \theta_j - s ) }{ ( \theta_j - \theta_{ j - 1 } ) }
  \left\| f(s) - f( \theta_{j-1} ) \right\|_E
  +
  \tfrac{ ( s - \theta_{ j - 1 } ) }{ ( \theta_j - \theta_{ j - 1 } ) }
  \left\| f(s) - f( \theta_j ) \right\|_E
\\ &
  \leq  
  \tfrac{ ( \theta_j - s ) }{ ( \theta_j - \theta_{ j - 1 } ) }
  \big( 
    s - \theta_{j-1} 
  \big)^{ \alpha }
  \left| f \right|_{
    \calC^{ \alpha }([0,T],\left\| \cdot \right\|_E)
  }
  +
  \tfrac{ ( s - \theta_{ j - 1 } ) }{ ( \theta_j - \theta_{j-1} ) }
  \big( 
    \theta_j - s  
  \big)^\alpha
  \left| f \right|_{
    \calC^\alpha( [0,T], \left\| \cdot \right\|_E )
  }
\\ &
  =
  \big(
    \tfrac{
      ( \theta_j - s ) 
    }{
      ( \theta_j - \theta_{j-1} )
    }
    \big(
      \tfrac{ 
        ( s - \theta_{j-1} ) 
      }{
        ( \theta_j - \theta_{j-1} )
      } 
    \big)^{ \alpha } 
    + 
    \tfrac{ ( s - \theta_{j-1} ) }{ ( \theta_j - \theta_{j-1} ) }
    \big(
      \tfrac{ 
        ( \theta_j - s ) 
      }{
        ( \theta_j - \theta_{j-1} )
      }
    \big)^{ \alpha } 
  \big)
  \big(
    \theta_j - \theta_{j-1}
  \big)^{ \alpha }
  \left| f \right|_{ \calC^{ \alpha }([0,T],\left\| \cdot \right\|_E) }
\\ &
= 
  g\big( 
    \tfrac{ s - \theta_{ j - 1 } }{ \theta_j - \theta_{ j - 1 } } 
  \big)
  \,
  \big(
    \theta_j - \theta_{ j - 1 }
  \big)^\alpha
  \left| f \right|_{
    \calC^\alpha( [0,T] , \left\| \cdot \right\|_E )
  }
\leq
  \big(
    \tfrac{ \theta_j - \theta_{ j - 1 } }{ 2 } 
  \big)^\alpha
  \left| f \right|_{
    \calC^\alpha( [0,T], \left\| \cdot \right\|_E ) 
  }
  .
\end{aligned}
\end{equation}
The proof of Lemma~\ref{lem:estNaffine} is thus completed.
\end{proof}

The next result, Corollary~\ref{cor:HoelderEstFunc},
provides estimates for the H\"{o}lder norm differences 
of two functions by using the difference of the two functions
on suitable grid points.
Corollary~\ref{cor:HoelderEstFunc} 
is a consequence of Lemma~\ref{lem:hoelderconv}
and Lemma~\ref{lem:estNaffine}.

\begin{corollary}
\label{cor:HoelderEstFunc}
Consider the notation in Subsection~\ref{notation},
let $ T \in (0,\infty) $, 
$ \theta \in \calP_T $,
$ \beta \in [0,1] $,
$ \alpha \in [0,\beta] $,
let $ ( E, \left\| \cdot \right\|_E ) $ be a normed 
vector space,
and let $ f, g \colon [0,T] \rightarrow E $
be functions.
Then 
\begin{equation}
\begin{split}
 &
\label{eq:HoederEstFunc}
  \left| f - g \right|_{
    \calC^{ \alpha }( [0,T], \left\| \cdot \right\|_E )
  }
\\ &  
  \leq 
  \tfrac{
    2 
  }{
    |
      d_{ \max }( \theta )
    |^{ \alpha }
  }
  \bigg[
    \sup_{ 
      t \in \theta
    }
    \left\|
      f( 
        t
      )
      -
      g( 
        t
      )
    \right\|_E
  +
  \tfrac{
    |
      d_{ \max }( \theta )
    |^{ \beta }
  }{ 2^{ \beta } }
  \,
  \big(
    | f |_{ \calC^{ \beta }( [0,T], \left\| \cdot \right\|_E) }
    +
    | g |_{ \calC^{ \beta }( [0,T], \left\| \cdot \right\|_E) }
  \big)
  \bigg]
\end{split}
\end{equation}
and
\begin{align}
\label{eq:HoederEstFunc2}
&
  \left\| f - g \right\|_{
    \calC^{ \alpha }( [0,T], \left\| \cdot \right\|_E )
  }
\\ & \nonumber
  \leq 
  \Big[
    \tfrac{ 2 }{
      |
        d_{ \max }( \theta )
      |^{ \alpha }
    }
    + 1
  \Big]
  \bigg[
    \sup_{
      t \in \theta
    }
    \left\|
      f( 
        t
      )
      -
      g( 
        t
      )
    \right\|_E
  +
  \tfrac{
    |
      d_{ \max }( \theta )
    |^{ \beta }
  }{
    2^{ \beta }
  }
  \,
  \big(
    | f |_{
      \calC^{\beta}([0,T],\left\| \cdot \right\|_E)
    }
    +
    | g |_{
      \calC^{\beta}([0,T],\left\| \cdot \right\|_E)
    }
  \big)
  \bigg]
  .
\end{align}
\end{corollary}

\begin{proof}[Proof of Corollary~\ref{cor:HoelderEstFunc}]
Lemma~\ref{lem:hoelderconv}
and the triangle inequality
ensure that 
\begin{align*}
\label{eq:HoelderEsth1}
&  
  \left|
    f - g
  \right|_{
    \calC^{ \alpha }( [0,T], \left\| \cdot \right\|_E )
  }
\\ & \yesnumber
  \leq
  \max\!\left\{
    |
       d_{ \max }( \theta )
    |^{ - \alpha }
    \,
    | f - g |_{
      \calC^{ 0, ( d_{ \max }( \theta ), \infty ) }( [0,T], \left\| \cdot \right\|_E ) 
    }
    ,
    |
      d_{ \max }( \theta )
    |^{ \beta - \alpha }
    \,
    | f - g |_{
      \calC^{ \beta }( [0,T], \left\| \cdot \right\|_E ) 
    }
  \right\}
\\ &
  \leq
  \max\!\left\{ 
    2
    \,
    |
      d_{ \max }( \theta )
    |^{ - \alpha }
    \,
    \| f - g \|_{ C( [0,T], \left\| \cdot \right\|_E ) }
  ,
  |
    d_{ \max }( \theta )
  |^{ \beta - \alpha }
  \left(
    | f |_{
      \calC^{ \beta }( [0,T], \left\| \cdot \right\|_E ) 
    }
    +
    | g |_{
      \calC^{ \beta }( [0,T], \left\| \cdot \right\|_E )
    }
  \right)
  \right\}
  .
\end{align*}
In addition, observe that
Lemma~\ref{lem:estNaffine}
and the triangle inequality
assure that 
\begin{align}\label{eq:HoelderEsth2}
 & \nonumber
    \left\| 
      f - g 
    \right\|_{ C( [0,T], \left\| \cdot \right\|_E ) }
  \leq    
    \big\| 
      f - [ f ]_{ \theta } 
    \big\|_{ C( [0,T], \left\| \cdot \right\|_E ) }
    +
    \big\| 
      [ f ]_{ \theta } - [ g ]_{ \theta } 
    \big\|_{ C( [0,T], \left\| \cdot \right\|_E ) }
    +
    \big\| 
      [ g ]_{ \theta } - g 
    \big\|_{ C( [0,T], \left\| \cdot \right\|_E ) }
\\ &
  \leq
  \sup_{ 
    t \in \theta
  }
    \|
      f( t )
      -
      g( t )
    \|_E
    +
    \big| 
      \tfrac{
        d_{ \max }(\theta)
      }{ 2 }
    \big|^{ \beta }
    \left(
      | f |_{
        \calC^{ \beta }( [0,T], \left\| \cdot \right\|_E )
      }
      +
      | g |_{
        \calC^{ \beta }( [0,T], \left\| \cdot \right\|_E )
      }
    \right).
\end{align}
Inserting \eqref{eq:HoelderEsth2} 
into~\eqref{eq:HoelderEsth1} yields inequality~\eqref{eq:HoederEstFunc}.
Moreover, adding~inequality~\eqref{eq:HoederEstFunc} and~\eqref{eq:HoelderEsth2} 
results in inequality~\eqref{eq:HoederEstFunc2}.
This finishes the proof of Corollary~\ref{cor:HoelderEstFunc}.
\end{proof}

\begin{lemma}
\label{lem:estcagrid}
Consider the notation in Subsection~\ref{notation},
let $( E, \left\| \cdot \right\|_E ) $ 
be a normed vector space, 
let $ T, c \in (0,\infty) $, 
$ \alpha \in [0,1] $,
$ \theta \in \mathscr{P}_T $,
$ N \in \N $,
$ 
  \theta_0, \dots, \theta_N \in [0,T] 
$
satisfy
$  
  0 = \theta_0 < \ldots < \theta_N = T
$
and $ \theta = \{ \theta_0, \dots, \theta_N \} $,
and let 
$ 
  f 
  \colon 
  [0,T]
  \rightarrow E
$
be a function.
Then 
\begin{equation}
\label{eq:grid_linear_interpolation}
  \left|
    [ f ]_{ \theta }
  \right|_{
    \calC^{ \alpha, (0,c]
    }( [0,T], \left\| \cdot \right\|_E )
  }
 \leq
  \tfrac{ 
    c^{ 1 - \alpha } 
  }{ 
    d_{ \min }( \theta ) 
  }
  \big[
    \sup\nolimits_{ j \in \{ 1, 2, \ldots, N \} }
    \| f( \theta_j ) - f( \theta_{ j - 1 } ) \|_E
  \big] 
  .
\end{equation}
\end{lemma}

\begin{proof}[Proof
of Lemma~\ref{lem:estcagrid}]
Observe that it holds
for all $ s, t \in [0,T] $ with $ t - s \in (0, c] $ 
that
\begin{equation}
\begin{aligned}
&
  \frac{ 
    \| [ f ]_{ \theta }( t ) - [ f ]_{ \theta }( s ) \|_E
  }{
    | t-s |^{\alpha}
  }
=
  \frac{ 
    \bigl\|
      \int_{ (s,t) \setminus \theta } ( [ f ]_{ \theta } )'( u ) \ud u 
    \bigr\|_E
  }{
    | t-s |^{ \alpha }
  }
\\ & \leq 
  \frac{ 
    | t - s| 
    \,
    \bigl[
      \sup_{
        u \in 
        (s, t) \setminus 
        \theta
      } 
      \| 
        ( [ f ]_{ \theta } )'(u) 
      \|_E
    \bigr]
  }{
    | t - s |^{ \alpha }
  }
\\ &
  \leq 
  | t - s |^{ 1 - \alpha }
  \left[
    \sup_{
      j \in \{ 
        1, 2, \ldots, N
      \}
    }
    \frac{ 
      \| f( \theta_j ) - f( \theta_{ j - 1 } ) \|_E
    }{
      | \theta_j - \theta_{ j - 1 } |
    }
  \right]
\\ &
  \leq 
  \frac{   
    c^{ 1 - \alpha } 
  }{ 
    d_{ \min }( \theta ) 
  }
  \left[
    \sup_{  
      j \in \{ 1, 2, \ldots, N \} 
    }
    \| f( \theta_j ) - f( \theta_{ j - 1 } ) \|_E
  \right]
  .
\end{aligned} 
\end{equation}
This completes the proof of
Lemma~\ref{lem:estcagrid}.
\end{proof}

\begin{lemma}
\label{lem:est_hoelder_semi_norm}
Consider the notation in Subsection~\ref{notation},
let 
$ ( E, \left\| \cdot \right\|_E ) $ 
be a normed vector space,
let 
$ T \in (0,\infty) $, 
$ \alpha \in [0,1] $,
$ \theta \in \calP_T $, 
and let
$ 
  f \colon [0,T] \rightarrow E
$ 
be a function.
Then 
$
  \left| [ f ]_{ \theta } \right|_{
    \calC^{ \alpha }( [0,T] , \left\| \cdot \right\|_E )
  }
  \leq 
  \left| f \right|_{
    \calC^{ \alpha }( [0,T] , \left\| \cdot \right\|_E )
  }
$.
\end{lemma}

\begin{proof}[Proof of Lemma~\ref{lem:est_hoelder_semi_norm}]
Throughout this proof let 
$ N \in \N $,
$ \theta_0, \theta_1, \dots, \theta_N \in [0,T] $
be the real numbers which satisfy
$
  0 = \theta_0 < \theta_1 < \ldots < \theta_N = T
$
and 
$
  \theta = \{ \theta_0, \theta_1, \dots, \theta_N \}
$
and let
$ n \colon [0,T] \to \N $
and 
$
  \rho \colon [0,T] \to [0,1] 
$
be the functions which satisfy
for all $ t \in [0,T] $ that
\begin{equation}
  n( t ) = \min\!\big\{ 
    k \in \{ 1, 2, \dots, N \} 
    \colon
    t \in [ \theta_{ k - 1 } , \theta_k ]
  \big\}
\qquad 
\text{and}
\qquad
  \rho( t ) = 
  \frac{ 
    t - \theta_{ n(t) - 1 } 
  }{
    \theta_{ n(t) } - \theta_{ n(t) - 1 }
  }
  .
\end{equation}
Note that it holds for all $ t \in [0,T] $ 
that
\begin{equation}
\label{eq:f_theta_rep}
  [ f ]_{ \theta }( t )
  =
  \left( 1 - \rho(t) \right)
  \cdot 
  f( \theta_{ n(t) - 1 } )
  +
  \rho(t) 
  \cdot f( \theta_{ n(t) } )
  =
  f( \theta_{ n(t) - 1 } )
  +
  \rho( t ) \cdot 
  \big(
    f( \theta_{ n(t) } )
    -
    f( \theta_{ n(t) - 1 } )
  \big)
  .
\end{equation}
Hence, we obtain for all $ t_1, t_2 \in [0,T] $
with $ t_1 < t_2 $
and $ n( t_1 ) = n( t_2 ) $
that
\begin{equation}
\begin{split}
\label{eq:compare_f_affine_same_grid}
&
  \left\|
    [ f ]_{ \theta }( t_1 )
    -
    [ f ]_{ \theta }( t_2 )
  \right\|_E
  =
  \bigl\|
    \left[
      \left( 1 - \rho( t_1 ) \right)
      \cdot 
      f( \theta_{ n(t_1) - 1 } )
      +
      \rho( t_1 )
      \cdot 
      f( \theta_{ n(t_1) } )
    \right]
    \\ & \hphantom{= \bigl\|}
    -
    \left[
      \left( 1 - \rho( t_2 ) \right)
      \cdot 
      f( \theta_{ n(t_1) - 1 } )
      +
      \rho( t_2 )
      \cdot 
      f( \theta_{ n(t_1) } )
    \right]
  \bigr\|_E
\\ & =
  \left\|
    \left( \rho( t_2 ) - \rho( t_1 ) \right)
    \cdot 
    f( \theta_{ n(t_1) - 1 } )
    +
    \left( \rho( t_1 ) - \rho( t_2 ) \right)
    \cdot 
    f( \theta_{ n(t_1) } )
  \right\|_E
\\ & =
  \left|
    \rho( t_1 ) - \rho( t_2 ) 
  \right|
  \cdot
  \left\|
    f( \theta_{ n(t_1) - 1 } )
    -
    f( \theta_{ n(t_1) } )
  \right\|_E
\\ &
\leq 
  \left|
    \rho( t_1 ) - \rho( t_2 ) 
  \right|
  \left| f \right|_{
    \calC^{ \alpha }( [0,T], \left\| \cdot \right\|_E )
  }
  \left|
    \theta_{ n(t_1) - 1 } 
    -
    \theta_{ n(t_1) }
  \right|^{ \alpha }
\\ & =
  \left|
    \rho( t_1 ) - \rho( t_2 ) 
  \right|^{ 1 - \alpha }
  \left| f \right|_{
    \calC^{ \alpha }( [0,T], \left\| \cdot \right\|_E )
  }
  \left|
    \left( 
      \rho( t_1 ) - \rho( t_2 ) 
    \right)
    \cdot
    \left(
        \theta_{ n(t_1) } - \theta_{ n(t_1) - 1 }
    \right)
  \right|^{ \alpha }
\\ & \leq
  \left| f \right|_{
    \calC^{ \alpha }( [0,T], \left\| \cdot \right\|_E )
  }
  \left|
    \left( 
      \rho( t_1 ) - \rho( t_2 ) 
    \right)
    \cdot
    \left(
        \theta_{ n(t_1) } - \theta_{ n(t_1) - 1 }
    \right)
  \right|^{ \alpha }
\\ & =
  \left| f \right|_{
    \calC^{ \alpha }( [0,T], \left\| \cdot \right\|_E )
  }
  \left|
      t_1 - \theta_{ n(t_1) - 1 } - \left( t_2 - \theta_{ n(t_1) - 1 } \right)
  \right|^{ \alpha }
\\ & =
  \left| f \right|_{
    \calC^{ \alpha }( [0,T], \left\| \cdot \right\|_E )
  }
  \left|
    t_1 - t_2
  \right|^{ \alpha }
  .
\end{split}
\end{equation}
Moreover, \eqref{eq:f_theta_rep} ensures 
for all $ t_1, t_2 \in [0,T] $
with $ n( t_1 ) < n( t_2 ) $
that
\begin{align}
& \nonumber
  \left\|
    [ f ]_{ \theta }( t_1 )
    -
    [ f ]_{ \theta }( t_2 )
  \right\|_E
=
  \bigl\|
    \left[
      \left( 1 - \rho( t_1 ) \right)
      \cdot 
      f( \theta_{ n(t_1) - 1 } )
      +
      \rho( t_1 )
      \cdot 
      f( \theta_{ n(t_1) } )
    \right]
\\ \nonumber & \hphantom{= \bigl\|}
    -
    \left[
      \left( 1 - \rho( t_2 ) \right)
      \cdot 
      f( \theta_{ n(t_2) - 1 } )
      +
      \rho( t_2 )
      \cdot 
      f( \theta_{ n(t_2) } )
    \right]
  \bigr\|_E
\\ \nonumber & \leq
  \left( 1 - \rho( t_1 ) \right)
  \left( 1 - \rho( t_2 ) \right)
  \left\|
    f( \theta_{ n(t_1) - 1 } )
    -
    f( \theta_{ n(t_2) - 1 } )
  \right\|_E
  +
  \rho( t_1 ) 
  \,
  \rho( t_2 ) 
  \left\|
    f( \theta_{ n(t_1) } )
    -
    f( \theta_{ n(t_2) } )
  \right\|_E
\\ \nonumber &
  +
  \left( 1 - \rho( t_1 ) \right)
  \rho( t_2 ) 
  \left\|
    f( \theta_{ n(t_1) - 1 } )
    -
    f( \theta_{ n(t_2) } )
  \right\|_E
  +
  \rho( t_1 )
  \left( 1 - \rho( t_2 ) \right)
  \left\|
    f( \theta_{ n(t_1) } )
    -
    f( \theta_{ n(t_2) - 1 } )
  \right\|_E
\\ \nonumber & \leq 
  \left| f \right|_{
    \calC^{ \alpha }( [0,T], \left\| \cdot \right\|_E )
  }
  \bigl\{
  \left( 1 - \rho( t_1 ) \right)
  \left( 1 - \rho( t_2 ) \right)
  |
    \theta_{ n(t_1) - 1 } 
    -
    \theta_{ n(t_2) - 1 } 
  |^{ \alpha }
  +
  \rho( t_1 ) 
  \,
  \rho( t_2 ) 
  \,
  |
    \theta_{ n(t_1) } 
    -
    \theta_{ n(t_2) } 
  |^{ \alpha }
\\ & \quad
    +
    \left( 1 - \rho( t_1 ) \right)
    \rho( t_2 ) 
    \,
    |
      \theta_{ n(t_1) - 1 } 
      -
      \theta_{ n(t_2) } 
    |^{ \alpha }
    +
    \rho( t_1 )
    \left( 1 - \rho( t_2 ) \right)
    |
      \theta_{ n(t_1) } 
      -
      \theta_{ n(t_2) - 1 } 
    |^{ \alpha }
  \bigr\}
  .
\end{align}
The concavity of the function 
$
  ( - \infty, 0 ] \ni x \mapsto | x |^{ \alpha } \in \R
$
hence proves 
for all $ t_1, t_2 \in [0,T] $
with $ n( t_1 ) < n( t_2 ) $
that
\begin{align}
& \nonumber
  \left\|
    [ f ]_{ \theta }( t_1 )
    -
    [ f ]_{ \theta }( t_2 )
  \right\|_E
\\ \nonumber & \leq 
  \left| f \right|_{
    \calC^{ \alpha }( [0,T], \left\| \cdot \right\|_E )
  }
  \bigl|
    \left( 1 - \rho( t_1 ) \right)
    \left( 1 - \rho( t_2 ) \right)
    \left(
      \theta_{ n(t_1) - 1 } 
      -
      \theta_{ n(t_2) - 1 } 
    \right)
    +
    \rho( t_1 ) 
    \,
    \rho( t_2 ) 
    \left(
      \theta_{ n(t_1) } 
      -
      \theta_{ n(t_2) } 
    \right)
\\ \nonumber & \quad
    +
    \left( 1 - \rho( t_1 ) \right)
    \rho( t_2 ) 
    \left(
      \theta_{ n(t_1) - 1 } 
      -
      \theta_{ n(t_2) } 
    \right)
    +
    \rho( t_1 )
    \left( 1 - \rho( t_2 ) \right)
    \left(
      \theta_{ n(t_1) } 
      -
      \theta_{ n(t_2) - 1 } 
    \right)
  \bigr|^{ \alpha }
\\ \nonumber & =
  \left| f \right|_{
    \calC^{ \alpha }( [0,T], \left\| \cdot \right\|_E )
  }
  \bigl|
    \left( 1 - \rho( t_1 ) \right)
    \theta_{ n(t_1) - 1 } 
    +
    \rho( t_1 )
    \,
    \theta_{ n(t_1) } 
    -
    \left( 1 - \rho( t_2 ) \right)
    \theta_{ n(t_2) - 1 } 
    -
    \rho( t_2 ) 
    \,
    \theta_{ n(t_2) } 
  \bigr|^{ \alpha }
\\ & =
  \left| f \right|_{
    \calC^{ \alpha }( [0,T], \left\| \cdot \right\|_E )
  }
  \\ \nonumber & \quad \cdot
  \left|
    \left\{ 
      \theta_{ n(t_1) - 1 } 
      +
      \rho( t_1 )
      \left[ 
        \theta_{ n(t_1) } 
        -
        \theta_{ n(t_1) - 1 } 
      \right]
    \right\}
    -
    \left\{ 
      \theta_{ n(t_2) - 1 } 
      +
      \rho( t_2 ) 
      \left[
        \theta_{ n(t_2) } 
        -
        \theta_{ n(t_2) - 1 } 
      \right]
    \right\}
  \right|^{ \alpha }
\\ \nonumber & =
  \left| f \right|_{
    \calC^{ \alpha }( [0,T], \left\| \cdot \right\|_E )
  }
  \left|
    t_1
    -
    t_2
  \right|^{ \alpha }
  .
\end{align}
Combining this and \eqref{eq:compare_f_affine_same_grid} completes the proof 
of Lemma~\ref{lem:est_hoelder_semi_norm}.
\end{proof}

\begin{lemma}[Approximations by piecewise affine linear functions]
\label{lem:hoelderconv2}
Consider the notation in Subsection~\ref{notation},
let 
$ ( E, \left\| \cdot \right\|_E ) $ 
be a normed vector space,
let 
$ T \in (0,\infty) $, 
$ \alpha \in [0,1] $,
$ \beta \in [\alpha,1] $, 
$ \theta \in \calP_T $, 
and let
$ 
  f, g \colon [0,T] \rightarrow E
$ 
be functions.
Then 
\begin{equation}
\label{eq:hoelderconv1}
  \left|
    f - [ g ]_{ \theta }
  \right|_{
    \calC^{ \alpha }( [0,T], \left\| \cdot \right\|_E )
  }
\leq
  \tfrac{ 
    2
    \,
    \left| 
      d_{ \max }( \theta ) 
    \right|^{ 1 - \alpha } 
  }{ 
    d_{ \min }( \theta ) 
  }
  \sup_{ t \in \theta }
    \left\|
      f( t )
      -
      g( t )
    \right\|_E
    +
  2
  \left|
    d_{ \max }( \theta )
  \right|^{ \beta - \alpha }
  \left|
    f
  \right|_{
    \calC^{ \beta }( [0,T], \left\| \cdot \right\|_E )
  }
\end{equation}
and
\begin{equation}
\begin{split}
\label{eq:hoelderconv2}
&
  \left\|
    f - [ g ]_{ \theta }
  \right\|_{
    \calC^{ \alpha }( [0,T], \left\| \cdot \right\|_E )
  }
\\ & 
\leq
  \left(
    \tfrac{ 
      2
      \,
      \left| 
        d_{ \max }( \theta ) 
      \right|^{ 1 - \alpha } 
    }{ 
      d_{ \min }( \theta ) 
    }
    + 1
  \right)
    \sup_{ t \in \theta }
    \left\|
      f( t )
      -
      g( t )
    \right\|_E
  +
  \left(
    \tfrac{
      2
    }{
      \left|
        d_{\max}(\theta)
      \right|^{ \alpha }
    }
    + 
    \tfrac{ 1 }{ 
      2^{ \beta }
    }
  \right)
  \left|
    d_{\max}(\theta)
  \right|^{ \beta }
  \left|
    f
  \right|_{
    \calC^{ \beta }( [0,T], \left\| \cdot \right\|_E )
  }
  .
\end{split}
\end{equation}
\end{lemma}
\begin{proof}[Proof
of Lemma~\ref{lem:hoelderconv2}]
Throughout this proof let 
$ N \in \N $,
$ \theta_0, \theta_1, \dots, \theta_N \in [0,T] $
be the real numbers which satisfy
$
  0 = \theta_0 < \theta_1 < \ldots < \theta_N = T
$
and 
$
  \theta = \{ \theta_0, \theta_1, \dots, \theta_N \}
$.
Note that
Lemma~\ref{lem:hoelderconv} implies that
\begin{equation}\label{eq:split2}
\begin{split}
  \left|
    f - [ g ]_{ \theta }
  \right|_{
    \calC^{ \alpha }( [0,T], \left\| \cdot \right\|_E )
  }
& \leq
  \max\Bigl\{
    \left| 
      d_{ \max }( \theta ) 
    \right|^{ - \alpha }
    \left|
      f - [ g ]_{ \theta }
    \right|_{
      \calC^{ 
        0, ( d_{ \max }( \theta ) , \infty )
      }( [0,T], \left\| \cdot \right\|_E )
    }
  ,
  \\ & \hphantom{\leq\max\Big\{}
   \left| 
     d_{ \max }( \theta ) 
   \right|^{ \beta - \alpha }
    \left| 
      f - [ g ]_{ \theta }
    \right|_{
      \calC^{ 
        \beta, ( 0, d_{ \max }( \theta ) ]
      }( 
        [0,T], \left\| \cdot \right\|_E
      )
    }
  \Bigr\} 
  .
\end{split}
\end{equation}
Next note that Lemma~\ref{lem:estNaffine} 
ensures that
\begin{equation}
\label{eq:splitparttwo}
\begin{aligned}
&
  \left|
    f - [ g ]_{ \theta }
  \right|_{
    \calC^{ 0, ( d_{\max}(\theta),\infty )
    }( [0,T], \left\| \cdot \right\|_E )
  }
  \leq
    2 
    \left\|
      f - [ g ]_{ \theta }
    \right\|_{
      C( [0,T], \left\| \cdot \right\|_E )
    }
\\ & \leq
    2 
    \left\|
      f - [ f ]_{ \theta }
    \right\|_{
      C( [0,T], \left\| \cdot \right\|_E )
    }
    +
    2 
    \left\|
      [ f ]_{ \theta } 
      - 
      [ g ]_{ \theta }
    \right\|_{
      C( [0,T], \left\| \cdot \right\|_E )
    }
\\ & \leq
   2 \,
   \big| 
     \tfrac{  
       d_{ \max }( \theta ) 
     }{ 2 } 
   \big|^{ \beta }
   \left|
     f
   \right|_{
     \calC^{ \beta }( [0,T], \left\| \cdot \right\|_E )
   }
   +
   2 \cdot
   \sup_{ 
     t \in \theta
   }
   \left\|
     f( 
       t
     )
     -
     g(
       t
     )
   \right\|_E
\\ & \leq
   2 \,
   \big| 
     d_{ \max }( \theta ) 
   \big|^{ \beta }
   \left|
     f
   \right|_{
     \calC^{ \beta }( [0,T], \left\| \cdot \right\|_E )
   }
   +
   2
   \cdot 
   \frac{ 
     d_{ \max }( \theta ) 
   }{
     d_{ \min }( \theta )
   }
   \cdot
   \sup_{ 
     t \in \theta
   }
   \left\|
     f( 
       t
     )
     -
     g(
       t
     )
   \right\|_E
  .
\end{aligned}
\end{equation}
Moreover, observe that 
Lemma~\ref{lem:estcagrid} 
and Lemma~\ref{lem:est_hoelder_semi_norm}
imply that
\begin{align}
& \nonumber
  \left| 
    f - [ g ]_{ \theta }
  \right|_{
    \calC^{ \beta, ( 0,d_{\max}(\theta) ]
    }( [0,T], \left\| \cdot \right\|_E )
  }
  \leq
    \left| 
      f - [ f ]_{ \theta }
    \right|_{
      \calC^{ \beta }( [0,T], \left\| \cdot \right\|_E
      )
    }
    +
    \left| 
      [ f - g ]_{ \theta }
    \right|_{
      \calC^{ 
        \beta ,
        ( 0, d_{\max}(\theta) ] 
      }( [0,T], \left\| \cdot \right\|_E )
    }
\\ \nonumber & \leq
      \left| 
        f
      \right|_{
        \calC^{ \beta }( [0,T], \left\| \cdot \right\|_E
        )
      }
      +
      \left| 
        [ f ]_{ \theta }
      \right|_{
        \calC^{ \beta }( [0,T], \left\| \cdot \right\|_E
        )
      }
\\ \label{eq:splitpartone} & \quad
      +
      \tfrac{ 
	\left| d_{\max}(\theta) \right|^{1-\beta}
      }{ d_{\min}(\theta) }
      \left[
        \sup_{
          j \in \{ 1, 2, \ldots, N \} 
        }
        \big\| 
          \big[ 
            f( \theta_j ) - g( \theta_j ) 
          \big]
          -
          \big[ 
            f( \theta_{ j - 1 } ) - g( \theta_{ j - 1 } ) 
          \big]
        \big\|_E
      \right]
\\ \nonumber & \leq
      2
      \left| 
        f
      \right|_{
        \calC^{ \beta }( [0,T], \left\| \cdot \right\|_E
        )
      }
      +
      \frac{ 2 }{
        \left|
          d_{ \max }( \theta )
        \right|^{ \beta }
      }
        \cdot
        \frac{ 
          d_{ \max }( \theta )  
        }{ 
          d_{ \min }( \theta ) 
        }
        \cdot 
        \sup_{
          t \in \theta 
        }
        \left\| 
          f( t ) - g( t ) 
        \right\|_E
      .
\end{align}
Substituting \eqref{eq:splitpartone} and \eqref{eq:splitparttwo} into 
\eqref{eq:split2} proves \eqref{eq:hoelderconv1}. It thus remains 
to prove estimate~\eqref{eq:hoelderconv2}. For this note that 
Lemma~\ref{lem:estNaffine} yields that
\begin{equation}
\label{eq:hoelderconvuse2}
\begin{aligned}
  \left\|
    f 
    - 
    [ g ]_{ \theta }
  \right\|_{
    C( [0,T], \left\| \cdot \right\|_E )
  }
& \leq
  \left\|
    f 
    - 
    [ f ]_{ \theta }
  \right\|_{
    C( [0,T], \left\| \cdot \right\|_E )
  }
  +
  \left\|
    [ f ]_{ \theta }    
    - 
    [ g ]_{ \theta }
  \right\|_{
    C( [0,T], \left\| \cdot \right\|_E )
  }
\\ & \leq
  \big|
    \tfrac{ d_{ \max }( \theta ) }{ 2 }
  \big|^{ \beta }
  \left|
    f
  \right|_{
    \calC^{ \beta }(
      [0,T], \left\| \cdot \right\|_E
    )
  }
  +
  \sup_{
    t \in \theta
  }
    \left\| f( t ) - g( t ) \right\|_E .
\end{aligned}
\end{equation}
Combining
\eqref{eq:hoelderconv1}
and
\eqref{eq:hoelderconvuse2}
shows \eqref{eq:hoelderconv2}.
The proof 
of Lemma~\ref{lem:hoelderconv2}
is thus completed.
\end{proof}

\subsection{Upper error bounds for stochastic processes
with H\"{o}l\-der continuous sample paths}
\label{subsec:holder2}

We now turn to the result announced in the introduction which provides 
convergence of stochastic processes in H\"older norms given convergence 
on the grid points. For this we first recall the Kolmogorov-Chentsov continuity theorem, cf., e.g.,
Revuz \& Yor~\cite[Theorem~I.2.1 and its proof]{RevuzYor1999}.

\begin{theorem}[Kolmogorov-Chentsov continuity theorem]
\label{thm:Kolmogorov}
Consider the notation in Subsection~\ref{notation}.
There exists a function
$
  \Xi = 
  ( \Xi_{ T, p, \alpha, \beta } )_{
    T, p, \alpha, \beta \in \R
  }
  \colon \R^4 \to \R
$
such that
for every
$ T \in [0,\infty) $, $ p \in (1,\infty) $,
$ \beta \in ( \nicefrac{ 1 }{ p } , 1 ] $,
every Banach space 
$ ( E, \left\| \cdot \right\|_E ) $,
every probability space
$
  ( \Omega, \mathscr{F}, \P )
$,
and every 
$ X \in \calC^{ \beta }( [0,T], \left\| \cdot \right\|_{  \mathscr{L}^p( \P; \left\| \cdot \right\|_E )  } ) $
there exists an $ ( \mathscr{F}, \left\| \cdot \right\|_E ) $-strongly measurable stochastic process
$
  Y \colon [0,T] \times \Omega \to E
$
with continuous sample paths
such that it holds for every
$ \alpha \in [ 0, \beta - \nicefrac{ 1 }{ p } ) $
that
\begin{equation}
\begin{split}
&
\Bigl( \E \Bigl[
\| Y \|_{ \calC^{ \alpha }( [0,T], \left\| \cdot \right\|_E )  }^p
\Bigr] \Bigr)^{ \nicefrac{1}{p} }
  \leq 
  \Xi_{ T, p, \alpha, \beta } 
  \,
  \| X \|_{
    \calC^{ \beta }( [0,T] , \left\| \cdot \right\|_{  \mathscr{L}^p( \P; \left\| \cdot \right\|_E )  } ) 
  }
  < \infty
\quad 
\text{and} \\
& \,  \forall \, t \in [0,T] \colon
  \P( X_t = Y_t ) = 1
  .
\end{split}
\end{equation}
\end{theorem}

The next result, Corollary~\ref{cor:hoelder1}, follows directly
from Corollary~\ref{cor:HoelderEstFunc}
(with
$ T = T $,
$ \theta = \theta $,
$ \beta = \gamma $,
$ \alpha = \beta $,
$ E = L^p ( \P; \left\| \cdot \right\|_E ) $,
$ f = ( [ 0, T ] \ni t \mapsto \{ Z \in \mathscr{L}^0( \P; \left\| \cdot \right\|_E ) \colon
\P( Z \neq X_t - Y_t ) = 0 \} \in L^p ( \P; \left\| \cdot \right\|_E ) ) $,
$ g = 0 $
for
$ p \in [1,\infty) $,
$ \beta \in [0,1] $,
$ \gamma \in [ \beta, 1 ] $
and $ ( \mathscr{F}, \left\| \cdot \right\|_E ) $-strongly measurable stochastic processes
$ X, Y \colon [0,T] \times \Omega \rightarrow E $
with
$ \forall \, t \in [ 0, T ] \colon \| X_t - Y_t \|_{ \mathscr{L}^p( \P; \left\| \cdot \right\|_E ) }
< \infty
$
in the notation of Corollary~\ref{cor:HoelderEstFunc})
and 
the Kolmogorov-Chentsov continuity theorem 
(see Theorem~\ref{thm:Kolmogorov} above).

\begin{corollary}[Grid point approximations]
\label{cor:hoelder1}
Consider the notation in Subsection~\ref{notation},
let $ T \in (0,\infty) $, 
$ \theta \in \calP_{T} $,
let
$
  \left( \Omega, \mathscr{F},
  \P \right) 
$
be a probability space,
and let
$
  \left( 
    E, \left\| \cdot \right\|_E
  \right)
$
be a Banach space.
Then 
\begin{enumerate}[(i)]
\item 
it holds for all 
$ p \in [1,\infty) $,
$ \beta \in [0,1] $, 
$ \gamma \in [ \beta , 1 ] $
and all $ ( \mathscr{F}, \left\| \cdot \right\|_E ) $-strongly measurable stochastic processes
$ 
  X, Y \colon [0,T] \times
  \Omega \rightarrow E 
$
that
\begin{align}
\nonumber
&  \left\|
    X - Y
  \right\|_{ 
    \calC^{ \beta }(
      [0,T],
      \left\| \cdot \right\|_{  \mathscr{L}^p( \P; \left\| \cdot \right\|_E )  }
    )
  }
 \leq
 \left(
    2
    \left| 
      d_{\max}( \theta )
    \right|^{ - \beta }
    + 
    1
 \right)
  \Big[
    \sup\nolimits_{ t \in \theta }
    \|
      X_t
      -
      Y_t
    \|_{ \mathscr{L}^p( \P; \left\| \cdot \right\|_E ) }
\\ & \quad
    +
      \left| 
        d_{\max}( \theta ) 
      \right|^{ \gamma }
    \big(
    \left| 
      X 
    \right|_{
      \calC^{ \gamma }( [0,T], \left\| \cdot \right\|_{  \mathscr{L}^p( \P; \left\| \cdot \right\|_E )  } )
    }
    +
    \left| 
      Y 
    \right|_{
      \calC^{ \gamma }( [0,T], \left\| \cdot \right\|_{  \mathscr{L}^p( \P; \left\| \cdot \right\|_E )  } )
    }
    \big)
  \Big]
\end{align}
\item 
\label{item:II_cor_hoelder2}
and it holds for all 
$ p \in (1,\infty) $,
$ \beta \in ( \nicefrac{ 1 }{ p } , 1 ] $,
$ \alpha \in [ 0, \beta - \nicefrac{ 1 }{ p } ) $,
$ \gamma \in [ \beta, 1 ] $
and all $ ( \mathscr{F}, \left\| \cdot \right\|_E ) $-strongly measurable stochastic processes
$ 
  X, Y \colon [0,T] \times \Omega
  \rightarrow E 
$ 
with continuous 
sample paths that
\begin{equation}
\begin{split}
&
\Bigl( \E \Bigl[
    \| X - Y \|_{ \calC^{ \alpha }( [0,T], \left\| \cdot \right\|_E )  }^p
\Bigr] \Bigr)^{ \nicefrac{1}{p} }
\leq
  \Xi_{ T, p, \alpha, \beta }
  \left\|
    X - Y
  \right\|_{ 
    \calC^{ \beta }(
      [0,T],
      \left\| \cdot \right\|_{  \mathscr{L}^p( \P; \left\| \cdot \right\|_E )  }
    )
  } \\
& \leq
  \Xi_{ T, p, \alpha, \beta }
 \left(
      2 
      \left| 
        d_{ \max }( \theta )
      \right|^{ - \beta }
    + 
    1
 \right)
  \Big[
    \sup\nolimits_{ t \in \theta }
    \|
      X_t
      -
      Y_t
    \|_{ \mathscr{L}^p( \P; \left\| \cdot \right\|_E ) }
\\ & \quad
    +
      \left| 
        d_{\max}( \theta ) 
      \right|^{ \gamma }
    \big(
    \left| 
      X 
    \right|_{
      \calC^{ \gamma }( [0,T], \left\| \cdot \right\|_{  \mathscr{L}^p( \P; \left\| \cdot \right\|_E )  } )
    }
    +
    \left| 
      Y 
    \right|_{
      \calC^{ \gamma }( [0,T], \left\| \cdot \right\|_{  \mathscr{L}^p( \P; \left\| \cdot \right\|_E )  } )
    }
    \big)
  \Big]
  .
\end{split}
\end{equation}
\end{enumerate}
\end{corollary}

The next result, Corollary~\ref{cor:hoelder2},
follows directly 
from Lemma~\ref{lem:hoelderconv2}
(with
$ E = L^p ( \P; \left\| \cdot \right\|_E ) $,
$ T = T $,
$ \alpha = \beta $,
$ \beta = \gamma $,
$ \theta = \theta $,
$ f = ( [ 0, T ] \ni t \mapsto \{ Z \in \mathscr{L}^0( \P; \left\| \cdot \right\|_E ) \colon
\P( Z \neq X_t - X_0 ) = 0 \} \in L^p ( \P; \left\| \cdot \right\|_E ) ) $,
$ g = ( [ 0, T ] \ni t \mapsto \{ Z \in \mathscr{L}^0( \P; \left\| \cdot \right\|_E ) \colon
\P( Z \neq [ Y ]_{ \theta }( t ) - X_0 ) = 0 \} \in L^p ( \P; \left\| \cdot \right\|_E ) ) $
for
$ p \in [1,\infty) $,
$ \beta \in [0,1] $,
$ \gamma \in [ \beta , 1 ] $
and $ ( \mathscr{F}, \left\| \cdot \right\|_E ) $-strongly measurable stochastic processes
$ X, Y \colon [0,T] \times \Omega \rightarrow E $
with
$ \sup_{ t \in \theta } \| X_t - Y_t \|_{ \mathscr{L}^p( \P; \left\| \cdot \right\|_E ) }
+ | X |_{ \calC^{ \gamma }( [0,T], \left\| \cdot \right\|_{  \mathscr{L}^p( \P; \left\| \cdot \right\|_E )  } ) }
< \infty
$
in the notation of Lemma~\ref{lem:hoelderconv2})
and the Kolmogorov-Chentsov continuity theorem
(see Theorem~\ref{thm:Kolmogorov} above).

\begin{corollary}[Piecewise affine linear stochastic processes]
\label{cor:hoelder2}
Consider the notation in Subsection~\ref{notation},
let $ T \in (0,\infty) $,
$
  \theta \in \mathscr{P}_T
$,
let
$
  \left( \Omega, \mathscr{F},
  \P \right) 
$
be a probability space,
and let
$
  \left( 
    E, \left\| \cdot \right\|_E
  \right)
$
be a Banach space.
Then 
\begin{enumerate}[(i)]
\item 
it holds 
for all 
$ p \in [1,\infty) $,
$ \beta \in [0,1] $,
$ \gamma \in [ \beta , 1 ] $
and all $ ( \mathscr{F}, \left\| \cdot \right\|_E ) $-strongly measurable stochastic processes
$ 
  X, Y \colon [0,T] \times
  \Omega \rightarrow E 
$
that
\begin{equation}
\begin{aligned}
&  \left\|
    X - [ Y ]_{ \theta }
  \right\|_{ 
    \calC^{ \beta }(
      [0,T],
      \left\| \cdot \right\|_{  \mathscr{L}^p( \P; \left\| \cdot \right\|_E )  }
    )
  }
 \leq
  \Big[
  \tfrac{ 
    2
    \,
    \left| 
      d_{ \max }( \theta ) 
    \right|^{ 1 - \beta } 
  }{ 
    d_{ \min }( \theta ) 
  }
  + 1
  \Big]
    \sup_{ t \in \theta }
    \|
      X_t
      -
      Y_t
    \|_{ \mathscr{L}^p( \P; \left\| \cdot \right\|_E ) }
\\ & \quad
  +
  \left[
    2 \, | d_{ \max }( \theta ) |^{ - \beta }
    + 2^{ - \gamma }
  \right]
  \left| 
    d_{ \max }( \theta ) 
  \right|^{ \gamma }
  \left| 
    X 
  \right|_{
    \calC^{ \gamma }( [0,T], \left\| \cdot \right\|_{  \mathscr{L}^p( \P; \left\| \cdot \right\|_E )  } )
  }
\end{aligned}
\end{equation}
\item 
\label{item:2_cor_hoelder2}
and it holds 
for all 
$ p \in (1,\infty) $,
$ \beta \in ( \nicefrac{ 1 }{ p }, 1 ] $,
$ \alpha \in [ 0, \beta - \nicefrac{ 1 }{ p } ) $,
$ \gamma \in [ \beta, 1 ] $
and all $ ( \mathscr{F}, \left\| \cdot \right\|_E ) $-strongly measurable stochastic processes
$ 
  X, Y \colon [0,T] \times
  \Omega \rightarrow E 
$
with continuous sample paths that
\begin{align}
& 
\Bigl( \E \Bigl[
    \| X - [ Y ]_{ \theta } \|_{ \calC^{ \alpha }( [0,T], \left\| \cdot \right\|_E )  }^p
\Bigr] \Bigr)^{ \nicefrac{1}{p} }
 \leq
  \Xi_{ T, p, \alpha, \beta }
  \Big(
  \Big[
  \tfrac{ 
    2
    \,
    \left| 
      d_{ \max }( \theta ) 
    \right|^{ 1 - \beta } 
  }{ 
    d_{ \min }( \theta ) 
  }
  + 1
  \Big]
    \sup_{ t \in \theta }
    \|
      X_t
      -
      Y_t
    \|_{ \mathscr{L}^p( \P; \left\| \cdot \right\|_E ) }
\nonumber \\ & \quad
  +
  \left[
    2 \, | d_{ \max }( \theta ) |^{ - \beta }
    + 2^{ - \gamma }
  \right]
  \left| 
    d_{ \max }( \theta ) 
  \right|^{ \gamma }
  \left| 
    X 
  \right|_{
    \calC^{ \gamma }( [0,T], \left\| \cdot \right\|_{  \mathscr{L}^p( \P; \left\| \cdot \right\|_E )  } )
  }
  \Big)
  . \label{eq:cor_hoelder2}
\end{align}
\end{enumerate}
\end{corollary}
In \eqref{eq:cor_hoelder2}
in Corollary~\ref{cor:hoelder2} 
we assume beside other assumptions that
$ \alpha $ is strictly smaller than $ \gamma $.
In general, this assumption cannot be omitted.
To give an example, let $ ( \Omega, \calF, \P) $ be a probability space
and let
$
  W \colon [0,1] \times \Omega \to \R
$ 
be
a one-dimensional standard Brownian motion
with continuous sample paths.
Then it clearly holds for all $ p \in [1,\infty) $
that 
$
  \| W \|_{ \calC^{ \nicefrac{1}{2} }( [0,1], \left\| \cdot \right\|_{ \mathscr{L}^p( \P; | \cdot | ) } ) } < \infty
$.
However, the fact that the sample paths of the Brownian motion 
are $ \P $-a.s.\
not $ \nicefrac{ 1 }{ 2 } $-H\"older continuous
(cf., e.g.,
Revuz \& Yor~\cite[Theorem~I.2.7]{RevuzYor1999}
and, e.g.,
Arcones~\cite[Corollary~3.1]{Arcones1995})
ensures that
it holds for all $ \theta \in \mathscr{P}_1 $, $ p \in (0,\infty) $
that
$
\E \bigl[
    \| W - [ W ]_{ \theta } \|_{ \calC^{ \nicefrac{1}{2} }( [0,1], \left| \cdot \right| )  }^p
\bigr]
= \infty
$.
The following corollary is related to
Lemma~A1 in Bally, Millet, \& Sanz-Sol\'e~\cite{BallyMilletSanzSole1995}.

\begin{corollary}[$ \mathscr{L}^p $-convergence in 
H\"older norms for a fixed $ p \in [ 1,\infty) $]
\label{cor:hoelder23}
Consider the notation in Subsection~\ref{notation},
let 
$ 
  T \in (0,\infty)
$,
$ p \in [1,\infty) $,
$
  \beta \in [ 0 , 1 ]
$,
let
$
  \left( \Omega, \mathscr{F},
  \P \right) 
$
be a probability space,
let
$
  \left( 
    E, \left\| \cdot \right\|_E
  \right)
$
be a Banach
space,
and let 
$
  Y^N 
  \colon [0,T] \times
  \Omega \rightarrow E 
$,
$ N \in \N_0 $,
be $ ( \mathscr{F}, \left\| \cdot \right\|_E ) $-strongly measurable stochastic processes
with continuous sample paths
which satisfy
$
  \limsup_{ N \to \infty } 
  | 
    Y^N
  |_{
    \calC^{ \beta }( [0,T], \left\| \cdot \right\|_{  \mathscr{L}^p( \P; \left\| \cdot \right\|_E )  } )
  }
  < \infty
$
and
$
  \forall \, t \in [0,T] \colon
  \limsup_{ N \to \infty }
  \|
    Y^0_t -
    Y_t^N
  \|_{ \mathscr{L}^p( \P; \left\| \cdot \right\|_E ) }
  = 0
$.
Then 
\begin{enumerate}[(i)]
\item 
\label{item:i_a_priori}
it holds that
$
    | Y^0 |_{
      \calC^{ \beta }( [0,T], \left\| \cdot \right\|_{  \mathscr{L}^p( \P; \left\| \cdot \right\|_E )  } ) 
    } 
  \leq 
    \limsup_{ N \to \infty } 
    |
      Y^N
    |_{
      \calC^{ \beta }( [0,T], \left\| \cdot \right\|_{  \mathscr{L}^p( \P; \left\| \cdot \right\|_E )  } ) 
    }
  < 
  \infty
$,
\item 
\label{item:ii_convergence}
it holds for all 
$ 
  \alpha \in [0,1]
  \cap 
  ( - \infty , \beta ) 
$
that
$
  \limsup_{ N \to \infty }
  \|
    Y^0 - Y^N
  \|_{ 
    \calC^{ \alpha }( [0,T], \left\| \cdot \right\|_{  \mathscr{L}^p( \P ; \left\| \cdot \right\|_E )  } )
  }
  = 0
$,
\item 
\label{item:iii_convergence}
and it holds for all 
$ 
  \alpha \in 
  [0,1]
  \cap
  ( - \infty , \beta - \nicefrac{ 1 }{ p } ) 
$
that
\begin{equation}
  \limsup_{ N \to \infty }
  \E \Bigl[
      \| Y^0 - Y^N \|_{ \calC^{ \alpha }( [0,T], \left\| \cdot \right\|_E )  }^p
  \Bigr]
  = 0.
\end{equation}
\end{enumerate}
\end{corollary}
\begin{proof}[Proof of Corollary~\ref{cor:hoelder23}]
Throughout this proof 
let $ \theta^n \in \mathscr{P}_T $,
$ n \in \N $,
be the sequence which satisfies for all
$ n \in \N $ that
$ 
  \theta^n
  = 
  \{ 0, \frac{ T }{ n } , \frac{ 2 T }{ n } , \ldots, 
  \frac{ ( n - 1 ) T }{ n } , T \}
  \in \mathscr{P}_T
$.
Observe that the assumption that
$
  \forall \, t \in [0,T] \colon
  \limsup_{ N \to \infty }
  \|
    Y^0_t -
    Y_t^N
  \|_{ \mathscr{L}^p( \P; \left\| \cdot \right\|_E ) }
  = 0
$
and the assumption that 
$
  \limsup_{ N \to \infty }
  | 
    Y^N
  |_{
    \calC^{ \beta }( [0,T], \left\| \cdot \right\|_{  \mathscr{L}^p( \P; \left\| \cdot \right\|_E )  } )
  }
  < \infty
$
ensure that
\begin{equation}  
\label{eq:XinCbeta}
\begin{split} 
&
    | Y^0 |_{
      \calC^{ \beta }( [0,T], \left\| \cdot \right\|_{  \mathscr{L}^p( \P; \left\| \cdot \right\|_E )  } ) 
    }
    =
    \sup_{ 
      \substack{ 
        s, t \in [0,T], 
      \\ 
        s \neq t 
      }
    }
    \left[
      \tfrac{
        \| Y^0_t - Y^0_s \|_{ \mathscr{L}^p( \P; \left\| \cdot \right\|_E ) }
      }{
        | t - s |^{ \beta } 
      }
    \right]
\\ & =
    \sup_{ 
      \substack{ 
        s, t \in [0,T], 
      \\ 
        s \neq t 
      }
    }
    \left[
      \tfrac{
        \limsup_{ N \to \infty }
        \| 
          ( Y_t^N - Y_s^N )
          +
          ( Y^0_t - Y^N_t ) 
          +
          ( Y^N_s - Y^0_s )
        \|_{ \mathscr{L}^p( \P; \left\| \cdot \right\|_E ) }
      }{
        | t - s |^{ \beta } 
      }
    \right]
  \\ &
    \leq
    \sup_{ 
      \substack{ 
        s, t \in [0,T], 
      \\ 
        s \neq t 
      }
    }
    \limsup_{ N \to \infty }
    \left[
      \tfrac{ 
        \| 
          Y_t^N - Y_s^N 
        \|_{ \mathscr{L}^p( \P; \left\| \cdot \right\|_E ) }
      }{
        | t - s |^{ \beta }
      }
    \right]
    \leq
    \limsup_{ N \to \infty }
    |
      Y^N
    |_{
      \calC^{ \beta }( [0,T], \left\| \cdot \right\|_{  \mathscr{L}^p( \P; \left\| \cdot \right\|_E )  } ) 
    }
    < \infty
    .
  \end{split}     
\end{equation}
This establishes Item~\eqref{item:i_a_priori}.
In the next step we prove Item~\eqref{item:ii_convergence}.
We apply Item~(i) 
in Corollary~\ref{cor:hoelder1} 
to obtain for all
$
  \alpha \in [0, \beta ]
$,
$ n, N \in \N $
that
\begin{align}
& \nonumber \|
    Y^0 - Y^N
  \|_{ 
    \calC^{ \alpha }(
      [0,T],
      \left\| \cdot \right\|_{  \mathscr{L}^p( \P; \left\| \cdot \right\|_E )  }
    )
  }
  \leq
 \left(
    2
    \left| 
      d_{\max}( \theta^n )
    \right|^{ - \alpha }
    + 
    1
 \right)
  \Big[
    \sup\nolimits_{ t \in \theta^n }
    \|
      Y^0_t
      -
      Y^N_t
    \|_{ \mathscr{L}^p( \P ; \left\| \cdot \right\|_E ) }
\\ & \nonumber \quad
    +
      \left| 
        d_{\max}( \theta^n ) 
      \right|^{ \beta }
    \big(
    | 
      Y^0 
    |_{
      \calC^{ \beta }( [0,T], \left\| \cdot \right\|_{  \mathscr{L}^p( \P; \left\| \cdot \right\|_E )  } )
    }
    +
    | 
      Y^N
    |_{
      \calC^{ \beta }( [0,T], \left\| \cdot \right\|_{  \mathscr{L}^p( \P; \left\| \cdot \right\|_E )  } )
    }
    \big)
  \Big]
\\ &
  \leq
 \left(
    \tfrac{ 2 \, T^{ - \alpha } }{ n^{ - \alpha } }
    + 
    1
 \right)
    \sup\nolimits_{ t \in \theta^n }
    \|
      Y^0_t
      -
      Y^N_t
    \|_{ \mathscr{L}^p( \P ; \left\| \cdot \right\|_E ) }
\\ & \nonumber \quad
    +
 \left(
    \tfrac{ 2 \, T^{ \beta - \alpha } }{ n^{ \beta - \alpha } }
    + 
    \tfrac{ T^{ \beta } }{ n^{ \beta } }
 \right)
    \big(
    | 
      Y^0 
    |_{
      \calC^{ \beta }( [0,T], \left\| \cdot \right\|_{  \mathscr{L}^p( \P; \left\| \cdot \right\|_E )  } )
    }
    +
    | 
      Y^N
    |_{
      \calC^{ \beta }( [0,T], \left\| \cdot \right\|_{  \mathscr{L}^p( \P; \left\| \cdot \right\|_E )  } )
    }
    \big)
    .
\end{align}
Item~\eqref{item:i_a_priori}
and the assumption that
$
  \forall \, t \in [0,T] \colon
  \limsup_{ N \to \infty }
  \|
    Y^0_t -
    Y_t^N
  \|_{ \mathscr{L}^p( \P; \left\| \cdot \right\|_E ) }
  = 0
$
hence imply 
for all
$
  \alpha \in [0, \beta ]
$,
$ n \in \N $
that
\begin{align}
& \nonumber  \limsup_{ N \to \infty }
  \|
    Y^0 - Y^N
  \|_{ 
    \calC^{ \alpha }(
      [0,T],
      \left\| \cdot \right\|_{  \mathscr{L}^p( \P; \left\| \cdot \right\|_E )  }
    )
  }
  \leq
 \Big[
    \tfrac{ 2 \, T^{ - \alpha } }{ n^{ - \alpha } }
    + 
    1
 \Big]
 \biggl[
    \limsup_{ N \to \infty }
    \sup_{ t \in \theta^n }
    \|
      Y^0_t
      -
      Y^N_t 
    \|_{ \mathscr{L}^p( \P ; \left\| \cdot \right\|_E ) }
  \biggr]
\\ & \quad
    +
 \Big[
    \tfrac{ 4 \, T^{ \beta - \alpha } }{ n^{ \beta - \alpha } }
    + 
    \tfrac{ 2 \, T^{ \beta } }{ n^{ \beta } }
 \Big]
   \limsup_{ N \to \infty }
    | 
      Y^N
    |_{
      \calC^{ \beta }( [0,T], \left\| \cdot \right\|_{  \mathscr{L}^p( \P; \left\| \cdot \right\|_E )  } )
    }
\\ \nonumber & =
 \Big[
    \tfrac{ 4 \, T^{ \beta - \alpha } }{ n^{ \beta - \alpha } }
    + 
    \tfrac{ 2 \, T^{ \beta } }{ n^{ \beta } }
 \Big]
   \limsup_{ N \to \infty }
    | 
      Y^N
    |_{
      \calC^{ \beta }( [0,T], \left\| \cdot \right\|_{  \mathscr{L}^p( \P; \left\| \cdot \right\|_E )  } )
    }
  < \infty
    .
\end{align}
Hence, we obtain 
for all
$
  \alpha \in [0,1] \cap ( - \infty , \beta )
$
that
\begin{align}
& \nonumber
  \limsup_{ N \to \infty }
  \|
    Y^0 - Y^N
  \|_{ 
    \calC^{ \alpha }(
      [0,T],
      \left\| \cdot \right\|_{  \mathscr{L}^p( \P; \left\| \cdot \right\|_E )  }
    )
  }
=
  \limsup_{ n \to \infty }
  \limsup_{ N \to \infty }
  \|
    Y^0 - Y^N
  \|_{ 
    \calC^{ \alpha }(
      [0,T],
      \left\| \cdot \right\|_{  \mathscr{L}^p( \P; \left\| \cdot \right\|_E )  }
    )
  }
\\ & \leq 
 \biggl[
    \limsup_{ n \to \infty }
    \tfrac{ 4 \, T^{ \beta - \alpha } }{ n^{ \beta - \alpha } }
    + 
    \limsup_{ n \to \infty }
    \tfrac{ 2 \, T^{ \beta } }{ n^{ \beta } }
 \biggr]
   \limsup_{ N \to \infty }
    | 
      Y^N
    |_{
      \calC^{ \beta }( [0,T], \left\| \cdot \right\|_{  \mathscr{L}^p( \P; \left\| \cdot \right\|_E )  } )
    }
  = 0
    .
\end{align}
This shows Item~\eqref{item:ii_convergence}.
It thus remains to establish Item~\eqref{item:iii_convergence}
to complete the proof of Corollary~\ref{cor:hoelder23}.
For this we apply the first inequality in Item~(ii) in Corollary~\ref{cor:hoelder1}
to obtain 
for all
$ r \in ( \nicefrac{ 1 }{ p } , \infty ) \cap ( - \infty , \beta ] $,
$ \alpha \in [ 0, r - \nicefrac{ 1 }{ p } ) $,
$ N \in \N $
that
\begin{align}
\Bigl( \E \Bigl[
    \| Y^0 - Y^N \|_{ \calC^{ \alpha }( [0,T], \left\| \cdot \right\|_E )  }^p
\Bigr] \Bigr)^{ \nicefrac{1}{p} }
  \leq
  \Xi_{ T, p, \alpha, r }
  \,
  \|
    Y^0 - Y^N
  \|_{ 
    \calC^r( [0,T] , \left\| \cdot \right\|_{  \mathscr{L}^p( \P; \left\| \cdot \right\|_E )  } )
  }
  .
\end{align}
This and Item~\eqref{item:ii_convergence} 
imply 
for all
$ r \in ( \nicefrac{ 1 }{ p } , \infty ) \cap ( - \infty , \beta ) $,
$ \alpha \in [ 0, r - \nicefrac{ 1 }{ p } ) $
that
\begin{align}
  \limsup_{ N \to \infty }
  \E \Bigl[
      \| Y^0 - Y^N \|_{ \calC^{ \alpha }( [0,T], \left\| \cdot \right\|_E )  }^p
  \Bigr]
  \leq
  ( \Xi_{ T, p, \alpha, r } )^p
  \,
  \limsup_{ N \to \infty }
  \|
    Y^0 - Y^N
  \|_{ 
    \calC^r( [0,T] , \left\| \cdot \right\|_{  \mathscr{L}^p( \P; \left\| \cdot \right\|_E )  } )
  }^p
  = 0
  .
\end{align}
This establishes Item~\eqref{item:iii_convergence}.
The proof of Corollary~\ref{cor:hoelder23}
is thus completed.
\end{proof}

The next result, Corollary~\ref{cor:hoelder3} below, 
is a consequence from Corollary~\ref{cor:hoelder1}
and Corollary~\ref{cor:hoelder2}.

\begin{corollary}[Convergence rates with respect to H\"older norms]
\label{cor:hoelder3}
Consider the notation in Subsection~\ref{notation},
let 
$ 
  T \in (0,\infty)
$,
$
  p \in (1,\infty)
$,
$
  \beta \in ( \nicefrac{ 1 }{ p } , 1 ]
$,
$ ( \theta^N )_{ N \in \N } \subseteq \calP_T $
satisfy 
$
  \limsup_{ N \to \infty }
  d_{ \max }( \theta^N ) = 0 
$,
let
$
  \left( \Omega, \mathscr{F},
  \P \right) 
$
be a probability space,
let
$
  \left( 
    E, \left\| \cdot \right\|_E
  \right)
$
be a Banach
space,
let 
$
  Y^N 
  \colon [0,T] \times
  \Omega \rightarrow E 
$,
$ N \in \N_0 $,
be $ ( \mathscr{F}, \left\| \cdot \right\|_E ) $-strongly measurable stochastic processes
with continuous sample paths
which satisfy
$ Y_0^0 \in \mathscr{L}^p( \P; \left\| \cdot \right\|_E ) $
and
\begin{equation}
\label{eq:ass_convergencegrid1}
  | Y^0 |_{
    \calC^{ \beta }( [0,T] , \left\| \cdot \right\|_{  \mathscr{L}^p( \P; \left\| \cdot \right\|_E )  } )
  }
  +
  \sup_{ N \in \N }
  \left[
  \left| 
    d_{ \max }( \theta^N )
  \right|^{ - \beta }
  \sup\nolimits_{ 
    t \in \theta^N
  }
  \|
    Y_t^0 -
    Y_t^N
  \|_{ \mathscr{L}^p( \P; \left\| \cdot \right\|_E ) }
  \right]
  < \infty , 
\end{equation}
and assume  
$
\big(
  \big[
    \sup_{ N \in \N }
    | 
      Y^N
    |_{
      \calC^{ \beta }( [0,T], \left\| \cdot \right\|_{  \mathscr{L}^p( \P; \left\| \cdot \right\|_E )  } )
    }
    < \infty
  \big]
$
or 
$
\big[
  \sup_{ N \in \N } 
  \nicefrac{ d_{ \max }( \theta^N ) }{ d_{ \min }( \theta^N ) }
  < \infty
$
and 
$
  \forall \, N \in \N \colon
  Y^N 
  =
  [ Y^N ]_{
    \theta^N
  } 
\big]
\big)
$.
Then it holds for all
$ \alpha \in [ 0, \beta - \nicefrac{ 1 }{ p } ) $,
$ \eps \in (0,\infty) $
that
\begin{equation}
\label{eq:convHoelder1}
  \sup_{ N \in \N }
  \left[
  \E \Bigl[
      \| Y^N \|_{ \calC^{ \alpha }( [0,T], \left\| \cdot \right\|_E )  }^p
  \Bigr]
  +
  \bigl|
    d_{ \max }( \theta^N )
  \bigr|^{ - ( \beta - \alpha - \nicefrac{1}{p} - \eps ) }
  \Bigl( \E \Bigl[
      \| Y^0 - Y^N \|_{ \calC^{ \alpha }( [0,T], \left\| \cdot \right\|_E )  }^p
  \Bigr] \Bigr)^{ \nicefrac{1}{p} }
  \right]
  < \infty .
\end{equation}
\end{corollary}

\begin{proof}[Proof of Corollary~\ref{cor:hoelder3}]
Throughout this proof let $ c_0 \in [0,\infty) $,
$ c_1, c_2 \in [0,\infty] $ be the extended real numbers given by
\begin{equation}
\begin{split}
  c_0 
& =
  | Y^0 |_{
    \calC^{ \beta }( [0,T] , \left\| \cdot \right\|_{  \mathscr{L}^p( \P; \left\| \cdot \right\|_E )  } )
  }
  +
  \sup_{ N \in \N }
  \left[
  \left| 
    d_{ \max }( \theta^N )
  \right|^{ - \beta }
  \sup\nolimits_{ 
    t \in \theta^N
  }
  \|
    Y_t^0 -
    Y_t^N
  \|_{ \mathscr{L}^p( \P; \left\| \cdot \right\|_E ) }
  \right]
  ,
\\ 
  c_1
& =
  \sup_{ N \in \N } 
  \left[
    \frac{ d_{ \max }( \theta^N ) }{ d_{ \min }( \theta^N ) }
  \right]
  ,
\qquad 
  \text{and}
\qquad
  c_2
  =
  \sup_{ N \in \N }
  | 
    Y^N
  |_{
    \calC^{ \beta }( [0,T], \left\| \cdot \right\|_{  \mathscr{L}^p( \P; \left\| \cdot \right\|_E )  } )
  }
  \,
  .
\end{split}
\end{equation}
Next we observe that
Item~\eqref{item:II_cor_hoelder2} in Corollary~\ref{cor:hoelder1}
ensures 
for all 
$ r \in ( \nicefrac{ 1 }{ p } , \beta ] $,
$ \alpha \in [ 0, r - \nicefrac{ 1 }{ p } ) $,
$ N \in \N $
that
\begin{align}
\nonumber &
  \Bigl( \E \Bigl[
      \| Y^0 - Y^N \|_{ \calC^{ \alpha }( [0,T], \left\| \cdot \right\|_E )  }^p
  \Bigr] \Bigr)^{ \nicefrac{1}{p} }
\leq
  \Xi_{ T, p, \alpha, r }
 \left(
      2 
      \left| 
        d_{ \max }( \theta^N )
      \right|^{ - r }
    + 
    1
 \right)
  \bigg[
    \sup_{ t \in \theta^N }
    \|
      Y^0_t
      -
      Y^N_t
    \|_{ \mathscr{L}^p( \P; \left\| \cdot \right\|_E ) }
\\ \nonumber & \quad
    +
      \left| 
        d_{ \max }( \theta^N ) 
      \right|^{ \beta }
    \big(
    | 
      Y^0 
    |_{
      \calC^{ \beta }( [0,T], \left\| \cdot \right\|_{  \mathscr{L}^p( \P ; \left\| \cdot \right\|_E )  } )
    }
    +
    | 
      Y^N 
    |_{
      \calC^{ \beta }( [0,T], \left\| \cdot \right\|_{  \mathscr{L}^p( \P ; \left\| \cdot \right\|_E )  } )
    }
    \big)
  \bigg]
\\ & \leq
  \Xi_{ T, p, \alpha, r }
 \left(
      2 
      \left| 
        d_{ \max }( \theta^N )
      \right|^{ ( \beta - r ) }
    + 
      \left| 
        d_{ \max }( \theta^N )
      \right|^{ \beta }
 \right)
  \Big[
    c_0
    +
    | 
      Y^N 
    |_{
      \calC^{ \beta }( [0,T], \left\| \cdot \right\|_{  \mathscr{L}^p( \P ; \left\| \cdot \right\|_E )  } )
    }
  \Big]
\\ \nonumber & \leq
  \Xi_{ T, p, \alpha, r }
 \left(
      2 
      \left| 
        d_{ \max }( \theta^N )
      \right|^{ ( \beta - r ) }
    + 
      \left| 
        d_{ \max }( \theta^N )
      \right|^{ \beta }
 \right)
  \left[
    c_0
    +
    c_2
  \right]
\\ \nonumber & =
  \Xi_{ T, p, \alpha, r }
 \left(
      2 
    + 
      \left| 
        d_{ \max }( \theta^N )
      \right|^{ r }
 \right)
      \left| 
        d_{ \max }( \theta^N )
      \right|^{ ( \beta - r ) }
  \left[
    c_0
    +
    c_2
  \right]
  .
\end{align}
This implies 
for all 
$ r \in ( \nicefrac{ 1 }{ p } , \beta ] $,
$ \alpha \in [ 0, r - \nicefrac{ 1 }{ p } ) $
that
\begin{equation}
\begin{split}
  \sup_{ N \in \N }
  \left[
      \bigl| 
        d_{ \max }( \theta^N )
      \bigr|^{ - ( \beta - r ) }
  \Bigl( \E \Bigl[
      \| Y^0 - Y^N \|_{ \calC^{ \alpha }( [0,T], \left\| \cdot \right\|_E )  }^p
  \Bigr] \Bigr)^{ \nicefrac{1}{p} }
  \right]
\leq
  \Xi_{ T, p, \alpha, r }
 \left(
      2 
    + 
      T^r
 \right)
  \left[
    c_0
    +
    c_2
  \right]
  .
\end{split}
\end{equation}
Hence, we obtain 
for all 
$ \alpha \in [ 0, \beta - \nicefrac{ 1 }{ p } ) $,
$ r \in ( \alpha + \nicefrac{ 1 }{ p } , \beta ] $
that
\begin{equation}
\begin{split}
&
  \sup_{ N \in \N }
  \left[
      \bigl| 
        d_{ \max }( \theta^N )
      \bigr|^{ - ( \beta - \alpha - \nicefrac{1}{p} - [ r - \alpha - \nicefrac{1}{p} ] ) }
  \Bigl( \E \Bigl[
      \| Y^0 - Y^N \|_{ \calC^{ \alpha }( [0,T], \left\| \cdot \right\|_E )  }^p
  \Bigr] \Bigr)^{ \nicefrac{1}{p} }
  \right]
\\ & \leq
  \Xi_{ T, p, \alpha, \alpha + \nicefrac{1}{p} + [ r - \alpha - \nicefrac{1}{p} ] }
 \left(
      3 
    + 
      T
 \right)
  \left(
    c_0
    +
    c_2
  \right)
  .
\end{split}
\end{equation}
This shows 
for all 
$ \alpha \in [ 0, \beta - \nicefrac{ 1 }{ p } ) $,
$ \varepsilon \in ( 0 , \beta - \alpha - \nicefrac{ 1 }{ p } ] $
that
\begin{equation}
\label{eq:RATE_first_limit}
\begin{split}
&  \sup_{ N \in \N }
  \left[
      \bigl| 
        d_{ \max }( \theta^N )
      \bigr|^{ - ( \beta - \alpha - \nicefrac{1}{p} - \varepsilon ) }
  \Bigl( \E \Bigl[
      \| Y^0 - Y^N \|_{ \calC^{ \alpha }( [0,T], \left\| \cdot \right\|_E )  }^p
  \Bigr] \Bigr)^{ \nicefrac{1}{p} }
  \right]
\\ & \leq
  \Xi_{ T, p, \alpha, \alpha + \nicefrac{1}{p} + \varepsilon }
 \left(
      3 
    + 
      T
 \right)
  \left(
    c_0
    +
    c_2
  \right)
  .
\end{split}
\end{equation}
In the next step we note that
Item~\eqref{item:2_cor_hoelder2} in Corollary~\ref{cor:hoelder2}
proves 
for all 
$ r \in ( \nicefrac{ 1 }{ p }, \beta ] $,
$ \alpha \in [ 0, r - \nicefrac{ 1 }{ p } ) $,
$ N \in \N $
that
\begin{align}
& \nonumber
\Bigl( \E \Bigl[
    \bigl\| Y^0 - [ Y^N ]_{ \theta^N } \bigr\|_{ \calC^{ \alpha }( [0,T], \left\| \cdot \right\|_E )  }^p
\Bigr] \Bigr)^{ \nicefrac{1}{p} }
\leq
  \Xi_{ T, p, \alpha, r }
  \biggl(
  \Big[
  \tfrac{ 
    2
    \,
    | 
      d_{ \max }( \theta^N ) 
    |^{ 1 - r } 
  }{ 
    d_{ \min }( \theta^N ) 
  }
  + 1
  \Big]
    \sup_{ t \in \theta^N }
    \|
      Y^0_t
      -
      Y^N_t
    \|_{ \mathscr{L}^p( \P; \left\| \cdot \right\|_E ) }
\\ & \quad
+
  \left[
    2 \, | d_{ \max }( \theta^N ) |^{ - r }
    + 2^{ - \beta }
  \right]
  \left| 
    d_{ \max }( \theta^N ) 
  \right|^{ \beta }
  | Y^0 |_{
    \calC^{ \beta }( [0,T], \left\| \cdot \right\|_{  \mathscr{L}^p( \P; \left\| \cdot \right\|_E )  } )
  }
  \biggr)
  .
\end{align}
This implies 
for all 
$ r \in ( \nicefrac{ 1 }{ p }, \beta ] $,
$ \alpha \in [ 0, r - \nicefrac{ 1 }{ p } ) $,
$ N \in \N $
that
\begin{equation}
\begin{split}
&
\Bigl( \E \Bigl[
    \bigl\| Y^0 - [ Y^N ]_{ \theta^N } \bigr\|_{ \calC^{ \alpha }( [0,T], \left\| \cdot \right\|_E )  }^p
\Bigr] \Bigr)^{ \nicefrac{1}{p} }
\\ & \leq
  c_0 
  \left|
    d_{ \max }( \theta^N )
  \right|^{ \beta }
  \Xi_{ T, p, \alpha, r }
  \Big(
  \tfrac{ 
    2
    \,
    | 
      d_{ \max }( \theta^N ) 
    |^{ 1 - r } 
  }{ 
    d_{ \min }( \theta^N ) 
  }
  + 1
  +
    2 \, | d_{ \max }( \theta^N ) |^{ - r }
    + 2^{ - \beta }
  \Big)
\\ & \leq
  2 \, c_0 
  \left|
    d_{ \max }( \theta^N )
  \right|^{ \beta }
  \Xi_{ T, p, \alpha, r }
  \left(
    \left[ c_1 + 1 \right]
    \left| 
      d_{ \max }( \theta^N ) 
    \right|^{ - r } 
    + 1
  \right)
  .
\end{split}
\end{equation}
Hence, we obtain 
for all 
$ r \in ( \nicefrac{ 1 }{ p }, \beta ] $,
$ \alpha \in [ 0, r - \nicefrac{ 1 }{ p } ) $
that
\begin{equation}
\begin{split}
&
  \sup_{ N \in \N }
  \left[
  \bigl|
    d_{ \max }( \theta^N )
  \bigr|^{ - ( \beta - r ) }
  \Bigl( \E \Bigl[
      \bigl\| Y^0 - [ Y^N ]_{ \theta^N } \bigr\|_{ \calC^{ \alpha }( [0,T], \left\| \cdot \right\|_E )  }^p
  \Bigr] \Bigr)^{ \nicefrac{1}{p} }
  \right]
\\ & \leq
  2 \, c_0 \,
  \Xi_{ T, p, \alpha, r }
  \left(
    c_1 + 1 
    + T^r
  \right)
\leq
  2 \, c_0 \,
  \Xi_{ T, p, \alpha, r }
  \left(
    2 + T + c_1  
  \right)
  .
\end{split}
\end{equation}
This shows
for all 
$ \alpha \in [ 0, \beta - \nicefrac{ 1 }{ p } ) $,
$ r \in ( \alpha + \nicefrac{ 1 }{ p }, \beta ] $
that
\begin{equation}
\begin{split}
&
  \sup_{ N \in \N }
  \left[
  \bigl|
    d_{ \max }( \theta^N )
  \bigr|^{ - ( \beta - \alpha - \nicefrac{1}{p} - [ r - \alpha - \nicefrac{1}{p} ] ) }
  \Bigl( \E \Bigl[
      \bigl\| Y^0 - [ Y^N ]_{ \theta^N } \bigr\|_{ \calC^{ \alpha }( [0,T], \left\| \cdot \right\|_E )  }^p
  \Bigr] \Bigr)^{ \nicefrac{1}{p} }
  \right]
\\ & \leq
  2 \, c_0 \,
  \Xi_{ T, p, \alpha, \alpha + \nicefrac{1}{p} + [ r - \alpha - \nicefrac{1}{p} ] }
  \left(
    2 + T + c_1  
  \right)
  .
\end{split}
\end{equation}
This establishes
for all 
$ \alpha \in [ 0, \beta - \nicefrac{ 1 }{ p } ) $,
$ \varepsilon \in ( 0 , \beta - \alpha - \nicefrac{ 1 }{ p } ] $
that
\begin{equation}
\label{eq:RATE_second_limit}
\begin{split}
&
  \sup_{ N \in \N }
  \left[
  \bigl|
    d_{ \max }( \theta^N )
  \bigr|^{ - ( \beta - \alpha - \nicefrac{1}{p} - \varepsilon ) }
  \Bigl( \E \Bigl[
      \bigl\| Y^0 - [ Y^N ]_{ \theta^N } \bigr\|_{ \calC^{ \alpha }( [0,T], \left\| \cdot \right\|_E )  }^p
  \Bigr] \Bigr)^{ \nicefrac{1}{p} }
  \right]
\\ & \leq
  2 \, c_0 \,
  \Xi_{ T, p, \alpha, \alpha + \nicefrac{1}{p} + \varepsilon }
  \left(
    2 + T + c_1  
  \right)
  .
\end{split}
\end{equation}
Combining \eqref{eq:RATE_first_limit} and \eqref{eq:RATE_second_limit}
assures for all
$ \alpha \in [ 0, \beta - \nicefrac{ 1 }{ p } ) $,
$ \eps \in (0,\infty) $
that
\begin{equation}
\label{eq:RATE_both_limits}
  \sup_{ N \in \N } 
  \left[
  \bigl|
    d_{ \max }( \theta^N )
  \bigr|^{ - ( \beta - \alpha - \nicefrac{1}{p} - \eps ) }
  \Bigl( \E \Bigl[
      \| Y^0 - Y^N \|_{ \calC^{ \alpha }( [0,T], \left\| \cdot \right\|_E )  }^p
  \Bigr] \Bigr)^{ \nicefrac{1}{p} }
  \right]
  < \infty .
\end{equation}
In addition, note that the assumption that
$ Y_0^0 \in \mathscr{L}^p( \P; \left\| \cdot \right\|_E ) $,
the assumption that
$
  | Y^0 |_{
    \calC^{ \beta }( [0,T] , \left\| \cdot \right\|_{  \mathscr{L}^p( \P; \left\| \cdot \right\|_E )  } )
  }
  < \infty
$,
the assumption that 
$ Y^0 $ has continuous sample paths, 
and Theorem~\ref{thm:Kolmogorov}
ensure for all
$ \alpha \in [ 0, \beta - \nicefrac{ 1 }{ p } ) $
that
$
\E \bigl[
    \| Y^0 \|_{ \calC^{ \alpha }( [0,T], \left\| \cdot \right\|_E )  }^p
\bigr]
  < \infty
$.
This and \eqref{eq:RATE_both_limits}
complete the proof of Corollary~\ref{cor:hoelder3}.
\end{proof}

The next result, Corollary~\ref{cor:eulermethod} below, 
illustrates
Corollary~\ref{cor:hoelder3}
through 
a simple example.
For this note that
standard results for the Euler-Maruyama method
show
under suitable hypotheses
for every 
$ p \in [2,\infty) $,
$ \beta \in [0,\nicefrac{ 1 }{ 2 } ] $
that
condition~\eqref{eq:ass_convergencegrid1}
in Corollary~\ref{cor:hoelder3}
with uniform time steps
is satisfied
(cf., e.g., Section~10.6 in~Kloeden \& Platen~\cite{kp92}).
The convergence rate established in Corollary~\ref{cor:eulermethod} (see~\eqref{eq:convSup5} below)
is essentially sharp; see Proposition~\ref{prop:brownian} below.
Corollary~\ref{cor:eulermethod} is related to
Theorem~1.2 in~\cite{CoxVanNeerven2010} 
and 
Theorem~1.1 in~\cite{CoxVanNeerven2013}.

\begin{corollary}[Euler-Maruyama 
method]
\label{cor:eulermethod}
Consider the notation in Subsection~\ref{notation},
let $ T \in (0,\infty) $, $ d, m \in \N $,
let
$
  \left( \Omega, \mathscr{F},
  \P \right) 
$
be a probability space
with a normal filtration
$
  \left( \mathscr{F}_t \right)_{
    t \in [0,T]
  }
$,
let
$
  W \colon [0,T] \times
  \Omega \rightarrow \mathbb{R}^m
$
be
an $m$-dimensional standard $ ( \mathscr{F}_t )_{ t \in [0,T] } $-Brownian motion
with continuous sample paths,
let
$
  \mu \colon \mathbb{R}^d
  \rightarrow \mathbb{R}^d
$
and
$
  \sigma \colon \mathbb{R}^d
  \rightarrow 
  \mathbb{R}^{ d \times m }
$
be globally Lipschitz continuous
functions,
let
$ 
  X \colon [0,T] \times
  \Omega \rightarrow \mathbb{R}^d
$
be an 
$
( \mathscr{F}_t )_{ t \in [0,T] }
$/$ \mathscr{B}( \mathbb{R}^d ) $-adapted
stochastic process
with continuous sample paths
which satisfies
$ 
  \forall \, p \in [1,\infty) 
  \colon
  \E\!\left[ 
    \left\| X_0 
    \right\|^p_{ \mathbb{R}^d }
  \right]
  < \infty
$
and which satisfies
for all $ t \in [0,T] $
that
\begin{equation}
[ X_t ]_{ \P, \mathscr{B}( \R^d ) } = \biggl[ X_0+\int_0^t\mu( X_s ) \ds \biggr]_{ \P, \mathscr{B}( \R^d ) }
            +\int_0^t \sigma( X_s ) \dWs,
\end{equation}
and let
$ 
  Y^N \colon [0,T] \times
  \Omega \rightarrow \mathbb{R}^d 
$,
$ N \in \N $,
be mappings which satisfy
for all 
$ N \in \N $,
$ n \in \{ 0, 1, \dots, N-1 \} 
$,
$ 
  t \in 
  [ 
    \frac{ n T }{ N }, 
    \frac{ (n+1) T }{ N }
  ] 
$
that
$ Y^N_0 = X_0 $ 
and
\begin{equation}
  Y_t^N
  = 
  Y_{ \frac{nT}{N} }^N
  +
  \left(
    t - \tfrac{ n T }{ N }
  \right) \cdot
  \mu\big(
    Y_{
      \frac{nT}{N}
    }^N
  \big)
  +
  \left(
    \tfrac{ t N }{ T } - n
  \right) \cdot
  \sigma\big(
    Y_{
      \frac{nT}{N}
    }^N
  \big)
  \big(
    W_{ 
      \frac{(n + 1) T}{N}
    }
    -
    W_{ 
      \frac{nT}{N}
    }
  \big).
\end{equation}
Then it holds
for all 
$ \alpha \in [ 0, \nicefrac{ 1 }{ 2 } ) $,
$ \eps \in (0,\infty) $,
$ p \in [1,\infty) $
that
\begin{equation}  
\begin{split}
\label{eq:convSup5} 
  \sup_{ N \in \N }
  \left[
  N^{ 
    \nicefrac{ 1 }{ 2 } - \alpha - \eps
  }
  \Bigl( \E \Bigl[
      \| X - Y^N \|_{ \calC^{ \alpha }( [0,T], \left\| \cdot \right\|_{\mathbb{R}^d} )  }^p
  \Bigr] \Bigr)^{ \nicefrac{1}{p} }
  \right]
  < \infty
  .
\end{split}     \end{equation}
\end{corollary}

\subsection{Lower error bounds for stochastic processes
with H\"{o}lder continuous
sample paths}
\label{sec:holder3}

In this subsection we comment on
the optimality of the convergence rate provided by 
Corollary~\ref{cor:hoelder3}
and Corollary~\ref{cor:eulermethod}, respectively.
In particular, in the setting of 
Corollary~\ref{cor:eulermethod},
Theorem~3 in
M\"{u}ller-Gronbach~\cite{m02} shows in the case $ \alpha = 0 $ 
that there exists a class of SDEs for which 
the factors
$ N^{ \nicefrac{1}{2} - \eps } $, $ N \in \N $,
on the left hand side of the estimate~\eqref{eq:convSup5}
can at best -- up to a constant -- 
be replaced by the factors 
$
  \frac{ N^{ \nicefrac{1}{2} } }{ \log( N ) } 
$,
$ N \in \N $.
In Proposition~\ref{prop:brownian} below 
we show for every $ \alpha \in [ 0, \nicefrac{ 1 }{ 2 } ) $ 
in the simple case of $ \mu = 0 $ 
and $ \sigma = ( \mathbb{R} \ni x \mapsto 1 \in \mathbb{R} ) $ 
in Corollary~\ref{cor:eulermethod} that
the factors
$ N^{ \nicefrac{1}{2} - \alpha - \eps } $, $ N \in \N $,
on the left hand side of the estimate~\eqref{eq:convSup5}
can at best -- up to a constant -- 
be replaced by the factors 
$
  N^{ \nicefrac{1}{2} - \alpha } 
$,
$ N \in \N $.
\iftoggle{arXiv:v3}{%
Our proof of Proposition~\ref{prop:brownian}
uses the following elementary lemma.
\begin{lemma}
\label{lem:holderprocesses}
Consider the notation in Subsection~\ref{notation},
let $ T \in (0,\infty) $,
$ p \in [1,\infty) $,
$ \alpha \in [0,1] $,
let 
$
  \left( \Omega, \mathscr{F},
  \mathbb{P} \right)
$
be a probability space,
let
$
  \left(
    E, \left\| \cdot \right\|_E
  \right)
$
be a normed vector space,
and let
$ 
  X \colon [0,T] \times
  \Omega \rightarrow E 
$
be an $ ( \mathscr{F}, \left\| \cdot \right\|_E ) $-strongly measurable stochastic process
with continuous sample paths.
Then
\begin{equation}
\max\!\big\{
  | X |_{
    \calC^{ \alpha }( [0,T], 
      \left\| \cdot \right\|_{\mathscr{L}^p( \P;  \left\| \cdot \right\|_E )}
    )
  }
  ,
  2^{
    ( \nicefrac{1}{p} - 1 )
  }
  \,
  \| X \|_{
    \calC^{ \alpha }( [0,T], 
      \left\| \cdot \right\|_{\mathscr{L}^p( \P;  \left\| \cdot \right\|_E )}
    )
  }
\big\}
\leq
\Bigl( \E \Bigl[
    \| X \|_{ \calC^{ \alpha }( [0,T], \left\| \cdot \right\|_E )  }^p
\Bigr] \Bigr)^{ \nicefrac{1}{p} }
.
\end{equation}
\end{lemma}

The proof of 
Lemma~\ref{lem:holderprocesses}
is clear.
Instead we now present the promised proposition
on the optimality of the convergence rate estimate 
in Corollary~\ref{cor:eulermethod}.
}{%
Proposition~\ref{prop:brownian} is proved
in Subsection~2.3 in~\cite{CoxHutzenthalerJentzenVanNeervenWelti2016arXiv}
(see \cite[Proposition~2.14]{CoxHutzenthalerJentzenVanNeervenWelti2016arXiv}).%
}

\begin{proposition}
\label{prop:brownian}
Consider the notation in Subsection~\ref{notation},
let $ T \in (0,\infty) $,
let 
$
  \left( \Omega, \mathscr{F},
  \mathbb{P} \right)
$
be a probability space,
let
$ 
  W \colon [0,T] \times
  \Omega \rightarrow \mathbb{R} 
$
be
a one-dimensional standard Brownian motion
with continuous sample paths,
and let
$ 
  W^N \colon [0,T] \times
  \Omega \rightarrow \mathbb{R} 
$,
$ N \in \mathbb{N} $,
be mappings which satisfy
for all 
$ N \in \mathbb{N} $,
$ n \in \{ 0, 1, \dots, N-1 \} $,
$ 
  t \in 
  \big[ 
    \frac{ n T }{ N }, 
    \frac{ (n+1) T }{ N }
  \big] 
$
that
\begin{equation}
\label{eq:def_WN_t}
  W^N_t 
  = 
  \left(
    n + 1 - \tfrac{ t N }{ T }
  \right) \cdot
  W_{ \frac{ n T }{ N } }
  +
  \left(
    \tfrac{ t N }{ T } - n
  \right) \cdot
  W_{ \frac{ (n+1) T }{ N } }
  .
\end{equation}
Then it holds 
for all $ \alpha \in [0, \nicefrac{ 1 }{ 2 } ] $,
$ p \in [1,\infty) $,
$ N \in \{ 2, 3, 4, \ldots \} $
that
\begin{equation}
\label{eq:lower_bound_I}
  \| 
    W - W^N
  \|_{
    C( 
      [0,T], 
      \left\| \cdot \right\|_{  \mathscr{L}^p( \P;  \left| \cdot \right| )  }
    )
  }
=
  \tfrac{
    \left\| W_T \right\|_{ 
      \mathscr{L}^p( \P;  \left| \cdot \right| ) 
    }
  }{
    2 \sqrt{ N }
  } ,
\end{equation}
\begin{equation}
\label{eq:lower_bound_II}
  \tfrac{
    | 
      W - W^N
    |_{
      \calC^{ \alpha }( 
        [0,T], 
        \left\| \cdot \right\|_{  \mathscr{L}^p( \P;  \left| \cdot \right| )  } 
      )
    }
  }{
      N^{ 
        \left( \alpha - \nicefrac{ 1 }{ 2 } \right)
      } 
      \,
      T^{ - \alpha }
      \,
      \left\| W_T \right\|_{
        \mathscr{L}^p( \P;  \left| \cdot \right| )
      }
  }
  =
 \tfrac{ 
    \left( \nicefrac{1}{2} - \alpha \right)^{ 
      \left( \nicefrac{ 1 }{ 2 } - \alpha \right)
    } 
  }{
    2^{ \alpha }
    \,
    \left( 
      1 - \alpha
    \right)^{ \left( 1 - \alpha \right) }
  } 
  \in
  \left[
    \tfrac{ 1 }{ \sqrt{2} }, 1
  \right] ,
\end{equation}
\begin{equation}
\label{eq:lower_bound_III}
  \tfrac{
    \| 
      W - W^N
    \|_{
      \calC^{ \alpha }( 
        [0,T], 
        \left\| \cdot \right\|_{  \mathscr{L}^p( \P;  \left| \cdot \right| )  } 
      )
    }
  }{
    N^{ 
      \left( \alpha - \nicefrac{ 1 }{ 2 } \right) 
    }
    \,
    T^{ - \alpha }
    \,
    \left\| W_T \right\|_{ 
      \mathscr{L}^p( \P;  \left| \cdot \right| ) 
    }
  }
  =
  \tfrac{
    T^{ \alpha }
  }{
    2 N^{ \alpha }
  } 
  +
 \tfrac{ 
    \left( \nicefrac{ 1 }{ 2 } - \alpha \right)^{ 
      \left( \nicefrac{ 1 }{ 2 } - \alpha \right)
    } 
  }{
    2^{ \alpha }
    \,
    \left( 
      1 - \alpha
    \right)^{ \left( 1 - \alpha \right) }
  } 
  \in
  \left[
    \tfrac{ 1 }{ \sqrt{2} }
    ,
    \tfrac{ 
      2 + T^{ \alpha }
    }{2} 
  \right]
  ,
\end{equation}
\begin{equation}
\label{eq:lower_bound_IIII}
  \tfrac{
    \left( \E \left[
        \| W - W^N \|_{ \calC^{ \alpha }( [0,T], \left| \cdot \right| )  }^p
    \right] \right)^{ \nicefrac{1}{p} }
  }{
    N^{ 
      \left( \alpha - \nicefrac{ 1 }{ 2 } \right) 
    }
    \,
    T^{ - \alpha }
    \,
    \left\| W_T \right\|_{ 
      \mathscr{L}^p( \P;  \left| \cdot \right| ) 
    }
  }
  \geq 
 \tfrac{ 
    \left( \nicefrac{1}{2} - \alpha \right)^{ 
      \left( \nicefrac{ 1 }{ 2 } - \alpha \right)
    } 
  }{
    2^{ \alpha }
    \,
    \left( 
      1 - \alpha
    \right)^{ \left( 1 - \alpha \right) }
  } 
  \geq
  \tfrac{ 1 }{ \sqrt{2} }
  .
\end{equation}
\end{proposition}
\iftoggle{arXiv:v3}{%
\begin{proof}[Proof
of Proposition~\ref{prop:brownian}]
Throughout this proof let 
$ f \colon [ 0, \nicefrac{ 1 }{ 2 } ] \to (0,\infty) $ 
be the function which satisfies
for all $ x \in [ 0, \nicefrac{ 1 }{ 2 } ] $ that
$
  f( x ) = 
  \frac{
    \left( 
      \nicefrac{ 1 }{ 2 } - x 
    \right)^{ ( \nicefrac{ 1 }{ 2 } - x ) }
  }{
    2^x \,
    \left( 
      1 - x 
    \right)^{ ( 1 - x ) }
  }
$
and let 
$ g_{ \alpha } \colon (0,1]^2 \to \R $,
$ \alpha \in [ 0, \nicefrac{ 1 }{ 2 } ] $,
be the functions which satisfy
for all 
$
  x, y \in (0,1]
$,
$ \alpha \in [ 0, \nicefrac{ 1 }{ 2 } ] $
that
\begin{equation}
  g_{ \alpha }( x, y )
  =
  \frac{
    x \left( 1 - x \right)
    +
    y \left( 1 - y \right)
  }{
    \left( x + y \right)^{ 2 \alpha }
  }
  .
\end{equation}
We first prove \eqref{eq:lower_bound_I}.
For this observe that it holds
for all 
$ N \in \mathbb{N} $,
$ n \in \{ 0, 1, \dots, N-1 \} $,
$ 
  t \in 
  \big[ 
    \frac{ n T }{ N }, 
    \frac{ (n+1) T }{ N }
  \big] 
$
that
\begin{equation}
\label{eq:dif_rep}
\begin{split}
  W_t - W^N_t
& =
  W_t -
  \left(
    n + 1 - \tfrac{ t N }{ T }
  \right) \cdot
  W_{ \frac{ n T }{ N } }
  -
  \left(
    \tfrac{ t N }{ T } - n
  \right) \cdot
  W_{ \frac{ (n+1) T }{ N } }
\\ & =
  \left(  
    n - \tfrac{ t N }{ T }
  \right)
  \cdot
  \left(
    W_{ \frac{ (n+1) T }{ N } } - W_t
  \right)
  +
  \left(
    n + 1 - \tfrac{ t N }{ T }
  \right)
  \cdot
  \left( 
    W_t -
    W_{ \frac{ n T }{ N } }
  \right).
\end{split}
\end{equation}
This and the fact that
\begin{equation}
 \forall \, N \in \N \colon \forall \, t \in \bigl( 0, \tfrac{T}{N} \bigr) \colon \forall \, p \in [1,\infty) \colon
    \Bigl\| \tfrac{W_t - W^N_t}{\| W_t - W^N_t \|_{ \mathscr{L}^2( \P;  \left| \cdot \right| ) } } \Bigr\|_{ \mathscr{L}^p( \P;  \left| \cdot \right| ) }
    = \tfrac{ \| W_T \|_{ \mathscr{L}^p( \P;  \left| \cdot \right| ) } }{\sqrt{ T }} 
\end{equation}
imply that it holds
for all $ N \in \N $, $ p \in [1,\infty) $
that
\begin{equation}
\label{eq:WN1}
\begin{split}
&
  \left\| 
    W - W^N
  \right\|_{
    C( 
      [0,T], 
      \left\| \cdot \right\|_{  \mathscr{L}^p( \P;  \left| \cdot \right| )  } 
    )
  }
=
  \sup\nolimits_{ t \in [0,T] }
  \left\|
    W_t - 
    W^N_t
  \right\|_{ \mathscr{L}^p( \P;  \left| \cdot \right| ) 
  }
  \\ &
=
  \sup\nolimits_{ t \in [0,\frac{T}{N}] }
  \left\|
    W_t - 
    W^N_t
  \right\|_{ \mathscr{L}^p( \P;  \left| \cdot \right| ) 
  }
= 
  \tfrac{ 
    \| W_T \|_{ \mathscr{L}^p( \P;  \left| \cdot \right| ) }
  }{
    \sqrt{ T }
  }
  \left[
  \sup\nolimits_{ t \in [0,\frac{T}{N}] }
  \left\|
    W_t - 
    W^N_t
  \right\|_{ \mathscr{L}^2( \P;  \left| \cdot \right| ) 
  }
  \right]
\\ & =
  \tfrac{ 
    \| W_T \|_{ \mathscr{L}^p( \P;  \left| \cdot \right| ) }
  }{
    \sqrt{ T }
  }
  \left[
  \sup\nolimits_{ t \in [0,\frac{T}{N}] }
  \left\|
    \tfrac{ t N }{ T }
    \cdot
    \left(
      W_t - W_{ \frac{ T }{ N } } 
    \right)
    +
    \left(
      1 - \tfrac{ t N }{ T }
    \right)
    \cdot
      W_t 
  \right\|_{ \mathscr{L}^2( \P;  \left| \cdot \right| ) 
  }
  \right]
\\ & =
  \tfrac{ 
    \| W_T \|_{ \mathscr{L}^p( \P;  \left| \cdot \right| ) }
  }{
    \sqrt{ T }
  }
  \left[
  \sup\nolimits_{ t \in [0,\frac{T}{N}] }
  \left[
    \left(
      \tfrac{ t N }{ T }
    \right)^2 \cdot  
    \left(
      \tfrac{ T }{ N } -
      t
    \right)
    +
    \left( 
      1 - \tfrac{ t N }{ T }
    \right)^2 \cdot
    t \,
  \right]^{ \nicefrac{ 1 }{ 2 } }
  \right]
\\ & =
  \tfrac{ 
    \| W_T \|_{ \mathscr{L}^p( \P;  \left| \cdot \right| ) }
  }{
    \sqrt{N}
  }
  \biggl[
    \sup_{ t \in [0,1] }
    \sqrt{ \left(
      t^2 \cdot  
      \left( 
        1 -
        t 
      \right)
      +
      \left( 
        1 - t
      \right)^2 \cdot
      t
    \right) }
  \biggr]
\\ & =
  \tfrac{ 
    \| W_T \|_{ \mathscr{L}^p( \P;  \left| \cdot \right| ) }
  }{
    \sqrt{N}
  }
  \biggl[
    \sup_{ t \in [0,1] }
    \sqrt{
      t \cdot
      \left( 
        1 - t
      \right)
    }
  \biggr]
  =
  \tfrac{ 
    \| W_T \|_{ \mathscr{L}^p( \P;  \left| \cdot \right| ) }
  }{
    2 \sqrt{N}
  }
  .
\end{split}
\end{equation}
This establishes \eqref{eq:lower_bound_I}.
In the next step we prove \eqref{eq:lower_bound_II}.
For this observe that
\eqref{eq:dif_rep}
shows for all
$ 
  N \in \{ 2, 3, 4, \ldots \}
$,
$ 
  n \in \{ 1, 2, \ldots, N - 1 \}
$,
$
  t_1 \in [0, \frac{T}{N}]
$,
$
  t_2 \in [ \frac{ n T }{ N },
  \frac{ (n+1) T }{ N } ]
$,
$ p \in [1,\infty) $
that
\begin{equation}
\label{eq:brownian_use1}
\begin{split}
&
    \left\|
      \left(
        W_{ t_2 } - W^N_{ t_2 }
      \right)
      -
      \left(
        W_{ t_1 } - W^N_{ t_1 }
      \right)
    \right\|_{ 
      \mathscr{L}^p( \P;  \left| \cdot \right| ) 
    }
\\ & =
  \Big\|
  \left(  
    n - \tfrac{ t_2 N }{ T }
  \right)
  \cdot
    \left(
      W_{ \frac{ (n+1) T }{ N } } - W_{ t_2 }
    \right)
    +
    \left(
      n + 1 - \tfrac{ t_2 N }{ T }
    \right)
    \cdot
    \left( 
      W_{ t_2 } -
      W_{ \frac{ n T }{ N } }
    \right)
\\ & \quad
    +
    \tfrac{ t_1 N }{ T }
    \cdot
    \left(
      W_{ \frac{ T }{ N } } - W_{ t_1 }
    \right)
    +
    \left(
      \tfrac{ t_1 N }{ T } - 1
    \right)
    \cdot
    W_{ t_1 }
  \Big\|_{ 
    \mathscr{L}^p( \P;  \left| \cdot \right| ) 
  }
\\ & =  
  \tfrac{
    \| W_T \|_{
      \mathscr{L}^p( \P;  \left| \cdot \right| )
    }
  }{
    \sqrt{ T }
  }
  \Bigl[
    \left(  
      n - \tfrac{ t_2 N }{ T }
    \right)^2
    \left(
      \tfrac{ (n+1) T }{ N } - t_2 
    \right)
    \\ & \hphantom{=\tfrac{\| W_T \|_{\mathscr{L}^p( \P;  \left| \cdot \right| )}}{\sqrt{ T }}\Bigl[}
    +
    \left(
      n + 1 - \tfrac{ t_2 N }{ T }
    \right)^2
    \left( 
      t_2 -
      \tfrac{ n T }{ N } 
    \right)    
    +
    \tfrac{ ( t_1 )^2 N^2 }{ T^2 }
    \left(
      \tfrac{ T }{ N } - t_1 
    \right)
    +
    \left( 
      \tfrac{ t_1 N }{ T } - 1
    \right)^2
    t_1
  \Bigr]^{ 
    \frac{ 1 }{ 2 } 
  }
\\ & =
  \tfrac{ 
    \| W_T \|_{
      \mathscr{L}^p( \P;  \left| \cdot \right| )
    }
  }{  
    \sqrt{N} 
  }
  \Big[ 
    \left(  
      \tfrac{ t_2 N }{ T } - n
    \right)
    \left(
      n + 1 - \tfrac{ t_2 N }{ T }
    \right)
    +
    \tfrac{ t_1 N }{ T }
    \left(
      1 - \tfrac{ t_1 N }{ T }
    \right)
  \Big]^{ 
    \frac{ 1 }{ 2 } 
  }
  .
\end{split}
\end{equation}
Moreover, 
\eqref{eq:def_WN_t}
ensures 
for all
$ N \in \mathbb{N} $,
$ 
  t_1, t_2 \in [0,\frac{T}{N} ] 
$,
$ p \in [1,\infty) $
with $ t_1 < t_2 $
that
\begin{align}
& \nonumber
    \left\|
      \left(
        W_{ t_2 } - W^N_{ t_2 }
      \right)
      -
      \left(
        W_{ t_1 } - W^N_{ t_1 }
      \right)
    \right\|_{ 
      \mathscr{L}^p( \P;  \left| \cdot \right| ) 
    }
\\ \nonumber & =
  \left\|
    \left(
      W_{ t_2 } - 
      \tfrac{ t_2 N }{ T } 
      \cdot
      W_{ \frac{ T }{ N } }
    \right)
    -
    \left(
      W_{ t_1 } 
      - 
      \tfrac{ t_1 N }{ T } \cdot
      W_{ \frac{ T }{ N } }
    \right)
  \right\|_{ 
    \mathscr{L}^p( \P;  \left| \cdot \right| ) 
  }
\\ \nonumber & =
  \tfrac{
    \| W_T \|_{
      \mathscr{L}^p( \P;  \left| \cdot \right| )
    }
  }{
    \sqrt{ T }
  }
  \left\|
    W_{ t_2 } - 
    W_{ t_1 } +
    \tfrac{ (t_1 - t_2) N }{ T } 
    \cdot
    W_{ \frac{ T }{ N } }
  \right\|_{ 
    \mathscr{L}^2( \P;  \left| \cdot \right| ) 
  }
\\ \label{eq:brownian_use2} & =
  \tfrac{
    \| W_T \|_{
      \mathscr{L}^p( \P;  \left| \cdot \right| )
    }
  }{
    \sqrt{ T }
  }
\\ \nonumber & \quad \cdot
  \left\|
    \left(
      1 +
      \tfrac{ (t_1 - t_2) N }{ T } 
    \right)
    \cdot
    \left(
      W_{ t_2 } - 
      W_{ t_1 } 
    \right)
    +
    \tfrac{ (t_1 - t_2) N }{ T } 
    \cdot
    \left(
      W_{ \frac{ T }{ N } } - W_{ t_2 }
    \right)
    +
    \tfrac{ (t_1 - t_2) N }{ T } 
    \cdot
    W_{ t_1 }
  \right\|_{ 
    \mathscr{L}^2( \P;  \left| \cdot \right| ) 
  }
\\ \nonumber & =
  \tfrac{
    \| W_T \|_{
      \mathscr{L}^p( \P;  \left| \cdot \right| )
    }
  }{
    \sqrt{ T }
  }
  \left[
    \left|
      1 +
      \tfrac{ (t_1 - t_2) N }{ T } 
    \right|^2
    \cdot
    \left( t_2 - t_1 \right)
    +
    \tfrac{ | t_1 - t_2 |^2 N^2 }{ T^2 } 
    \cdot
    \left(
      \tfrac{ T }{ N } - t_2 
    \right)
    +
    \tfrac{ |t_1 - t_2|^2 N^2 }{ T^2 } 
    \cdot
    t_1 
  \right]^{ \nicefrac{ 1 }{ 2 } }
\\ \nonumber & =
  \tfrac{
    \| W_T \|_{
      \mathscr{L}^p( \P;  \left| \cdot \right| )
    }
  }{
    \sqrt{ T }
  }
  \cdot
  \left[
    1 
    +
      \tfrac{ 2 (t_1 - t_2) N }{ T } 
    +
      \tfrac{ (t_1 - t_2)^2 N^2 }{ T^2 } 
    +
    \tfrac{ | t_1 - t_2 | \, N^2 }{ T^2 } 
    \cdot
    \left(
      \tfrac{ T }{ N } + t_1 - t_2 
    \right)
  \right]^{ \nicefrac{ 1 }{ 2 } }
  \cdot
  \left( t_2 - t_1 \right)^{ \nicefrac{ 1 }{ 2 } }
\\ \nonumber & =
  \tfrac{
    \| W_T \|_{
      \mathscr{L}^p( \P;  \left| \cdot \right| )
    }
  }{
    \sqrt{ T }
  }
  \cdot
  \left(
    1 + \tfrac{ (t_1 - t_2) N }{ T }
  \right)^{ \nicefrac{ 1 }{ 2 } }
  \cdot
  \left( t_2 - t_1 \right)^{ 
    \nicefrac{ 1 }{ 2 } 
  }
  .
\end{align}
Combining 
\eqref{eq:brownian_use1}
and
\eqref{eq:brownian_use2}
proves
for all $ N \in \{ 2, 3, 4, \ldots \} $,
$ \alpha \in [ 0, \nicefrac{ 1 }{ 2 } ] $,
$ p \in [1,\infty) $
that
\begin{equation}
\begin{split}
&
  \left| 
    W - W^N
  \right|_{
    \calC^{ \alpha }( 
      [0,T], 
      \left\| \cdot \right\|_{  \mathscr{L}^p( \P;  \left| \cdot \right| )  } 
    )
  }
  =
  \sup_{
    t_1, t_2 \in [0, T ] ,
    \,
    t_1 < t_2
  }
  \frac{
    \left\| 
      \left( W_{ t_2 } - W^N_{ t_2 } \right)
      -
      \left( W_{ t_1 } - W^N_{ t_1 } \right)
    \right\|_{
      \mathscr{L}^p( \P;  \left| \cdot \right| )
    }
  }{
    | t_1 - t_2 |^{ \alpha }
  }
\\ & =
  \sup_{
    t_1 \in [ 0, \frac{ T }{ N } ],
    \, t_2 \in [ 0, T ],
    \, t_1 < t_2
  }
  \frac{
    \left\| 
      \left( W_{ t_2 } - W^N_{ t_2 } \right)
      -
      \left( W_{ t_1 } - W^N_{ t_1 } \right)
    \right\|_{
      \mathscr{L}^p( \P;  \left| \cdot \right| )
    }
  }{
    | t_1 - t_2 |^{ \alpha }
  }
\\ & = 
  \tfrac{
    \| W_T \|_{ \mathscr{L}^p( \P;  \left| \cdot \right| ) }
  }{
    \sqrt{ T }
  }
  \left|
  \max\!\left\{
    \sup_{
      \substack{
        t_1, t_2 \in [0, \frac{T}{N} ] ,
      \\
        t_1 < t_2
      }
    }
    \tfrac{
      \left(
        1 + \frac{ (t_1 - t_2) N }{ T }
      \right)
    }{
      \left( t_2 - t_1 \right)^{ (2 \alpha - 1) }
    }
  ,
    \sup_{
      \substack{
        t_1 \in [0, \frac{T}{N} ] ,
      \\
        t_2 \in ( \frac{T}{N} , \frac{ 2 T }{ N } ]
      }
    }
    \tfrac{
      T \left[
      \left(  
        \tfrac{ t_2 N }{ T } - 1
      \right)
      \left(
        2 - \tfrac{ t_2 N }{ T }
      \right)
      +
      \tfrac{ t_1 N }{ T }
      \left(
        1 - \tfrac{ t_1 N }{ T }
      \right)
      \right]
    }{
      N \,
      \left( t_2 - t_1 \right)^{ 2 \alpha }
    }
  \right\}
  \right|^{ \frac{1}{2} }
  .
\end{split}
\end{equation}
This implies 
for all $ N \in \{ 2, 3, 4, \ldots \} $,
$ \alpha \in [ 0, \nicefrac{ 1 }{ 2 } ] $,
$ p \in [1,\infty) $
that
\begin{equation}
\begin{split}
&
  \left| 
    W - W^N
  \right|_{
    \calC^{ \alpha }( 
      [0,T], 
      \left\| \cdot \right\|_{  \mathscr{L}^p( \P;  \left| \cdot \right| )  } 
    )
  }
\\ & = 
  \tfrac{
    \| W_T \|_{ \mathscr{L}^p( \P;  \left| \cdot \right| ) }
  }{
    \sqrt{ T }
  }
  \left|
  \left|
    \tfrac{ T }{ N }
  \right|^{ 
    \left( 1 - 2 \alpha \right)
  }
  \max\!\left\{
    \sup_{
      \substack{
        x \in (0, 1 ] 
      }
    }
    \frac{ 
      \left(
        1 - x
      \right)
    }{
      x^{ (2 \alpha - 1) }
    }
  ,
    \sup_{
      \substack{
        x \in [0, 1 ] , \\
        y \in (1, 2]
      }
    }
    \frac{
      \left(  
        y - 1
      \right)
      \left(
        2 - y
      \right)
      +
      x
      \left(
        1 - x
      \right)
    }{
      \left( y - x \right)^{ 2 \alpha }
    }
  \right\}
  \right|^{ \frac{ 1 }{ 2 } }
\\ & = 
  \tfrac{
    \| W_T \|_{ \mathscr{L}^p( \P;  \left| \cdot \right| ) }
  }{
    \sqrt{ T }
  }
  \left|
    \tfrac{ T }{ N }
  \right|^{ 
    \left( \frac{1}{2} - \alpha \right)
  }
  \left[
    \max\!\left\{
      \sup_{
        \substack{
          x \in (0, 1 ] 
        }
      }
      \frac{
        x
        \left(
          1 - x
        \right)
      }{
        x^{ 2 \alpha }
      }
      ,
    \sup_{
      \substack{
        x \in [0, 1], \\
        y \in (0, 1]
      }
    }
    \frac{
      y
      \left(
        2 - ( y + 1 )
      \right)
      +
      x
      \left(
        1 - x
      \right)
    }{
      \left( [ y + 1 ] - [ 1 - x ] \right)^{ 2 \alpha }
    }
    \right\}
  \right]^{ \! \frac{1}{2} }
  .
\end{split}
\end{equation}
Hence, we obtain 
for all $ N \in \{ 2, 3, 4, \ldots \} $,
$ \alpha \in [ 0, \nicefrac{ 1 }{ 2 } ] $,
$ p \in [1,\infty) $
that
\begin{align}
\label{eq:WN2}
& \nonumber
  \left| 
    W - W^N
  \right|_{
    \calC^{ \alpha }( 
      [0,T], 
      \left\| \cdot \right\|_{  \mathscr{L}^p( \P;  \left| \cdot \right| )  }
    )
  }
\\ \nonumber & = 
  \tfrac{
    \| W_T \|_{ \mathscr{L}^p( \P;  \left| \cdot \right| ) }
  }{
    \sqrt{ T }
  }
  \left|
    \tfrac{ T }{ N }
  \right|^{ 
    \left( \frac{1}{2} - \alpha \right)
  }
  \left[
    \max\!\left\{
      \sup_{
        \substack{
          y \in (0, 1 ] 
        }
      }
      \frac{
        y
        \left(
          1 - y
        \right)
      }{
        y^{ 2 \alpha }
      }
      ,
    \sup_{
      x \in [0,1]
    }
    \sup_{
      y \in (0,1]
    }
    \frac{
      x
      \left(
        1 - x
      \right)
      +
      y
      \left(
        1 - y
      \right)
    }{
      \left( x + y \right)^{ 2 \alpha }
    }
    \right\}
  \right]^{ \! \frac{1}{2} }
\\ & =
  \tfrac{
    \| W_T \|_{ \mathscr{L}^p( \P;  \left| \cdot \right| ) }
  }{
    \sqrt{ T }
  }
  \left|
    \tfrac{ T }{ N }
  \right|^{ 
    \left( \frac{ 1 }{ 2 } - \alpha \right)
  }
  \left[
    \sup_{
      x \in [0,1]
    }
    \sup_{
      y \in (0,1]
    }
    \frac{
      x
      \left(
        1 - x
      \right)
      +
      y
      \left(
        1 - y
      \right)
    }{
      \left( x + y \right)^{ 2 \alpha }
    }
  \right]^{ \nicefrac{1}{2} }
\\ \nonumber & =
    N^{ ( \alpha - \nicefrac{ 1 }{ 2 } ) }
    \,
    T^{ - \alpha }
    \, 
    \| W_T \|_{ \mathscr{L}^p( \P;  \left| \cdot \right| ) } 
  \left[
    \sup_{
      x, y \in (0,1]
    }
    g_{ \alpha }( x, y )
  \right]^{ \nicefrac{ 1 }{ 2 } }
  .
\end{align}
To complete the proof of \eqref{eq:lower_bound_II},
we study a few properties of
the functions $ g_{ \alpha } $, $ \alpha \in [0, \nicefrac{ 1 }{ 2 } ] $.
Note that it holds for all $ x, y \in [ 0, 1 ] $ that
\begin{equation}
\label{eq:g_estimate}
 x(1-x)+y(1-y)
 =
 ( x + y ) \Bigl( 1- \frac{x+y}{2} \Bigr) - \frac{(x-y)^2}{2}
 \leq
 2 \, \Bigl( \frac{x+y}{2} \Bigr) \Bigl( 1 -\frac{x+y}{2} \Bigr).
\end{equation}
In addition, observe that it holds
for all $ \alpha \in [ 0, \nicefrac{ 1 }{ 2 } ] $, $ z \in ( 0, 1 ) $
that
\begin{equation}
  \tfrac{ \partial }{ \partial z }
  \bigl(
    z^{ ( 1 - 2 \alpha ) } \left( 1 - z \right)
  \bigr)
  =
  \left( 1 - 2 \alpha \right)
  z^{ - 2 \alpha } 
  \left( 1 - z \right)
  -
  z^{ ( 1 - 2 \alpha ) } 
  =
  -
  2 \left( 1 - \alpha \right)
  \bigl[
    z
    -
    \tfrac{
      \nicefrac{1}{2} - \alpha
    }{
      1 - \alpha 
    }
  \bigr]
  z^{ - 2 \alpha }
  .
\end{equation}
Combining this with \eqref{eq:g_estimate}
ensures
for all $ \alpha \in [ 0, \nicefrac{ 1 }{ 2 } ] $, $ x, y \in ( 0, 1 ] $
that
\begin{equation}
\begin{split}
g_{\alpha} ( x, y )
& =
\frac{ x ( 1 - x ) + y ( 1 - y ) }{ ( x + y )^{ 2 \alpha } }
\leq
2^{ ( 1 - 2 \alpha ) }
\Bigl( \frac{ x + y }{2} \Bigr)^{ 1 - 2\alpha } \Bigl( 1 - \frac{ x + y }{2} \Bigr)
\\ & \leq
2^{ ( 1 - 2 \alpha ) }
\sup_{ z \in ( 0, 1 ] } \bigl[ z^{ ( 1 - 2 \alpha ) } ( 1 - z ) \bigr]
=
2^{ ( 1 - 2 \alpha ) }
\bigl[
    \tfrac{ \nicefrac{1}{2} - \alpha }{ 1 - \alpha }
\bigr]^{ ( 1 - 2 \alpha ) }
\bigl[
    1 -
    \tfrac{ \nicefrac{1}{2} - \alpha }{ 1 - \alpha }
\bigr]
\\ & =
2^{ - 2 \alpha }
\bigl[
    \tfrac{ \nicefrac{1}{2} - \alpha }{ 1 - \alpha }
\bigr]^{ ( 1 - 2 \alpha ) }
\bigl[
    \tfrac{ 1 }{ 1 - \alpha }
\bigr]
=
\Bigl[
  \tfrac{
    ( \nicefrac{ 1 }{ 2 } - \alpha )^{
      ( \nicefrac{ 1 }{ 2 } - \alpha )
    }
  }{
    2^{ \alpha } \, ( 1 - \alpha )^{ ( 1 - \alpha ) }
  }
\Bigr]^2
=
\left[ 
    f( \alpha )
\right]^2
.
\end{split}
\end{equation}
This proves
for all $ \alpha \in [ 0, \nicefrac{ 1 }{ 2 } ] $
that
\begin{equation}
\begin{split}
\left[ 
    f( \alpha )
\right]^2
=
\sup_{ z \in ( 0, 1 ] }
\bigl[ ( 2 \, z )^{ ( 1 - 2 \alpha ) } ( 1 - z ) \bigr]
=
\sup_{ x \in ( 0, 1 ] }
g_{ \alpha }( x, x )
\leq
\sup_{ x, y \in (0,1] }
g_{ \alpha }( x, y )
\leq
\left[ 
    f( \alpha )
\right]^2
.
\end{split}
\end{equation}
This shows
for all $ \alpha \in [ 0, \nicefrac{ 1 }{ 2 } ] $
that
\begin{equation}
\label{eq:g_identity}
  \sup_{ x, y \in (0,1] }
  g_{ \alpha }( x, y )
  =
  \sup_{ x \in (0,1] }
  g_{ \alpha }( x, x )
  =
  \left[ 
    f( \alpha )
  \right]^2
  .
\end{equation}
Next note that it holds for all 
$ \alpha \in ( 0, \nicefrac{ 1 }{ 2 } ) $
that
\begin{equation}
\label{eq:f_identity}
  f( \alpha )
  =
  \exp\!\left(
    \left( \tfrac{ 1 }{ 2 } - \alpha \right) 
    \cdot
    \ln\!\left( \tfrac{ 1 }{ 2 } - \alpha \right)
    +
    \left( \alpha - 1 \right)
    \cdot 
    \ln\!\left( 1 - \alpha \right)
    -
    \alpha
    \cdot 
    \ln( 2 )
  \right)
  .
\end{equation}
Moreover, observe that it holds
for all $ \alpha \in (0, \nicefrac{ 1 }{ 2 } ) $
that
\begin{equation}
\begin{split}
&
  \tfrac{ \partial }{ \partial \alpha }
  \left(
    \left( \tfrac{ 1 }{ 2 } - \alpha \right) 
    \cdot
    \ln\!\left( \tfrac{ 1 }{ 2 } - \alpha \right)
    +
    \left( \alpha - 1 \right)
    \cdot 
    \ln\!\left( 1 - \alpha \right)
    -
    \alpha
    \cdot 
    \ln( 2 )
  \right)
\\ & =
    - \ln\!\left( \tfrac{ 1 }{ 2 } - \alpha \right)
    - 1
    +
    \ln\!\left( 1 - \alpha \right)
    + 1
    -
    \ln( 2 )
  =
    \ln\!\left( 1 - \alpha \right)
    - \ln\!\left( \tfrac{ 1 }{ 2 } - \alpha \right)
    -
    \ln( 2 )
\\ & 
  =
  \ln\!\left(
    \tfrac{ 1 - \alpha }{ 1 - 2 \alpha }
  \right)
  > 0
  .
\end{split}
\end{equation}
This and \eqref{eq:f_identity} ensure that $ f $ is strictly increasing.
Equation~\eqref{eq:g_identity} hence proves for all $ \alpha \in [ 0, \nicefrac{1}{2} ] $ that
\begin{equation}
  \sup_{ x, y \in (0,1] }
  g_{ \alpha }( x, y )
  =
  \sup_{ x \in (0,1] }
  g_{ \alpha }( x, x )
  =
  \left[ 
    f( \alpha )
  \right]^2
  \in
  \left[ 
    \left| f(0) \right|^2
    ,
    \left| f( \tfrac{ 1 }{ 2 } ) \right|^2
  \right]
  =
  \left[ 
    \tfrac{ 1 }{ 2 }
    ,
    1
  \right]
  .
\end{equation}
Putting this into \eqref{eq:WN2}
establishes \eqref{eq:lower_bound_II}.
Combining \eqref{eq:lower_bound_I}
with \eqref{eq:lower_bound_II}
proves \eqref{eq:lower_bound_III}.
Moreover, \eqref{eq:lower_bound_II}
and Lemma~\ref{lem:holderprocesses}
imply \eqref{eq:lower_bound_IIII}.
The proof of Proposition~\ref{prop:brownian}
is thus completed.
\end{proof}
}{}

\iftoggle{arXiv:v3}{%
\section{Basic results for mild solutions of SEEs}%
In this section we collect a number of elementary results for mild solution processes of
SEEs, most of which are well-known.%
}{}
\iftoggle{arXiv:v3}{%
\subsection{Temporal regularity of solutions of SEEs}
\begin{proposition}
\label{prop:temporal_regularity_SPDE}
Consider the notation in Subsection~\ref{notation},
let 
$ ( H, \left< \cdot, \cdot \right>_H, \left\| \cdot \right\|_H ) $
and
$ ( U, \left< \cdot, \cdot \right>_U, \left\| \cdot \right\|_U ) $
be separable $ \R $-Hilbert spaces,
let $ \mathbb{H} \subseteq H $ be a non-empty orthonormal basis of $ H $,
let $ \lambda \colon \mathbb{H} \to \mathbb{R} $
be a function with 
$
  \sup_{ h \in \mathbb{H} }
  \lambda_h < 0
$,
let
$
  A \colon D(A) \subseteq H \rightarrow H
$
be the linear operator which satisfies
$
  D(A) 
  = 
    \bigl\{ 
      v \in H 
      \colon
      \sum_{ h \in \mathbb{H} } 
	\left| 
	  \lambda_h 
	  \langle h, v \rangle_H 
	\right|^2
	< \infty
    \bigr\}
$
and
$
\forall \, v \in D(A) \colon
    Av
   =
    \sum_{ h \in \mathbb{H} } 
    \lambda_h \langle h, v \rangle_H h
$,
let 
$ 
  ( H_r , \left< \cdot , \cdot \right>_{ H_r }, \left\| \cdot \right\|_{ H_r } ) 
$,
$ r \in \R $,
be a family of interpolation spaces associated to
$ - A $
(cf., e.g., \cite[Section~3.7]{sy02}),
let 
$ T \in (0,\infty) $,
$ p \in [2,\infty) $,
$ \gamma \in \R $,
$ \eta \in [0, 1) $, 
$ \beta \in [ \gamma - \nicefrac{ \eta }{ 2 }, \gamma ] $,
$ F \in C( H_{ \gamma } , H_{ \gamma - \eta } ) $,
$ B \in C( H_{ \gamma } , \mathrm{HS}( U, H_{ \beta } ) ) $
satisfy
$
| F |_{
    \calC^1( H_{ \gamma }, \left\| \cdot \right\|_{  H_{ \gamma - \eta }  } )
}
+
| B |_{
    \calC^1( H_{ \gamma }, \left\| \cdot \right\|_{  \mathrm{HS}( U, H_{ \beta } )  } )
}
< \infty
$,
let 
$ ( \Omega, \calF, \P ) $
be a probability space with a normal filtration 
$ ( \calF_t )_{ t \in [0,T] } $,
let 
$
  ( W_t )_{ t \in [0,T] } 
$ 
be an $ \operatorname{Id}_U $-cylindrical
$ ( \Omega, \calF, \P, ( \calF_t )_{ t \in [0,T] } ) $-Wiener
process,
and let 
$ X \colon [0,T] \times \Omega \to H_{ \gamma } $
be an
$ ( \mathscr{F}_t )_{ t \in [0,T] } $/$ \mathscr{B}( H_{ \gamma } ) $-predictable stochastic process 
which satisfies
for all $ t \in [ 0, T ] $ that
$ 
\sup_{ s \in [0,T] } 
\| X_s \|_{ \mathscr{L}^p( \P ; \| \cdot \|_{ H_{ \gamma } } ) } 
< \infty 
$
and
\begin{equation}
\begin{split}
[ X_t ]_{ \P, \mathscr{B}( H_\gamma ) }
& = 
\biggl[
    e^{ t A } X_0
    +
    \int_0^t
        \mathbbm{1}_{
          \{ 
            \int_0^t
            \| e^{ ( t - u ) A } F( X_u ) \|_{ H_{ \gamma } } \dd u
            < \infty
          \} } \,
        e^{ ( t - s ) A }
        F( X_s )
    \ds
\biggr]_{ \P, \mathscr{B}( H_\gamma ) }
\\ & \quad +
\int_0^t
    e^{ ( t - s ) A }
    B( X_s )
\dWs
.
\end{split}
\end{equation}
Then it holds
for all 
$ 
r \in 
[
  \gamma ,
  \min\{ 
    1 + \gamma - \eta
    ,
    \nicefrac{ 1 }{ 2 } + \beta 
  \}
) 
$,
$
\varepsilon \in 
[
  0,
  \min\{ 
    1 + \gamma - \eta - r
    ,
    \nicefrac{ 1 }{ 2 } + \beta - r
  \}
)
$
that
$ \inf_{ s \in ( 0, T ] }
\P( X_s \in H_r ) = 1 $
and 
\begin{equation}
\begin{split}
& 
\sup_{ 
  \substack{ 
    t_1, t_2 \in [0,T] ,
  \\
    t_1 \neq t_2
  }
}
\Biggl(
\frac{ 
  \|
    (
      X_{ t_1 } - e^{ t_1 A } X_0
    ) \,
  \mathbbm{1}_{ \{ X_{ t_1 } \in H_r \} }
    -
    (
      X_{ t_2 } - e^{ t_2 A } X_0
    ) \,
  \mathbbm{1}_{ \{ X_{ t_2 } \in H_r \} }
  \|_{
    \mathscr{L}^p( 
      \P ; 
      \left\| \cdot \right\|_{ 
        H_{ r } 
      } 
    ) 
  }
}{
  \left| t_1 - t_2 \right|^{ \varepsilon }
}
\Biggr)
\\ & 
\leq 
\biggl[
  \sup_{ s \in [0,T] }
  \|
    F( X_s ) 
  \|_{
    \mathscr{L}^p( \P; \left\| \cdot \right\|_{ H_{ \gamma - \eta } } )
  }
\biggr]
  \frac{
    2
    \,
    T^{ 
      ( 1 + \gamma - \eta - r - \varepsilon ) 
    }
  }{
    ( 1 + \gamma - \eta - r - \varepsilon ) 
  }
\\ & \quad
+
\biggl[
  \sup_{ s \in [0,T] }
  \|
    B( X_s ) 
  \|_{
    \mathscr{L}^p( \P; \left\| \cdot \right\|_{ \mathrm{HS}( U, H_{ \beta } ) } )
  }
\biggr]
  \frac{ 
    \sqrt{ p \, ( p - 1 ) }
    \,
    T^{ ( \nicefrac{1}{2} + \beta - r - \varepsilon ) }
  }{
    ( 
      1 + 2 \beta - 2 r - 2 \varepsilon 
    )^{ \nicefrac{1}{2} }
  }
< \infty
.
\end{split}
\end{equation}
\end{proposition}
\begin{proof}[Proof of Proposition~\ref{prop:temporal_regularity_SPDE}]
Note that
the fact that it holds
for all $ u \in [ 0, 1 ] $ that
\begin{equation}
\label{eq:smoothing}
\sup_{ t \in ( 0, T ] }
    t^u \| ( - A )^u e^{ t A }
         \|_{ L( H ) }
\leq 1
\quad
\text{and}
\quad
\sup_{ t \in ( 0, T ] }
    t^{ -u } \| ( - A )^{ -u } ( e^{ tA } - \operatorname{Id}_H )
         \|_{ L( H ) }
\leq 1
\end{equation}
(cf., e.g., \cite[Lemma 12.36]{RenardyRogers2004})
implies that it holds
for all 
$ r \in [ \gamma, 1 + \gamma - \eta ) $,
$ \varepsilon \in [ 0, 1 + \gamma \linebreak - \eta - r ) $,
$ t_1 \in [0,T) $, $ t_2 \in ( t_1, T ] $
that
\begin{align*}
\label{eq:prop_temp_reg_F}
&
\biggl\|
  \int_0^{ t_1 }
    \mathbbm{1}_{
      \{ 
        \int_0^{ t_1 }
        \| e^{ ( t_1 - u ) A } F( X_u ) \|_{ H_{ r } } \dd u
        < \infty
      \} } \,
  e^{ ( t_1 - s ) A }
  F( X_s ) \ds
  \\
&\pushright{
  -
  \int_0^{ t_2 }
    \mathbbm{1}_{
      \{ 
        \int_0^{ t_2 }
        \| e^{ ( t_2 - u ) A } F( X_u ) \|_{ H_{ r } } \dd u
        < \infty
      \} } \,
  e^{ ( t_2 - s ) A }
  F( X_s ) \ds
\biggr\|_{
  \mathscr{L}^p( \P; \left\| \cdot \right\|_{ H_r } )
}
}
\\ & \leq 
\int_{ t_1 }^{ t_2 } 
\bigl\|
  e^{ ( t_2 - s ) A }
  F( X_s ) 
\bigr\|_{
  \mathscr{L}^p( \P; \left\| \cdot \right\|_{ H_r } )
}
\ds
+
\int_0^{ t_1 } 
\bigl\|
  e^{ ( t_1 - s ) A }
  \bigl(
    \operatorname{Id}_{ H_{ \gamma - \eta } }
    \mathop{-}
    e^{ ( t_2 - t_1 ) A }
  \bigr)
  F( X_s ) 
\bigr\|_{
  \mathscr{L}^p( \P; \left\| \cdot \right\|_{ H_r } )
}
\ds
\\ & \leq 
\biggl[
  \sup_{ s \in [0,T] }
  \|
    F( X_s ) 
  \|_{
    \mathscr{L}^p( \P; \left\| \cdot \right\|_{ H_{ \gamma - \eta } } )
  }
\biggr]
\biggl[ 
  \int_{ t_1 }^{ t_2 } 
  ( t_2 - s )^{ \gamma - \eta - r }
  \ds
  +
  \int_0^{ t_1 }
  ( t_1 - s )^{ \gamma - \eta - r - \varepsilon }
  ( t_2 - t_1 )^{ \varepsilon }
  \ds
\biggr]
\\ & = \yesnumber
\biggl[
  \sup_{ s \in [0,T] }
  \|
    F( X_s ) 
  \|_{
    \mathscr{L}^p( \P; \left\| \cdot \right\|_{ H_{ \gamma - \eta } } )
  }
\biggr]
\biggl[
  \frac{
    ( t_2 - t_1 )^{ 
      ( 1 + \gamma - \eta - r ) 
    }
  }{
    ( 1 + \gamma - \eta - r ) 
  }
  +
  \frac{ 
    ( t_2 - t_1 )^{ \varepsilon }
    ( t_1 )^{ ( 1 + \gamma - \eta - r - \varepsilon ) }
  }{
    ( 
      1 + \gamma - \eta - r - \varepsilon
    )
  }
\biggr] 
\\ & \leq 
\biggl[
  \sup_{ s \in [0,T] }
  \|
    F( X_s ) 
  \|_{
    \mathscr{L}^p( \P; \left\| \cdot \right\|_{ H_{ \gamma - \eta } } )
  }
\biggr]
\biggl[ 
  \frac{
    2 \,
    T^{ 
      ( 1 + \gamma - \eta - r - \varepsilon ) 
    }
    ( t_2 - t_1 )^{ \varepsilon }
  }{
    ( 1 + \gamma - \eta - r - \varepsilon ) 
  }
\biggr] 
.
\end{align*}
Moreover,
\eqref{eq:smoothing}
ensures
for all 
$ r \in [ \gamma, \nicefrac{ 1 }{ 2 } + \beta ) $,
$ \varepsilon \in [ 0, \nicefrac{ 1 }{ 2 } + \beta - r ) $,
$ t_1 \in [0,T) $, $ t_2 \in ( t_1, T ] $
that
\begin{equation}
\label{eq:prop_temp_reg_B}
\begin{split}
&
\biggl\|
  \int_0^{ t_1 } 
  e^{ ( t_1 - s ) A }
  B( X_s ) \dWs
  -
  \int_0^{ t_2 } 
  e^{ ( t_2 - s ) A }
  B( X_s ) \dWs
\biggr\|_{
  L^p( \P; \left\| \cdot \right\|_{ H_r } )
}^2
\\ & \leq 
\tfrac{ p \, ( p - 1 ) }{ 2 }
\int_{ t_1 }^{ t_2 } 
\bigl\|
  e^{ ( t_2 - s ) A }
  B( X_s ) 
\bigr\|_{
  \mathscr{L}^p( \P; \left\| \cdot \right\|_{ \mathrm{HS}( U, H_r ) } )
}^2
\ds
\\ & \quad
+
\tfrac{ p \, ( p - 1 ) }{ 2 }
\int_0^{ t_1 } 
\bigl\|
  e^{ ( t_1 - s ) A }
  \bigl(
    \operatorname{Id}_{ H_{ \beta } }
    \mathop{-}
    e^{ ( t_2 - t_1 ) A }
  \bigr)
  B( X_s ) 
\bigr\|_{
  \mathscr{L}^p( \P; \left\| \cdot \right\|_{ \mathrm{HS}( U, H_r ) } )
}^2
\ds
\\ & \leq 
\tfrac{ p \, ( p - 1 ) }{ 2 }
\biggl[
  \sup_{ s \in [0,T] }
  \|
    B( X_s ) 
  \|_{
    \mathscr{L}^p( \P; \left\| \cdot \right\|_{ \mathrm{HS}( U, H_{ \beta } ) } )
  }
\biggr]^2
\\ & \quad \cdot 
\biggl[ 
  \int_{ t_1 }^{ t_2 } 
  ( t_2 - s )^{ ( 2 \beta - 2 r ) }
  \ds
  +
  \int_0^{ t_1 }
  ( t_1 - s )^{ ( 2 \beta - 2 r - 2 \varepsilon ) }
  ( t_2 - t_1 )^{ 2 \varepsilon }
  \ds
\biggr]
\\ & \leq 
\tfrac{ p \, ( p - 1 ) }{ 2 }
\biggl[
  \sup_{ s \in [0,T] }
  \|
    B( X_s ) 
  \|_{
    \mathscr{L}^p( \P; \left\| \cdot \right\|_{ \mathrm{HS}( U, H_{ \beta } ) } )
  }
\biggr]^2
\biggl[
  \frac{
    2 \,
    T^{ ( 1 + 2 \beta - 2 r - 2 \varepsilon ) }
    ( t_2 - t_1 )^{ 2 \varepsilon }
  }{
    ( 
      1 + 2 \beta - 2 r - 2 \varepsilon 
    )
  }
\biggr]
.
\end{split}
\end{equation}
Combining \eqref{eq:prop_temp_reg_F}
and \eqref{eq:prop_temp_reg_B} completes the 
proof of Proposition~\ref{prop:temporal_regularity_SPDE}.
\end{proof}
\begin{corollary}
\label{cor:temporal_regularity_SPDE}
Consider the notation in Subsection~\ref{notation},
let 
$ ( H, \left< \cdot, \cdot \right>_H, \left\| \cdot \right\|_H ) $
and\linebreak
$ ( U, \left< \cdot, \cdot \right>_U, \left\| \cdot \right\|_U ) $
be separable $ \R $-Hilbert spaces,
let $ \mathbb{H} \subseteq H $ be a non-empty orthonormal basis of $ H $,
let $ \lambda \colon \mathbb{H} \to \mathbb{R} $
be a function with 
$
  \sup_{ h \in \mathbb{H} }
  \lambda_h < 0
$,
let
$
  A \colon D(A) \subseteq H \rightarrow H
$
be the linear operator which satisfies
$
  D(A) 
  = 
    \bigl\{ 
      v \in H 
      \colon
      \sum_{ h \in \mathbb{H} } 
	\left| 
	  \lambda_h 
	  \langle h, v \rangle_H 
	\right|^2
	< \infty
    \bigr\}
$
and
$
\forall \, v \in D(A) \colon
    Av
   =
    \sum_{ h \in \mathbb{H} } 
    \lambda_h \langle h, v \rangle_H h
$,
let 
$ 
  ( H_r , \left< \cdot , \cdot \right>_{ H_r }, \left\| \cdot \right\|_{ H_r } ) 
$,
$ r \in \R $,
be a family of interpolation spaces associated to
$ - A $,
let 
$ T \in (0,\infty) $,
$ p \in [2,\infty) $,
$ \gamma \in \R $,
$ \eta \in [0, 1) $, 
$ \beta \in [ \gamma - \nicefrac{ \eta }{ 2 }, \gamma ] $,
$ \delta \in [ \gamma, \infty) $,
$ F \in C( H_{ \gamma } , H_{ \gamma - \eta } ) $,
$ B \in C( H_{ \gamma } , \mathrm{HS}( U, H_{ \beta } ) ) $
satisfy
$
| F |_{
    \calC^1( H_{ \gamma }, \left\| \cdot \right\|_{  H_{ \gamma - \eta }  } )
}
+
| B |_{
    \calC^1( H_{ \gamma }, \left\| \cdot \right\|_{  \mathrm{HS}( U, H_{ \beta } )  } )
}
< \infty
$,
let 
$ ( \Omega, \calF, \P ) $
be a probability space with a normal filtration 
$ ( \calF_t )_{ t \in [0,T] } $,
let 
$
  ( W_t )_{ t \in [0,T] } 
$ 
be an $ \operatorname{Id}_U $-cylindrical
$ ( \Omega, \calF, \P, ( \calF_t )_{ t \in [0,T] } ) $-Wiener
process,
and let 
$ X \colon [0,T] \times \Omega \to H_{ \gamma } $
be an
$ ( \mathscr{F}_t )_{ t \in [0,T] } $/$ \mathscr{B}( H_{ \gamma } ) $-predictable stochastic process 
which satisfies
for all $ t \in [ 0, T ] $ that
$ 
\sup_{ s \in [0,T] } 
\| X_s \|_{ \mathscr{L}^p( \P ; \| \cdot \|_{ H_{ \gamma } } ) } 
< \infty 
$,
$ X_0( \Omega ) \subseteq H_\delta $,
$ \E [ \| X_0 \|_{ H_\delta }^p ] < \infty $,
and
\begin{equation}
\begin{split}
[ X_t ]_{ \P, \mathscr{B}( H_\gamma ) }
& = 
\biggl[
    e^{ t A } X_0
    +
    \int_0^t
        \mathbbm{1}_{
          \{ 
            \int_0^t
            \| e^{ ( t - u ) A } F( X_u ) \|_{ H_{ \gamma } } \dd u
            < \infty
          \} } \,
        e^{ ( t - s ) A }
        F( X_s )
    \ds
\biggr]_{ \P, \mathscr{B}( H_\gamma ) }
\\ & \quad +
\int_0^t
    e^{ ( t - s ) A }
    B( X_s )
\dWs
.
\end{split}
\end{equation}
Then it holds
for all 
$ 
r \in 
[
  \gamma ,
  \min\{ 
    1 + \gamma - \eta
    ,
    \nicefrac{ 1 }{ 2 } + \beta 
  \}
) 
$,
$
\varepsilon \in 
[
  0,
  \min\{ 
    1 + \gamma - \eta - r
    ,
    \nicefrac{ 1 }{ 2 } + \beta - r
  \}
)
$
that
$ \inf_{ s \in ( 0, T ] }
\P( X_s \in H_r ) = 1 $
and
\begin{equation}
\begin{split}
& 
\sup_{ 
  \substack{ 
    t_1, t_2 \in [0,T] ,
  \\
    t_1 \neq t_2
  }
}
\Biggl(
\frac{ 
  \left| \min\{ t_1, t_2 \} \right|^{
    \max\{ r + \varepsilon - \delta , 0 \}
  }
  \|
  \mathbbm{1}_{ \{ X_{ t_1 } \in H_r \} }
    X_{ t_1 } 
    -
  \mathbbm{1}_{ \{ X_{ t_2 } \in H_r \} }
    X_{ t_2 } 
  \|_{
    \mathscr{L}^p( 
      \P ; 
      \left\| \cdot \right\|_{ 
        H_{ r } 
      } 
    ) 
  }
}{
  \left| t_1 - t_2 \right|^{ \varepsilon }
}
\Biggr)
\\ & 
\leq 
  \| 
    X_0 
  \|_{
    \mathscr{L}^p( 
      \P ; 
      \left\| \cdot \right\|_{ 
        H_{ \min\{ \delta, r + \varepsilon \} } 
      } 
    ) 
  }
+
\biggl[
  \sup_{ s \in [0,T] }
  \|
    F( X_s ) 
  \|_{
    \mathscr{L}^p( \P; \left\| \cdot \right\|_{ H_{ \gamma - \eta } } )
  }
\biggr]
  \frac{
    2
    \,
    T^{ 
      ( 1 + \gamma - \eta - \min\{ \delta, r + \varepsilon \} ) 
    }
  }{
    ( 1 + \gamma - \eta - r - \varepsilon ) 
  }
\\ & \quad
+
\biggl[
  \sup_{ s \in [0,T] }
  \|
    B( X_s ) 
  \|_{
    \mathscr{L}^p( \P; \left\| \cdot \right\|_{ \mathrm{HS}( U, H_{ \beta } ) } )
  }
\biggr]
  \frac{ 
    \sqrt{ p \, ( p - 1 ) }
    \,
    T^{ ( \nicefrac{1}{2} + \beta - \min\{ \delta, r + \varepsilon \} ) }
  }{
    ( 
      1 + 2 \beta - 2 r - 2 \varepsilon 
    )^{ \nicefrac{1}{2} }
  }
< \infty
.
\end{split}
\end{equation}
\end{corollary}
\begin{proof}[Proof of Corollary~\ref{cor:temporal_regularity_SPDE}]
The fact that
$
\forall \, u \in [ 0, 1 ] \colon
\bigl(
\sup_{ t \in ( 0, T ] }
    t^u \| ( - A )^u e^{ t A }
         \|_{ L( H ) }
\leq 1
$
and
$
\sup_{ t \in ( 0, T ] }
    t^{ -u } \| ( - A )^{ -u } ( e^{ tA } - \operatorname{Id}_H )
         \|_{ L( H ) }
\leq 1
\bigr)
$
ensures 
for all 
$ 
r \in 
[
  \gamma ,
  \min\{ 
    1 + \gamma - \eta
    ,
    \nicefrac{ 1 }{ 2 } + \beta 
  \}
) 
$,
$
\varepsilon \in 
[
  0,
  \min\{ 
    1 + \gamma - \eta - r
    ,
    \nicefrac{ 1 }{ 2 } + \beta - r
  \}
)
$
that
\begin{align*}
& 
\sup_{ 
  \substack{ 
    t_1, t_2 \in [0,T] ,
  \\
    t_1 \neq t_2
  }
}
\Biggl(
\frac{ 
  \|
  | \min\{ t_1, t_2 \} |^{
    \max\{ r + \varepsilon - \delta , 0 \}
  }
  (
    e^{ t_1 A } X_0
    -
    e^{ t_2 A } X_0 
  )
  \|_{
    \mathscr{L}^p( 
      \P ; 
      \left\| \cdot \right\|_{ 
        H_{ r } 
      } 
    ) 
  }
}{
  \left| t_1 - t_2 \right|^{ \varepsilon }
}
\Biggr)
\\ & 
\leq 
\sup_{ 
  \substack{ 
    t_1, t_2 \in [0,T] ,
  \\
    t_1 < t_2
  }
}
\Biggl(
\frac{ 
  \|
    | t_1 |^{
      \max\{ r + \varepsilon - \delta , 0 \}
    }
    ( - A )^{
      r - \min\{ \delta, r + \varepsilon \}
    }
    (
      e^{ t_1 A } 
      -
      e^{ t_2 A }
    )
  \|_{
    L( H )
  }
  \,
  \| 
    X_0 
  \|_{
    \mathscr{L}^p( 
      \P ; 
      \left\| \cdot \right\|_{ 
        H_{ \min\{ \delta, r + \varepsilon \} } 
      } 
    ) 
  }
}{
  \left| t_1 - t_2 \right|^{ \varepsilon }
}
\Biggr)
\\ & \yesnumber
\leq 
\sup_{ 
  \substack{ 
    t_1, t_2 \in (0,T] ,
  \\
    t_1 < t_2
  }
}
\Bigl(
  | t_1 |^{
    \max\{ r + \varepsilon - \delta , 0 \}
  }
  \| 
    ( - A )^{
      r + \varepsilon - \min\{ \delta, r + \varepsilon \} 
    }
    \,
    e^{ t_1 A } 
  \|_{
    L( H )
  }
  \,
  \| 
    X_0 
  \|_{
    \mathscr{L}^p( 
      \P ; 
      \left\| \cdot \right\|_{ 
        H_{ \min\{ \delta, r + \varepsilon \} } 
      } 
    ) 
  }
\Bigr)
\\ & 
= 
\sup_{ 
    t_1 \in (0,T]
}
\Bigl(
  | t_1 |^{
    \max\{ r + \varepsilon - \delta , 0 \}
  }
  \| 
    ( - A )^{
      \max\{ r + \varepsilon - \delta, 0 \} 
    }
    \,
    e^{ t_1 A } 
  \|_{
    L( H )
  } \,
  \| 
    X_0 
  \|_{
    \mathscr{L}^p( 
      \P ; 
      \left\| \cdot \right\|_{ 
        H_{ \min\{ \delta, r + \varepsilon \} } 
      } 
    ) 
  }
\Bigr)
\\ & 
\leq 
  \| 
    X_0 
  \|_{
    \mathscr{L}^p( 
      \P ; 
      \left\| \cdot \right\|_{ 
        H_{ \min\{ \delta, r + \varepsilon \} } 
      } 
    ) 
  }
.
\end{align*}
Combining this with the triangle inequality
and Proposition~\ref{prop:temporal_regularity_SPDE} proves
for all 
$ 
r \in 
[
  \gamma ,
  \min\{ 
    1 $\linebreak$+ \gamma - \eta
    ,
    \nicefrac{ 1 }{ 2 } + \beta 
  \}
) 
$,
$
\varepsilon \in 
[
  0,
  \min\{ 
    1 + \gamma - \eta - r
    ,
    \nicefrac{ 1 }{ 2 } + \beta - r
  \}
)
$
that
$ \inf_{ s \in ( 0, T ] }
\P( X_s \in H_r ) = 1 $
and 
\begin{align*}
&
\sup_{ 
  \substack{ 
    t_1, t_2 \in [0,T] ,
  \\
    t_1 \neq t_2
  }
}
\Biggl(
\frac{ 
  \left| \min\{ t_1, t_2 \} \right|^{
    \max\{ r + \varepsilon - \delta , 0 \}
  }
  \|
  \mathbbm{1}_{ \{ X_{ t_1 } \in H_r \} }
    X_{ t_1 } 
    -
  \mathbbm{1}_{ \{ X_{ t_2 } \in H_r \} }
    X_{ t_2 } 
  \|_{
    \mathscr{L}^p( 
      \P ; 
      \left\| \cdot \right\|_{ 
        H_{ r } 
      } 
    ) 
  }
}{
  \left| t_1 - t_2 \right|^{ \varepsilon }
}
\Biggr)
\\ & \leq
\sup_{ 
  \substack{ 
    t_1, t_2 \in [0,T] ,
  \\
    t_1 \neq t_2
  }
}
\Biggl(
\frac{ 
  \|
  | \min\{ t_1, t_2 \} |^{
    \max\{ r + \varepsilon - \delta , 0 \}
  }
  (
    e^{ t_1 A } X_0
    -
    e^{ t_2 A } X_0 
  )
  \|_{
    \mathscr{L}^p( 
      \P ; 
      \left\| \cdot \right\|_{ 
        H_{ r } 
      } 
    ) 
  }
}{
  \left| t_1 - t_2 \right|^{ \varepsilon }
}
\Biggr)
+
T^{
    \max\{ r + \varepsilon - \delta , 0 \}
}
\\ & \quad \cdot \yesnumber
\sup_{ 
  \substack{ 
    t_1, t_2 \in [0,T] ,
  \\
    t_1 \neq t_2
  }
}
\Biggl(
\frac{ 
  \|
    (
      X_{ t_1 } - e^{ t_1 A } X_0
    ) \,
  \mathbbm{1}_{ \{ X_{ t_1 } \in H_r \} }
    -
    (
      X_{ t_2 } - e^{ t_2 A } X_0
    ) \,
  \mathbbm{1}_{ \{ X_{ t_2 } \in H_r \} }
  \|_{
    \mathscr{L}^p( 
      \P ; 
      \left\| \cdot \right\|_{ 
        H_{ r } 
      } 
    ) 
  }
}{
  \left| t_1 - t_2 \right|^{ \varepsilon }
}
\Biggr)
\\ & \leq
\| X_0 \|_{
\mathscr{L}^p( 
  \P ; 
  \left\| \cdot \right\|_{ 
    H_{ \min\{ \delta, r + \varepsilon \} } 
  } 
) 
}
+
\biggl[
  \sup_{ s \in [0,T] }
  \|
    F( X_s ) 
  \|_{
    \mathscr{L}^p( \P; \left\| \cdot \right\|_{ H_{ \gamma - \eta } } )
  }
\biggr]
  \frac{
    2
    \,
    T^{ 
      ( 1 + \gamma - \eta - \min\{ \delta, r + \varepsilon \} ) 
    }
  }{
    ( 1 + \gamma - \eta - r - \varepsilon ) 
  }
\\ & \quad
+
\biggl[
  \sup_{ s \in [0,T] }
  \|
    B( X_s ) 
  \|_{
    \mathscr{L}^p( \P; \left\| \cdot \right\|_{ \mathrm{HS}( U, H_{ \beta } ) } )
  }
\biggr]
  \frac{ 
    \sqrt{ p \, ( p - 1 ) }
    \,
    T^{ ( \nicefrac{1}{2} + \beta - \min\{ \delta, r + \varepsilon \} ) }
  }{
    ( 
      1 + 2 \beta - 2 r - 2 \varepsilon 
    )^{ \nicefrac{1}{2} }
  }
< \infty
.
\end{align*}
The proof 
of Corollary~\ref{cor:temporal_regularity_SPDE}
is thus completed.
\end{proof}
}{}
\iftoggle{arXiv:v3}{%
\subsection{A priori bounds for solutions of SEEs}
\begin{lemma}
\label{lem:generalized_exponentials}
Consider the notation in Subsection~\ref{notation}
and let
$ \mathbb{B} \colon ( 0, \infty )^2 \to ( 0, \infty ) $
and
$ \mathrm{E}_\eta \colon [ 0, \infty ) \to [ 0, \infty ) $, $ \eta \in ( -\infty, 1 ) $,
be the functions which satisfy
for all $ \eta \in ( -\infty, 1 ) $, $ x, y \in ( 0, \infty ) $, $ z \in [ 0, \infty ) $ that
$ \mathbb{B}( x, y )
= \int_{0}^{1}  t^{ ( x - 1 ) }  ( 1 - t )^{ ( y - 1 ) }  \ud t $
and
$ \mathrm{E}_\eta( z )
= 1 +
\sum_{n=1}^{\infty}
    z^n
    \prod_{k=0}^{n-1}
        \mathbb{B}( 1 - \eta, k ( 1 - \eta ) + 1 ) $.
Then it holds
for all $ \eta \in ( -\infty, 1 ) $, $ x \in [ 0, \infty ) $
that
$
\sqrt{ \mathrm{E}_\eta ( x^2 ) }
=
\mathscr{E}_{ ( 1 - \eta ) }( x )
$.
\end{lemma}
\begin{proof}[Proof of Lemma~\ref{lem:generalized_exponentials}]
Note that
the fact that
$
\forall \, x, y \in ( 0, \infty ) \colon
\mathbb{B}( x, y )
=
\frac{ \Gamma(x) \, \Gamma( y ) }{ \Gamma( x + y ) }
$
implies that it holds
for all $ \eta \in ( -\infty, 1 ) $, $ x \in [ 0, \infty ) $
that
\begin{equation}
\begin{split}
\mathrm{E}_\eta ( x^2 )
& =
1 +
\sum_{n=1}^{\infty}
    ( x^2 )^n
    \prod_{k=0}^{n-1}
        \mathbb{B}( 1 - \eta, k ( 1 - \eta ) + 1 )
\\ & =
1 +
\sum_{n=1}^{\infty}
    x^{2n}
    \prod_{k=0}^{n-1}
        \frac{ \Gamma( 1 - \eta ) \, \Gamma( k ( 1 - \eta ) + 1 )
            }{ \Gamma( ( k + 1 ) ( 1 - \eta ) + 1 ) }
=
1 +
\sum_{n=1}^{\infty}
    \frac{ x^{2n} [ \Gamma( 1 - \eta ) ]^n
        }{ \Gamma( n ( 1 - \eta ) + 1 )  }
\\ & =
\sum_{n=0}^{\infty}
    \frac{ x^{2n} [ \Gamma( 1 - \eta ) ]^n
        }{ \Gamma( n ( 1 - \eta ) + 1 )  }
=
[ \mathscr{E}_{ ( 1 - \eta ) }( x ) ]^2.
\end{split}
\end{equation}
The proof of Lemma~\ref{lem:generalized_exponentials} is thus completed.
\end{proof}
\begin{proposition}[A priori bounds]
\label{prop:regularity_SPDE2}
Consider the notation in Subsection~\ref{notation},
let\linebreak
$ ( H, \left< \cdot, \cdot \right>_H, \left\| \cdot \right\|_H ) $
and
$ ( U, \left< \cdot, \cdot \right>_U, \left\| \cdot \right\|_U ) $
be separable $ \R $-Hilbert spaces,
let $ \mathbb{H} \subseteq H $ be a non-empty orthonormal basis of $ H $,
let $ \lambda \colon \mathbb{H} \to \mathbb{R} $
be a function with 
$
  \sup_{ h \in \mathbb{H} }
  \lambda_h < 0
$,
let
$
  A \colon D(A) \subseteq H \rightarrow H
$
be the linear operator which satisfies
$
  D(A) 
  = 
    \bigl\{ 
      v \in H 
      \colon $\linebreak
      $
      \sum_{ h \in \mathbb{H} } 
	\left| 
	  \lambda_h 
	  \langle h, v \rangle_H 
	\right|^2
	< \infty
    \bigr\}
$
and
$
\forall \, v \in D(A) \colon
    Av
   =
    \sum_{ h \in \mathbb{H} } 
    \lambda_h \langle h, v \rangle_H h
$,
let 
$ 
  ( H_r , \left< \cdot , \cdot \right>_{ H_r }, \left\| \cdot \right\|_{ H_r } ) 
$,
$ r \in \R $,
be a family of interpolation spaces associated to
$ - A $,
let 
$ T \in (0,\infty) $,
$ p \in [2,\infty) $,
$ \gamma \in \R $,
$ \eta \in [0, 1) $,
$ F \in C( H_{ \gamma } , H_{ \gamma - \eta } ) $,
$ B \in C( H_{ \gamma } , \mathrm{HS}( U, H_{ \gamma - \nicefrac{ \eta }{ 2 } } ) ) $
satisfy
$
| F |_{
    \calC^1( H_{ \gamma }, \left\| \cdot \right\|_{  H_{ \gamma - \eta }  } )
}
+
| B |_{
    \calC^1( H_{ \gamma }, \left\| \cdot \right\|_{  \mathrm{HS}( U, H_{ \gamma - \nicefrac{ \eta }{ 2 } } )  } )
}
< \infty
$,
let 
$ ( \Omega, \calF, \P ) $
be a probability space with a normal filtration 
$ ( \calF_t )_{ t \in [0,T] } $,
let 
$
  ( W_t )_{ t \in [0,T] } 
$ 
be an $ \operatorname{Id}_U $-cylindrical
$ ( \Omega, \calF, \P, ( \calF_t )_{ t \in [0,T] } ) $-Wiener
process,
and let 
$ X \colon [0,T] \times \Omega \to H_{ \gamma } $
be an 
$ ( \mathscr{F}_t )_{ t \in [0,T] } $/$ \mathscr{B}( H_{ \gamma } ) $-predictable stochastic process
which satisfies
for all $ t \in [ 0, T ] $ that
$ 
\sup_{ s \in [0,T] } 
\| X_s \|_{ \mathscr{L}^p( \P ; \| \cdot \|_{ H_{ \gamma } } ) } < \infty 
$
and
\begin{equation}
\begin{split}
[ X_t ]_{ \P, \mathscr{B}( H_\gamma ) }
& = 
\biggl[
    e^{ t A } X_0
    +
    \int_0^t
        \mathbbm{1}_{
          \{ 
            \int_0^t
            \| e^{ ( t - u ) A } F( X_u ) \|_{ H_{ \gamma } } \dd u
            < \infty
          \} } \,
        e^{ ( t - s ) A }
        F( X_s )
    \ds
\biggr]_{ \P, \mathscr{B}( H_\gamma ) }
\\ & \quad +
\int_0^t
    e^{ ( t - s ) A }
    B( X_s )
\dWs
.
\end{split}
\end{equation}
Then
\begin{align}
&
\sup_{ t \in [0,T] }
\left\|
 \max \bigl\{ 1,
   \| X_t \|_{ H_{ \gamma } }
 \bigr\}
\right\|_{ \mathscr{L}^p( \P ; \left| \cdot \right| ) }
\leq
\sqrt{2}
\left\|
 \max \bigl\{ 1,
   \| X_0 \|_{ H_{ \gamma } }
 \bigr\}
\right\|_{ \mathscr{L}^p( \P ; \left| \cdot \right| ) }
\\ & \quad \cdot \nonumber
\mathscr{E}_{ (1 - \eta) } \biggl[ 
\tfrac{ 
  T^{ 1 - \eta }
  \sqrt{ 2 }
}{
  \sqrt{ 1 - \eta }
} \,
\biggl(
  \sup_{ v \in H_{ \gamma } }
  \tfrac{
    \| F( v ) \|_{ H_{ \gamma - \eta } }
  }{
    \max\{ 1, \| v \|_{ H_{ \gamma } } \}
  }
\biggr)
+
\sqrt{
  T^{ 1 - \eta } p ( p - 1 )
} \,
\biggl(
  \sup_{ v \in H_{ \gamma } }
  \tfrac{
    \| B( v ) \|_{ \mathrm{HS}( U, H_{ \gamma - \nicefrac{\eta}{2} } ) }
  }{
    \max\{ 1, \| v \|_{ H_{ \gamma } } \}
  }
\biggr)
\biggr]
< \infty
.
\end{align}
\end{proposition}
\begin{proof}[Proof of Proposition~\ref{prop:regularity_SPDE2}]
The Burkholder-Davis-Gundy-type inequality
in Lemma~7.7 in Da Prato \& Zabczyk~\cite{dz92},
the fact that
$
\forall \, u \in [ 0, 1 ] \colon
\sup_{ t \in ( 0, T ] }
    t^u \| ( - A )^u e^{ t A }
         \|_{ L( H ) }
\leq 1
$,
and
H\"older's inequality
imply that it holds
for all $ t \in [0,T] $
that
\begin{equation}
\begin{split}
&
\bigl\|
 \max \bigl\{ 1,
   \| X_t \|_{ H_{ \gamma } }
 \bigr\}
\bigr\|_{ \mathscr{L}^p( \P ; \left| \cdot \right| ) }
\\ & \leq
\bigl\|
 \max \bigl\{ 1,
   \| X_0 \|_{ H_{ \gamma } }
 \bigr\}
\bigr\|_{ \mathscr{L}^p( \P ; \left| \cdot \right| ) }
+
\int_0^t
\bigl\| 
 e^{ ( t - s ) A }
 F( X_s )
\bigr\|_{ 
 \mathscr{L}^p( \P ; \left\| \cdot \right\|_{ H_{ \gamma } } )
}
\ds
\\ & \quad +
\sqrt{ 
 \frac{ p \, ( p - 1 ) }{ 2 }
}
\biggl[
\int_0^t
\bigl\| 
 e^{ ( t - s ) A }
 B( X_s )
\bigr\|_{ 
 \mathscr{L}^p( \P ; \left\| \cdot \right\|_{ \mathrm{HS}( U, H_{ \gamma } ) } )
}^2
\ds
\biggr]^{ 
 \nicefrac{ 1 }{ 2 }
}
\\ & \leq
\bigl\|
 \max \bigl\{ 1,
   \| X_0 \|_{ H_{ \gamma } }
 \bigr\}
\bigr\|_{ \mathscr{L}^p( \P ; \left| \cdot \right| ) }
+
\biggl[
\frac{
 t^{ ( 1 - \eta ) }
}{
 ( 1 - \eta )
}
\int_0^t
( t - s )^{ - \eta } \,
\| 
 F( X_s )
\|_{ 
 \mathscr{L}^p( \P ; \left\| \cdot \right\|_{ H_{ \gamma - \eta } } )
}^2
\ds
\biggr]^{ \nicefrac{ 1 }{ 2 } }
\\ & \quad +
\sqrt{ 
 \frac{ p \, ( p - 1 ) }{ 2 }
}
\biggl[
\int_0^t
( t - s )^{
 - \eta
} \,
\| 
 B( X_s )
\|_{ 
 \mathscr{L}^p( \P ; \left\| \cdot \right\|_{ \mathrm{HS}( U, H_{ \gamma - \nicefrac{\eta}{2} } ) } )
}^2
\ds
\biggr]^{ 
 \nicefrac{ 1 }{ 2 }
}
\\ & \leq
\bigl\|
 \max \bigl\{ 1,
   \| X_0 \|_{ H_{ \gamma } }
 \bigr\}
\bigr\|_{ \mathscr{L}^p( \P ; \left| \cdot \right| ) }
+
\biggl[
\int_0^t
( t - s )^{ - \eta } \,
\bigl\|
 \max \bigl\{ 1,
   \| X_s \|_{ H_{ \gamma } }
 \bigr\}
\bigr\|_{ \mathscr{L}^p( \P ; \left| \cdot \right| ) }^2
\ds
\biggr]^{ \nicefrac{ 1 }{ 2 } }
\\ & \quad
\cdot
\biggl[
\sqrt{
\tfrac{
 T^{ ( 1 - \eta ) }
}{
 ( 1 - \eta )
}
} \,
\biggl(
  \sup_{ v \in H_{ \gamma } }
  \tfrac{
    \| F( v ) \|_{ H_{ \gamma - \eta } }
  }{
    \max\{ 1, \| v \|_{ H_{ \gamma } } \}
  }
\biggr)
+
\sqrt{ 
 \tfrac{ p \, ( p - 1 ) }{ 2 }
} \,
\biggl(
  \sup_{ v \in H_{ \gamma } }
  \tfrac{
    \| B( v ) \|_{ \mathrm{HS}( U, H_{ \gamma - \nicefrac{\eta}{2} } ) }
  }{
    \max\{ 1, \| v \|_{ H_{ \gamma } } \}
  }
\biggr)
\biggr]
.
\end{split}
\end{equation}
This and the fact that
$
\forall \, a, b \in \R \colon
( a + b )^2 \leq 2 \, ( a^2 + b^2 )
$
prove
for all $ t \in [0,T] $
that
\begin{equation}
\begin{split}
&
\bigl\|
 \max \bigl\{ 1,
   \| X_t \|_{ H_{ \gamma } }
 \bigr\}
\bigr\|_{ \mathscr{L}^p( \P ; \left| \cdot \right| ) }^2
\leq
2 \,
\bigl\|
 \max \bigl\{ 1,
   \| X_0 \|_{ H_{ \gamma } }
 \bigr\}
\bigr\|_{ \mathscr{L}^p( \P ; \left| \cdot \right| ) }^2
\\ & \quad
+
\int_0^t
( t - s )^{ - \eta } \,
\bigl\|
 \max \bigl\{ 1,
   \| X_s \|_{ H_{ \gamma } }
 \bigr\}
\bigr\|_{ \mathscr{L}^p( \P ; \left| \cdot \right| ) }^2
\ds
\\ & \quad
\cdot
\biggl[
\sqrt{
\tfrac{
 2 \, T^{ ( 1 - \eta ) }
}{
 ( 1 - \eta )
}
} \,
\biggl(
  \sup_{ v \in H_{ \gamma } }
  \tfrac{
    \| F( v ) \|_{ H_{ \gamma - \eta } }
  }{
    \max\{ 1, \| v \|_{ H_{ \gamma } } \}
  }
\biggr)
+
\sqrt{ 
 p \, ( p - 1 ) 
} \,
\biggl(
  \sup_{ v \in H_{ \gamma } }
  \tfrac{
    \| B( v ) \|_{ \mathrm{HS}( U, H_{ \gamma - \nicefrac{\eta}{2} } ) }
  }{
    \max\{ 1, \| v \|_{ H_{ \gamma } } \}
  }
\biggr)
\biggr]^2
.
\end{split}
\end{equation}
E.g., Lemma~2.6 in
Andersson, Jentzen, \& Kurniawan~\cite{AnderssonJentzenKurniawan2015arXiv}
and
Lemma~\ref{lem:generalized_exponentials}
hence complete the proof
of Proposition~\ref{prop:regularity_SPDE2}.
\end{proof}
}{}
\iftoggle{arXiv:v3}{%
\subsection{A strong perturbation estimate for SEEs}
\begin{proposition}[Perturbation estimate]
\label{prop:strong_perturbation}
Consider the notation in Subsection~\ref{notation},
let 
$ ( H, \left< \cdot, \cdot \right>_H, \left\| \cdot \right\|_H ) $
and
$ ( U, \left< \cdot, \cdot \right>_U, \left\| \cdot \right\|_U ) $
be separable $ \R $-Hilbert spaces,
let $ \mathbb{H} \subseteq H $ be a non-empty orthonormal basis of $ H $,
let $ \lambda \colon \mathbb{H} \to \mathbb{R} $
be a function with 
$
  \sup_{ h \in \mathbb{H} }
  \lambda_h < 0
$,
let
$
  A \colon D(A) \subseteq H \rightarrow H
$
be the linear operator which satisfies
$
  D(A) 
  = 
    \bigl\{ 
      v \in H 
      \colon $\linebreak
      $
      \sum_{ h \in \mathbb{H} } 
	\left| 
	  \lambda_h 
	  \langle h, v \rangle_H 
	\right|^2
	< \infty
    \bigr\}
$
and
$
\forall \, v \in D(A) \colon
    Av
   =
    \sum_{ h \in \mathbb{H} } 
    \lambda_h \langle h, v \rangle_H h
$,
let 
$ 
  ( H_r , \left< \cdot , \cdot \right>_{ H_r }, \left\| \cdot \right\|_{ H_r } ) 
$,
$ r \in \R $,
be a family of interpolation spaces associated to
$ - A $,
let 
$ T \in [0,\infty) $,
$ p \in [2,\infty) $,
$ \gamma \in \R $,
$ \eta \in [0, 1) $,
$ F \in C( H_{ \gamma } , H_{ \gamma - \eta } ) $,
$ B \in C( H_{ \gamma } , \mathrm{HS}( U, H_{ \gamma - \nicefrac{ \eta }{ 2 } } ) ) $
satisfy
$
| F |_{
    \calC^1( H_{ \gamma }, \left\| \cdot \right\|_{  H_{ \gamma - \eta }  } )
}
+
| B |_{
    \calC^1( H_{ \gamma }, \left\| \cdot \right\|_{  \mathrm{HS}( U, H_{ \gamma - \nicefrac{ \eta }{ 2 } } )  } )
}
< \infty
$,
let 
$ ( \Omega, \calF, \P ) $
be a probability space with a normal filtration 
$ ( \calF_t )_{ t \in [0,T] } $,
let 
$
  ( W_t )_{ t \in [0,T] } 
$ 
be an $ \operatorname{Id}_U $-cylindrical
$ ( \Omega, \calF, \P, ( \calF_t )_{ t \in [0,T] } ) $-Wiener
process,
and let 
$
X^1, X^2 \colon [0,T] \times \Omega \to H_{ \gamma }
$
be $ ( \mathscr{F}_t )_{ t \in [0,T] } $/$ \mathscr{B}( H_{ \gamma } ) $-predictable stochastic processes
which satisfy
$
\max_{ k \in \{ 1, 2 \} }
\sup_{ s \in [0,T] }
\| X_s^k \|_{
 \mathscr{L}^p( \P ; \left\| \cdot \right\|_{ H_{ \gamma } } )
}
< \infty
$. 
Then 
\begin{align}
\label{eq:perturbation}
& \nonumber
\sup_{ t \in [0,T] }
\left\|
 X^1_t - X^2_t
\right\|_{ 
 \mathscr{L}^p( 
   \P ; 
   \left\| \cdot \right\|_{ 
     H_{ \gamma } 
   } 
 ) 
}
\\ & \leq \nonumber
\mathscr{E}_{ (1 - \eta) }\!\left[ 
\tfrac{ 
 T^{ 1 - \eta }
 \sqrt{ 2 }
 \,
 \left| F \right|_{
   \calC^1( 
     H_{ \gamma }, \left\| \cdot \right\|_{ H_{ \gamma - \eta } }
   )
 }
}{ \sqrt{ 1 - \eta } }
 +
 \sqrt{ 
   T^{ 1 - \eta }
   p ( p - 1 ) 
 } 
\left| B \right|_{
 \calC^1( 
   H_{ \gamma }, \left\| \cdot \right\|_{ \mathrm{HS}( U, H_{ \gamma - \nicefrac{ \eta }{ 2 } } ) }
 )
}
\right]
\\ & \quad 
\cdot
\sqrt{2}
\,
\sup_{ t \in [0,T] }
\bigg\|
 \biggl[ X_t^1
 -
   \int_0^t
    \mathbbm{1}_{
      \{ 
        \int_0^t
        \| e^{ ( t - r ) A } F( X^1_r ) \|_{ H_{ \gamma } } \dd r
        < \infty
      \}
    }
   e^{ ( t - s ) A } F( X^1_s ) 
   \ds
 \biggr]_{ \P, \mathscr{B}( H_{ \gamma } ) }
\\ & \quad \nonumber
   -
   \int_0^t
   e^{ ( t - s ) A } B( X^1_s ) 
   \dWs
-
\biggl\{
 \biggl[ X_t^2
 -
 \int_0^t
 \mathbbm{1}_{
   \{ 
     \int_0^t
     \| e^{ ( t - r ) A } F( X^2_r ) \|_{ H_{ \gamma } } \dd r
     < \infty
   \}
 }
   e^{ ( t - s ) A } F( X^2_s ) 
 \ds
 \biggr]_{ \P, \mathscr{B}( H_{ \gamma } ) }
\\ & \quad \nonumber
 -
 \int_0^t
   e^{ ( t - s ) A } B( X^2_s ) 
 \dWs
\biggr\}
\biggr\|_{ 
 L^p( 
   \P ; 
   \left\| \cdot \right\|_{ H_{ \gamma } }
 ) 
}
< \infty
.
\end{align}
\end{proposition}
\begin{proof}[Proof of Proposition~\ref{prop:strong_perturbation}]
Throughout this proof we assume w.l.o.g.\ that $ T \neq 0 $
and throughout this proof let
$ \mathscr{A} \colon H_{ \gamma + 1 } \subseteq H_\gamma \rightarrow H_\gamma $
be the linear operator which satisfies
for all $ v \in H_{ \gamma + 1 } $ that
$
\mathscr{A} v
=
\sum_{ h \in \mathbb{H} } 
\lambda_h \langle ( - \lambda_h )^{ - \gamma } h, v \rangle_{ H_\gamma }
( - \lambda_h )^{ - \gamma } h
$.
Observe that
$ 
  ( H_{ r + \gamma } , \left< \cdot , \cdot \right>_{ H_{ r + \gamma }  }, \left\| \cdot \right\|_{ H_{ r + \gamma } } ) 
$,
$ r \in \R $,
is a family of interpolation spaces associated to
$ - \mathscr{A} $.
This,
Lemma~\ref{lem:generalized_exponentials},
and Proposition~2.7 in
Andersson, Jentzen, \& Kurniawan~\cite{AnderssonJentzenKurniawan2015arXiv}
show for all $ \varepsilon \in ( 0, \infty ) $ that
\begin{align}
& \nonumber
\sup_{ t \in ( 0,T] }
\left\|
 X^1_t - X^2_t
\right\|_{ 
 \mathscr{L}^p( 
   \P ; 
   \left\| \cdot \right\|_{ 
     H_{ \gamma } 
   } 
 ) 
}
\\ & \leq \nonumber
\mathscr{E}_{ (1 - \eta) } \biggl[ 
\tfrac{
 T^{ 1 - \eta }
 \sqrt{ 2 }
 \,
 \left| F \right|_{
   \calC^1( 
     H_{ \gamma }, \left\| \cdot \right\|_{ H_{ \gamma - \eta } }
   )
 }
}{ \sqrt{ 1 - \eta } }
\sup_{ t \in ( 0, T ] }
     t^\eta \| ( - \mathscr{A} )^\eta e^{ t \mathscr{A} }
          \|_{ L( H_{ \gamma } ) }
\\ & \quad \nonumber
 +
 \sqrt{ 
   T^{ 1 - \eta }
   p ( p - 1 ) 
 }
 \bigl(
\left| B \right|_{
 \calC^1( 
   H_{ \gamma }, \left\| \cdot \right\|_{ \mathrm{HS}( U, H_{ \gamma - \nicefrac{ \eta }{ 2 } } ) }
 )
}
+ \varepsilon
\bigr)
\sup_{ t \in ( 0, T ] }
     t^{ \nicefrac{\eta}{2} } \| ( - \mathscr{A} )^{ \nicefrac{\eta}{2} } e^{ t \mathscr{A} }
          \|_{ L( H_{ \gamma } ) }
\biggr]
\\ & \quad 
\cdot
\sqrt{2}
\,
\sup_{ t \in ( 0,T] }
\bigg\|
 \biggl[ X_t^1
 -
   \int_0^t
    \mathbbm{1}_{
      \{ 
        \int_0^t
        \| e^{ ( t - r ) \mathscr{A} } F( X^1_r ) \|_{ H_{ \gamma } } \dd r
        < \infty
      \}
    }
   e^{ ( t - s ) \mathscr{A} } F( X^1_s ) 
   \ds
 \biggr]_{ \P, \mathscr{B}( H_{ \gamma } ) }
\\ & \quad \nonumber
   -
   \int_0^t
   e^{ ( t - s ) \mathscr{A} } B( X^1_s ) 
   \dWs
-
\biggl\{
 \biggl[ X_t^2
 -
 \int_0^t
 \mathbbm{1}_{
   \{ 
     \int_0^t
     \| e^{ ( t - r ) \mathscr{A} } F( X^2_r ) \|_{ H_{ \gamma } } \dd r
     < \infty
   \}
 }
   e^{ ( t - s ) \mathscr{A} } F( X^2_s ) 
 \ds
 \biggr]_{ \P, \mathscr{B}( H_{ \gamma } ) }
\\ & \quad \nonumber
 -
 \int_0^t
   e^{ ( t - s ) \mathscr{A} } B( X^2_s ) 
 \dWs
\biggr\}
\biggr\|_{ 
 L^p( 
   \P ; 
   \left\| \cdot \right\|_{ H_{ \gamma } }
 ) 
}
.
\end{align}
The fact that
$
\forall \, u \in [ 0, 1 ] \colon
\sup_{ t \in ( 0, T ] }
    t^u \| ( - \mathscr{A} )^u e^{ t \mathscr{A} }
         \|_{ L( H_{ \gamma } ) }
\leq 1
$
hence proves \eqref{eq:perturbation}.
The proof of Proposition~\ref{prop:strong_perturbation} is thus completed.
\end{proof}
}{}
\iftoggle{arXiv:v3}{%
\subsection{Existence of continuous solutions of SEEs}
The next result, Proposition~\ref{prop:existence_continuous_SPDEs} below,
proves the existence of continuous solution processes of SPDEs
(see, e.g., Theorem~7.1 in
van~Neerven, Veraar, \& Weis~\cite{VanNeervenVeraarWeis2008}
for a similar result in a more general framework).
\begin{proposition}
\label{prop:existence_continuous_SPDEs}
Consider the notation in Subsection~\ref{notation},
let 
$
( 
 H,
 \langle \cdot, \cdot \rangle_H, 
 \left\| \cdot \right\|_H
)
$
and\linebreak
$
( 
 U ,
 \langle \cdot, \cdot \rangle_U , 
 \left\| \cdot \right\|_U
)
$
be separable $ \mathbb{R} $-Hilbert spaces, 
let $ \mathbb{H} \subseteq H $ be a non-empty orthonormal basis of $ H $,
let 
$ T \in (0,\infty) $,
$ p \in [2,\infty) $,
let 
$ ( \Omega, \mathscr{F}, \P ) $
be a probability space with a normal filtration 
$ ( \mathscr{F}_t )_{ t \in [0,T] } $, 
let 
$
( W_t )_{ t \in [0,T] } 
$ 
be an $ \operatorname{Id}_U $-cylindrical 
$ ( \Omega, \mathscr{F}, \P, ( \mathscr{F}_t )_{ t \in [0,T] } ) $-Wiener
process,
let $ \lambda \colon \mathbb{H} \to \mathbb{R} $
be a function with 
$
  \sup_{ h \in \mathbb{H} }
  \lambda_h < 0
$,
let
$
  A \colon D(A) \subseteq H \rightarrow H
$
be the linear operator which satisfies
$
  D(A) 
  = 
    \bigl\{ 
      v \in H 
      \colon
      \sum_{ h \in \mathbb{H} } 
	\left| 
	  \lambda_h 
	  \langle h, v \rangle_H 
	\right|^2
	< \infty
    \bigr\}
$
and
$
\forall \, v \in D(A) \colon
    Av
   =
    \sum_{ h \in \mathbb{H} } 
    \lambda_h \langle h, v \rangle_H h
$,
let 
$ 
( H_r , \left< \cdot , \cdot \right>_{ H_r }, \left\| \cdot \right\|_{ H_r } ) 
$,
$ r \in \R $,
be a family of interpolation spaces associated to
$ - A $,
and let 
$ \gamma \in \R $, $ \eta \in [0,1) $,
$ F \in C( H_{ \gamma } , H_{ \gamma - \eta } ) $,
$ B \in C( H_{ \gamma } , \mathrm{HS}( U, H_{ \gamma - \nicefrac{ \eta }{ 2 } } ) ) $,
$ \xi \in \mathscr{M}( \mathscr{F}_0, \mathscr{B}( H_{ \gamma } ) ) $
satisfy
$
| F |_{
    \calC^1( H_{ \gamma }, \left\| \cdot \right\|_{  H_{ \gamma - \eta }  } )
}
+
| B |_{
    \calC^1( H_{ \gamma }, \left\| \cdot \right\|_{  \mathrm{HS}( U, H_{ \gamma - \nicefrac{ \eta }{ 2 } } )  } )
}
< \infty
$.
Then there exists an 
$
( \mathscr{F}_t )_{ t \in [0,T] }
$/$ \mathscr{B}( H_{ \gamma } ) $-adapted stochastic process
$
X \colon [0,T] \times \Omega \to H_{ \gamma }
$
with continuous sample paths which satisfies
for all $ t \in [0,T] $ 
that
$
[ X_t ]_{ \P, \mathscr{B}( H_\gamma ) }
=
\bigl[
    e^{ t A } \xi
    +
    \int_0^t
        e^{ ( t - s ) A }
        F( X_s )
    \ds
\bigr]_{ \P, \mathscr{B}( H_\gamma ) }
+
\int_0^t
    e^{ ( t - s ) A }
    B( X_s ) 
\dWs
$
and
\begin{equation}
\begin{split}
&
\sup_{ t \in [0,T] }
\left\|
 \max \bigl\{ 1,
   \| X_t \|_{ H_{ \gamma } }
 \bigr\}
\right\|_{ \mathscr{L}^p( \P ; \left| \cdot \right| ) }
\leq
\sqrt{2}
\left\|
 \max \bigl\{ 1,
   \| \xi \|_{ H_{ \gamma } }
 \bigr\}
\right\|_{ \mathscr{L}^p( \P ; \left| \cdot \right| ) }
\\ & \cdot
\mathscr{E}_{ (1 - \eta) } \biggl[ 
\tfrac{ 
  T^{ 1 - \eta }
  \sqrt{ 2 }
}{
  \sqrt{ 1 - \eta }
} \,
\biggl(
  \sup_{ v \in H_{ \gamma } }
  \tfrac{
    \| F( v ) \|_{ H_{ \gamma - \eta } }
  }{
    \max\{ 1, \| v \|_{ H_{ \gamma } } \}
  }
\biggr)
+
\sqrt{
  T^{ 1 - \eta } p ( p - 1 )
} \,
\biggl(
  \sup_{ v \in H_{ \gamma } }
  \tfrac{
    \| B( v ) \|_{ \mathrm{HS}( U, H_{ \gamma - \nicefrac{\eta}{2} } ) }
  }{
    \max\{ 1, \| v \|_{ H_{ \gamma } } \}
  }
\biggr)
\biggr]
.
\end{split}
\end{equation}
\end{proposition}
\begin{proof}[Proof 
of Proposition~\ref{prop:existence_continuous_SPDEs}]
Throughout this proof
let 
$ \Omega_n \in \mathscr{F}_0 $, $ n \in \N_0 $,
be the sets which satisfy
for all $ n \in \N_0 $ that
$
\Omega_n = \{ \| \xi \|_{ H_{ \gamma } } < n \}
$
and let
$ \xi_n \colon \Omega \to H_{ \gamma } $,
$ n \in \N $,
be the mappings which satisfy
for all $ n \in \N $ that
$
\xi_n
=
\xi \, \mathbbm{1}_{ \Omega_n } 
$.
Note that it holds
for all $ q \in [0,\infty) $, $ n \in \N $
that
$
\E\bigl[ \| \xi_n \|^q_{ H_{ \gamma } } \bigr] \leq n^q < \infty
$.
E.g., Theorem~5.1 in Jentzen \& Kloeden~\cite{jk12},
Proposition~\ref{prop:temporal_regularity_SPDE},
and the Kolmogorov-Chentsov continuity theorem (see Theorem~\ref{thm:Kolmogorov})
hence ensure that there exist 
$ ( \mathscr{F}_t )_{ t \in [0,T] } $/$ \mathscr{B}( H_{ \gamma } ) $-adapted 
stochastic processes with continuous sample paths 
$
X^n \colon [0,T] \times \Omega \to H_{ \gamma }
$,
$ n \in \N $,
which satisfy
for all $ n \in \N $, $ t \in [0,T] $
that
$
\sup_{ s \in [0,T] }
\E\bigl[ \| X^n_s \|_{ H_{ \gamma } }^p \bigr]
< \infty
$
and
\begin{equation}
[ X^n_t ]_{ \P, \mathscr{B}( H_\gamma ) }
=
\biggl[
e^{ t A } \xi_n 
+ \int_0^t
e^{ ( t - s ) A } F( X^n_s ) \ds
\biggr]_{ \P, \mathscr{B}( H_\gamma ) }
+ \int_0^t
e^{ ( t - s ) A } B( X^n_s ) \dWs
.
\end{equation}
Observe that it holds
for all $ k \in \N $,
$ n, m \in \{ k, k + 1, \ldots \} $,
$ t \in [0,T] $
that
\begin{equation}
\begin{split}
[
( X^n_t - X^m_t ) \, \mathbbm{1}_{
 \Omega_k
}
]_{ \P, \mathscr{B}( H_\gamma ) }
& =
\biggl[
\int_0^t
e^{ ( t - s ) A } 
\bigl[ 
 F(   
   \mathbbm{1}_{
     \Omega_k
   }
   X^n_s 
 ) 
 -
 F(   
   \mathbbm{1}_{
     \Omega_k
   }
   X^m_s 
 ) 
\bigr]
\mathbbm{1}_{
 \Omega_k
} \ds
\biggr]_{ \P, \mathscr{B}( H_\gamma ) }
\\ & \quad
+
\int_0^t
e^{ ( t - s ) A } 
\bigl[ 
 B(   
   \mathbbm{1}_{
     \Omega_k
   }
   X^n_s 
 ) 
 -
 B(   
   \mathbbm{1}_{
     \Omega_k
   }
   X^m_s 
 ) 
\bigr] 
\mathbbm{1}_{
 \Omega_k
}
\dWs
.
\end{split}
\end{equation}
Proposition~2.1 
in Jentzen \& Kurniawan~\cite{JentzenKurniawan2015arXiv}
hence shows
for all $ k \in \N $,
$ n, m \in \{ k, k + 1, \ldots \} $
that
\begin{equation}
\sup_{ t \in [0,T] }
\left\|
 (
   X^n_t
   -
   X^m_t
 ) \,
  \mathbbm{1}_{
    \Omega_k
  }
\right\|_{
 \mathscr{L}^p( \P; \left\| \cdot \right\|_{ H_{ \gamma } } )
}
= 0 .
\end{equation}
This implies that
\begin{equation}
\label{eq:XX_measure_1}
\P \biggl(
 \forall \, k \in \N \colon
 \forall \, n, m \in \{ k , k + 1 , \ldots \} 
 \colon
  \mathbbm{1}_{
    \Omega_k
  }
 \biggl[
   \sup_{ t \in [0,T] }
   \left\| X^n_t - X^m_t \right\|_{ H_{ \gamma } }
 \biggr]
 = 0
\biggr) = 1 .
\end{equation}
Next let 
$ Y \colon [0,T] \times \Omega \to H_{ \gamma } $
be the mapping which satisfies
for all $ (t,\omega) \in [0,T] \times \Omega $
that
\begin{equation}
Y_t( \omega )
=
\sum_{ n = 1 }^{ \infty }
X_t^n( \omega )
\cdot
\mathbbm{1}_{
 \Omega_n \setminus \Omega_{ n - 1 }
}( \omega )
.
\end{equation}
Note that it holds for all $ n \in \N $ that
\begin{equation}
\begin{split}
&
\mathbbm{1}_{
    \Omega_n
}
\sup_{ t \in [0,T] }
\|
    Y_t - X^n_t
\|_{ H_{ \gamma } }
=
\sup_{ t \in [0,T] }
\|
Y_t \,
\mathbbm{1}_{
    \Omega_n
}
- 
X^n_t \,
\mathbbm{1}_{
    \Omega_n
}
\|_{ H_{ \gamma } }
\\ &
=
\sup_{ t \in [0,T] } \,
\Biggl\|
    \Biggl[
    \sum_{ k = 1 }^n
        X^k_t \,
        \mathbbm{1}_{
            \Omega_k
            \setminus
             \Omega_{ k - 1 }
        }
\Biggr]
-
X^n_t \,
\mathbbm{1}_{
    \Omega_n
}
\Biggr\|_{ H_{ \gamma } }
=
\sup_{ t \in [0,T] } \,
\Biggl\|
    \sum_{ k = 1 }^n
    (
     X^k_t 
     - 
     X^n_t
    ) \,
    \mathbbm{1}_{
        \Omega_k
        \setminus
        \Omega_{ k - 1 }
}
\Biggr\|_{ H_{ \gamma } }
\\ & 
=
\sum_{ k = 1 }^n
\bigg[
\mathbbm{1}_{
 \Omega_k
}
\sup_{ t \in [0,T] }
\|
 X^k_t 
 - 
 X^n_t
\|_{ H_{ \gamma } }
\biggr]
\mathbbm{1}_{
    \Omega_k
    \setminus
    \Omega_{ k - 1 }
} .
\end{split}
\end{equation}
This and \eqref{eq:XX_measure_1}
show that
\begin{equation}
\label{eq:X_Y_measure_1}
\P \biggl(
 \forall \, n \in \N 
 \colon
 \mathbbm{1}_{
   \Omega_n 
 }
 \sup_{ t \in [0,T] }
 \| Y_t - X^n_t \|_{ H_{ \gamma } } 
 = 0
\biggr) = 1 .
\end{equation}
Hence, we obtain
for all $ n \in \N $, $ t \in [0,T] $ that
\begin{equation}
\begin{split}
&
[
Y_t \,
\mathbbm{1}_{ \Omega_n }
]_{ \P, \mathscr{B}( H_\gamma ) }
=
[
X^n_t \,
\mathbbm{1}_{ \Omega_n }
]_{ \P, \mathscr{B}( H_\gamma ) }
\\ & =
\biggl(
\biggl[
e^{ t A } \xi_n 
+ \int_0^t
e^{ ( t - s ) A } F( X^n_s ) \ds
\biggr]_{ \P, \mathscr{B}( H_\gamma ) }
+ \int_0^t
e^{ ( t - s ) A } B( X^n_s ) \dWs
\biggr) \,
\mathbbm{1}_{
 \Omega_n
}
\\ & =
\biggl(
\biggl[
e^{ t A } \xi
+
\int_0^t
    e^{ ( t - s ) A } \,
    \mathbbm{1}_{
        \Omega_n
    }
    F( X^n_s )
\ds
\biggr]_{ \P, \mathscr{B}( H_\gamma ) }
+
\int_0^t
    e^{ ( t - s ) A } \,
    \mathbbm{1}_{
        \Omega_n
    }
    B( X^n_s )
\dWs
\biggr) \,
\mathbbm{1}_{
 \Omega_n
}
\\ & =
\biggl(
\biggl[
e^{ t A } \xi
+
\int_0^t
    e^{ ( t - s ) A }
    F( Y_s )
\ds
\biggr]_{ \P, \mathscr{B}( H_\gamma ) }
+
\int_0^t
    e^{ ( t - s ) A }
    B( Y_s )
\dWs
\biggr) \,
\mathbbm{1}_{
 \Omega_n
}
.
\end{split}
\end{equation}
This implies
for all $ t \in [0,T] $ that
\begin{equation}
[
Y_t
]_{ \P, \mathscr{B}( H_\gamma ) }
=
\biggl[
e^{ t A } \xi
+
\int_0^t
    e^{ ( t - s ) A }
    F( Y_s )
\ds
\biggr]_{ \P, \mathscr{B}( H_\gamma ) }
+
\int_0^t
    e^{ ( t - s ) A }
    B( Y_s )
\dWs
.
\end{equation}
Next note that 
\eqref{eq:X_Y_measure_1}
and Proposition~\ref{prop:regularity_SPDE2}
ensure for all $ n \in \N $ that
\begin{equation}
\begin{split}
&
\sup_{ t \in [0,T] }
\left\|
 \max \bigl\{ 1,
   \|  Y_t \, \mathbbm{1}_{ \Omega_n } \|_{ H_{ \gamma } }
 \bigr\}
\right\|_{ \mathscr{L}^p( \P ; \left| \cdot \right| ) }
=
\sup_{ t \in [0,T] }
\left\|
 \max \bigl\{ 1,
   \| X_t^n \, \mathbbm{1}_{ \Omega_n } \|_{ H_{ \gamma } }
 \bigr\}
\right\|_{ \mathscr{L}^p( \P ; \left| \cdot \right| ) }
\\ & 
\leq
\sup_{ t \in [0,T] }
\left\|
 \max \bigl\{ 1,
   \| X_t^n \|_{ H_{ \gamma } }
 \bigr\}
\right\|_{ \mathscr{L}^p( \P ; \left| \cdot \right| ) }
\leq
\sqrt{2}
\left\|
 \max \bigl\{ 1,
   \| \xi_n \|_{ H_{ \gamma } }
 \bigr\}
\right\|_{ \mathscr{L}^p( \P ; \left| \cdot \right| ) }
\\ & \cdot
\mathscr{E}_{ (1 - \eta) } \biggl[ 
\tfrac{ 
  T^{ 1 - \eta }
  \sqrt{ 2 }
}{
  \sqrt{ 1 - \eta }
} \,
\biggl(
  \sup_{ v \in H_{ \gamma } }
  \tfrac{
    \| F( v ) \|_{ H_{ \gamma - \eta } }
  }{
    \max\{ 1, \| v \|_{ H_{ \gamma } } \}
  }
\biggr)
+
\sqrt{
  T^{ 1 - \eta } p ( p - 1 )
} \,
\biggl(
  \sup_{ v \in H_{ \gamma } }
  \tfrac{
    \| B( v ) \|_{ \mathrm{HS}( U, H_{ \gamma - \nicefrac{\eta}{2} } ) }
  }{
    \max\{ 1, \| v \|_{ H_{ \gamma } } \}
  }
\biggr)
\biggr]
.
\end{split}
\end{equation}
This and Fatou's lemma imply
for all $ t \in [0,T] $ that
\begin{equation}
\begin{split}
&
\left\|
 \max \bigl\{ 1,
   \| Y_t \|_{ H_{ \gamma } }
 \bigr\}
\right\|_{ \mathscr{L}^p( \P ; \left| \cdot \right| ) }
=
\Bigl\|
 \liminf_{ n \to \infty }
 \max \bigl\{ 1,
   \| Y_t \, \mathbbm{1}_{ \Omega_n } \|_{ H_{ \gamma } }
 \bigr\}
\Bigr\|_{ \mathscr{L}^p( \P ; \left| \cdot \right| ) }
\\ & \leq
\liminf_{ n \to \infty }
\left\|
 \max \bigl\{ 1,
   \| Y_t \, \mathbbm{1}_{ \Omega_n } \|_{ H_{ \gamma } }
 \bigr\}
\right\|_{ \mathscr{L}^p( \P ; \left| \cdot \right| ) }
\leq
\sqrt{2}
\left\|
 \max \bigl\{ 1,
   \| \xi \|_{ H_{ \gamma } }
 \bigr\}
\right\|_{ \mathscr{L}^p( \P ; \left| \cdot \right| ) }
\\ & \cdot
\mathscr{E}_{ (1 - \eta) } \biggl[ 
\tfrac{ 
  T^{ 1 - \eta }
  \sqrt{ 2 }
}{
  \sqrt{ 1 - \eta }
} \,
\biggl(
  \sup_{ v \in H_{ \gamma } }
  \tfrac{
    \| F( v ) \|_{ H_{ \gamma - \eta } }
  }{
    \max\{ 1, \| v \|_{ H_{ \gamma } } \}
  }
\biggr)
+
\sqrt{
  T^{ 1 - \eta } p ( p - 1 )
} \,
\biggl(
  \sup_{ v \in H_{ \gamma } }
  \tfrac{
    \| B( v ) \|_{ \mathrm{HS}( U, H_{ \gamma - \nicefrac{\eta}{2} } ) }
  }{
    \max\{ 1, \| v \|_{ H_{ \gamma } } \}
  }
\biggr)
\biggr]
.
\end{split}
\end{equation}
The proof of 
Proposition~\ref{prop:existence_continuous_SPDEs}
is thus completed.
\end{proof}
}{}
\iftoggle{arXiv:v3}{%
\subsection{Uniqueness of left-continuous solutions of SEEs 
with semi-globally Lipschitz continuous coefficients}
The proof of the next result, Proposition~\ref{prop:uniqueness_c_local}, is similar to the proof of Theorem~7.4 in Da~Prato \& Zabczyk~\cite{dz92}
(also see, e.g., Lemma~8.2 in van~Neerven, Veraar, \& Weis~\cite{VanNeervenVeraarWeis2008} for an analogous result in a more general framework).
\begin{proposition}[Local solutions]
\label{prop:uniqueness_c_local}
Consider the notation in Subsection~\ref{notation},
let\linebreak
$ ( H, \left< \cdot, \cdot \right>_H, \left\| \cdot \right\|_H ) $
and
$ ( U, \left< \cdot, \cdot \right>_U, \left\| \cdot \right\|_U ) $
be separable $ \R $-Hilbert spaces,
let $ \mathbb{H} \subseteq H $ be a non-empty orthonormal basis of $ H $,
let $ \lambda \colon \mathbb{H} \to \mathbb{R} $
be a function with 
$
  \sup_{ h \in \mathbb{H} }
  \lambda_h < 0
$,
let
$
  A \colon D(A) \subseteq H \rightarrow H
$
be the linear operator which satisfies
$
  D(A) 
  = 
    \bigl\{ 
      v \in H 
      \colon $\linebreak
      $
      \sum_{ h \in \mathbb{H} } 
	\left| 
	  \lambda_h 
	  \langle h, v \rangle_H 
	\right|^2
	< \infty
    \bigr\}
$
and
$
\forall \, v \in D(A) \colon
    Av
   =
    \sum_{ h \in \mathbb{H} } 
    \lambda_h \langle h, v \rangle_H h
$,
let 
$ 
( 
  H_r 
  ,
  \left< \cdot, \cdot \right>_{ H_r }
  ,
  \left\| \cdot \right\|_{ H_r } 
) 
$, $ r \in \R $,
be a family of interpolation spaces associated 
to $ - A $,
let
$ T \in (0,\infty) $,
$ \gamma \in \R $, $ \eta \in [0,1) $,
$ F \in C( H_{ \gamma }, H_{ \gamma - \eta } ) $,
$ B \in C( H_{ \gamma }, \mathrm{HS}( U, H_{ \gamma - \nicefrac{\eta}{2} } ) ) $
satisfy for all bounded sets 
$ E \subseteq H_{ \gamma } $
that
$
| F|_E |_{
  \calC^1( E, \left\| \cdot \right\|_{ H_{ \gamma - \eta } } )
}
+
| B|_E |_{
  \calC^1( E, \left\| \cdot \right\|_{ \mathrm{HS}( U, H_{ \gamma - \nicefrac{\eta}{2} } ) } )
}
< \infty
$,
let 
$ ( \Omega, \mathscr{F}, \P ) $
be a probability space with a normal 
filtration $ ( \mathscr{F}_t )_{ t \in [0,T] } $,
let 
$
  ( W_t )_{ t \in [0,T] } 
$ 
be an $ \operatorname{Id}_U $-cylindrical
$ ( \Omega, \calF, \P, ( \calF_t )_{ t \in [0,T] } ) $-Wiener
process,
let 
$
\tau_k \colon \Omega \to [0,T]
$,
$ k \in \{ 1, 2 \} $,
be $ ( \mathscr{F}_t )_{ t \in [0,T] } $-stopping times,
and let 
$
X^k \colon [0,T] \times \Omega \to H_{ \gamma }
$,
$ k \in \{ 1, 2 \} $,
be $ ( \mathscr{F}_t )_{ t \in [0,T] } $/$ \mathscr{B}( H_{ \gamma } ) $-adapted
stochastic processes 
with left-continuous and bounded sample paths
which satisfy
for all $ k \in \{ 1, 2 \} $,
$ t \in [0,T] $ 
that
\begin{equation}
\begin{split}
\bigl[
  X^k_t \,
  \mathbbm{1}_{
    \{ t \leq \tau_k \}
  }
\bigr]_{
  \P, \mathscr{B}( H_{ \gamma } )
}
=
\biggl(
& \biggl[
e^{ t A } X^k_0
  + 
  \int_0^t 
  \mathbbm{1}_{
    \{ s < \tau_k \}
  }
  \,
  e^{ ( t - s ) A } F( X^k_s ) \ds
\biggr]_{ \P, \mathscr{B}( H_{ \gamma } ) }
\\ & + 
  \int_0^t 
  \mathbbm{1}_{
    \{ s < \tau_k \}
  }
  \,
  e^{ ( t - s ) A } 
  B( X^k_s ) \dWs
\biggr)
\mathbbm{1}_{
  \{ t \leq \tau_k \}
}
.
\end{split}
\end{equation}
Then
$
\P \bigl(
  \forall \, t \in [0,T] \colon 
  \mathbbm{1}_{
    \{ X^1_0 = X^2_0 \}
  } \,
  X^1_{ \min\{ t, \tau_1, \tau_2 \} }
  =
  \mathbbm{1}_{
    \{ X^1_0 = X^2_0 \}
  } \,
  X^2_{ \min\{ t, \tau_1, \tau_2 \} }
\bigr)
= 1
$.
\end{proposition}
Corollary~\ref{cor:uniqueness_c}
is an immediate consequence from 
Proposition~\ref{prop:uniqueness_c_local}.
\begin{corollary}[Continuous solutions]
\label{cor:uniqueness_c}
Consider the notation in Subsection~\ref{notation},
let\linebreak
$ ( H, \left< \cdot , \cdot \right>_H , \left\| \cdot \right\|_H ) $
and
$ ( U, \left< \cdot , \cdot \right>_U , \left\| \cdot \right\|_U ) $
be separable $ \R $-Hilbert spaces,
let $ \mathbb{H} \subseteq H $ be a non-empty orthonormal basis of $ H $,
let $ \lambda \colon \mathbb{H} \to \mathbb{R} $
be a function with 
$
  \sup_{ h \in \mathbb{H} }
  \lambda_h < 0
$,
let
$
  A \colon D(A) \subseteq H \rightarrow H
$
be the linear operator which satisfies
$
  D(A) 
  = 
    \bigl\{ 
      v \in H 
      \colon $\linebreak
      $
      \sum_{ h \in \mathbb{H} } 
	\left| 
	  \lambda_h 
	  \langle h, v \rangle_H 
	\right|^2
	< \infty
    \bigr\}
$
and
$
\forall \, v \in D(A) \colon
    Av
   =
    \sum_{ h \in \mathbb{H} } 
    \lambda_h \langle h, v \rangle_H h
$,
let 
$ 
( 
  H_r 
  ,
  \left< \cdot, \cdot \right>_{ H_r }
  ,
  \left\| \cdot \right\|_{ H_r } 
) 
$, $ r \in \R $,
be a family of interpolation spaces associated 
to $ - A $,
let
$ T \in (0,\infty) $,
$ \gamma \in \R $, $ \eta \in [0,1) $,
$ F \in C( H_{ \gamma }, H_{ \gamma - \eta } ) $,
$ B \in C( H_{ \gamma }, \mathrm{HS}( U, H_{ \gamma - \nicefrac{\eta}{2} } ) ) $
satisfy for all bounded sets 
$ E \subseteq H_{ \gamma } $
that
$
| F|_E |_{
  \calC^1( E, \left\| \cdot \right\|_{ H_{ \gamma - \eta } } )
}
+
| B|_E |_{
  \calC^1( E, \left\| \cdot \right\|_{ \mathrm{HS}( U, H_{ \gamma - \nicefrac{\eta}{2} } ) } )
}
< \infty
$,
let 
$ ( \Omega, \mathscr{F}, \P ) $
be a probability space with a normal 
filtration $ ( \mathscr{F}_t )_{ t \in [0,T] } $,
let 
$
  ( W_t )_{ t \in [0,T] } 
$ 
be an $ \operatorname{Id}_U $-cylindrical
$ ( \Omega, \calF, \P, ( \calF_t )_{ t \in [0,T] } ) $-Wiener
process,
and let 
$
X^k \colon [0,T] \times \Omega \to H_{ \gamma }
$,
$ k \in \{ 1, 2 \} $,
be\linebreak
$ ( \mathscr{F}_t )_{ t \in [0,T] } $/$ \mathscr{B}( H_{ \gamma } ) $-adapted 
stochastic processes with continuous sample paths
which satisfy
for all $ k \in \{ 1, 2 \} $,
$ t \in [0,T] $ that
\begin{equation}
[ X^k_t ]_{
  \P, \mathscr{B}( H_{ \gamma } )
}
=
\biggl[
e^{ t A } X^1_0 
+ \int_0^t e^{ ( t - s ) A } F( X^k_s ) \ds
\biggr]_{
  \P, \mathscr{B}( H_{ \gamma } )
}
+ \int_0^t e^{ ( t - s ) A } B( X^k_s ) \dWs
.
\end{equation}
Then 
$
\P\bigl(
  \forall \, t \in [0,T] \colon
  X_t^1 = X_t^2
\bigr)
= 1
$.
\end{corollary}
}{}
\section{Convergence in H\"{o}lder norms for Galerkin approximations}
\label{sec:Galerkin}

\subsection{Setting}
\label{ssec:setting}
%
%
%
%
Consider the notation in Subsection~\ref{notation},
let 
$
  ( 
    H,
    \langle \cdot, \cdot \rangle_H, 
    \left\| \cdot \right\|_H
  )
$
and
$
  ( 
    U ,
    \langle \cdot, \cdot \rangle_U , 
    \left\| \cdot \right\|_U
  )
$
be separable $ \mathbb{R} $-Hilbert spaces, 
let $ \mathbb{H} \subseteq H $ be a non-empty orthonormal basis of $ H $,
let 
$ T, \iota \in (0,\infty) $,
let 
$ ( \Omega, \calF, \P ) $
be a probability space with a normal filtration 
$ ( \calF_t )_{ t \in [0,T] } $,
let 
$
  ( W_t )_{ t \in [0,T] } 
$ 
be an $ \operatorname{Id}_U $-cylindrical
$ ( \Omega, \calF, \P, ( \calF_t )_{ t \in [0,T] } ) $-Wiener
process,
let $ \lambda \colon \mathbb{H} \to \mathbb{R} $
be a function with 
$
  \sup_{ h \in \mathbb{H} }
  \lambda_h < 0
$,
let
$
  A \colon D(A) \subseteq H \rightarrow H
$ 
be the linear operator which satisfies
\begin{equation}
  D(A) 
  = 
    \Biggl\{ 
      v \in H 
      \colon
      \sum_{ h \in \mathbb{H} } 
	\left| 
	  \lambda_h 
	  \langle h, v \rangle_H 
	\right|^2
	< \infty
    \Biggr\}
\end{equation}
and which satisfies for all $ v \in D(A) $ that
\begin{equation}
    Av
   =
    \sum_{ h \in \mathbb{H} } 
    \lambda_h \langle h, v \rangle_H h,
\end{equation}
let 
$ 
  ( H_r , \left< \cdot , \cdot \right>_{ H_r }, \left\| \cdot \right\|_{ H_r } ) 
$,
$ r \in \R $,
be a family of interpolation spaces associated to
$ - A $%
\iftoggle{arXiv:v3}{}{
(cf., e.g., \cite[Section~3.7]{sy02})},
let 
$ \gamma \in \R $, $ \alpha \in [0,1) $, 
$ \beta \in [ 0, \nicefrac{ 1 }{ 2 } ) $, 
$ \chi \in [ \beta, \nicefrac{ 1 }{ 2 } ) $,
$ F \in C( H_{ \gamma }, H_{ \gamma - \alpha } ) $,
$ B \in C( H_{ \gamma }, \mathrm{HS}( U, H_{ \gamma - \beta } ) ) $
satisfy for all bounded sets
$ E \subseteq H_{ \gamma } $
that
\begin{equation}
  | F|_E |_{
    \calC^1( E, \left\| \cdot \right\|_{  H_{ \gamma - \alpha }  } )
  }
  +
  | B|_E |_{
    \calC^1( E, \left\| \cdot \right\|_{  \mathrm{HS}( U, H_{ \gamma - \beta } )  } )
  }
  < \infty,
\end{equation}
let $ \mathbb{H}_N \subseteq \mathbb{H} $,
$ N \in \N_0 $,
be sets
which satisfy 
$
  \mathbb{H}_0 = \mathbb{H}
$
and
$
  \sup_{ N \in \mathbb{N} } 
  N^{ \iota } 
  \sup\!\big( 
    \{
      \nicefrac{ 1 }{ | \lambda_h | }
      \colon
      h \in \mathbb{H} \backslash \mathbb{H}_N 
    \} \cup \{ 0 \}
  \big)
  < \infty
$,
let 
$ P_N \in L( H_{ \min\{ 0, \gamma - 1 \} } ) $,
$ N \in \N_0 $,
and 
$ \mathscr{P}_N \in L( U ) $,
$ N \in \N_0 $,
be linear operators which satisfy 
for all $ N \in \N_0 $,
$ v \in H $
that
\begin{equation}
  P_N( v ) 
  = \sum_{ h \in \mathbb{H}_N } 
  \left< h, v \right>_H h,
\end{equation}
and let 
$
  X^N \colon [0,T] \times \Omega \rightarrow H_{ \gamma }
$,
$ N \in \N_0 $,
be
$
( \mathscr{F}_t )_{ t \in [0,T] }
$/$ \mathscr{B}( H_{ \gamma } ) $-adapted stochastic processes
with continuous sample paths which satisfy
for all $ N \in \N_0 $, $ t \in [0,T] $ 
that
\begin{equation}
\label{eq:solution}
\bigl[  X_t^N  \bigr]_{ \P, \mathscr{B}( H_\gamma ) }
  = 
\biggl[  e^{ t A } P_N X_0^0
  +
  \int_0^t
    e^{ ( t - s ) A }
    P_N F( X^N_s )
  \ds \biggr]_{ \P, \mathscr{B}( H_\gamma ) }
  +
  \int_0^t
    e^{ ( t - s ) A }
    P_N B( X_s^N ) \mathscr{P}_N
  \dWs
  .
\end{equation}
\subsection{Strong convergence in H\"{o}lder norms for Galerkin approximations of SEEs with globally Lipschitz continuous nonlinearities}

The next lemma, Lemma~\ref{lem:a_priori_XN} below,
follows directly from, e.g.,
\iftoggle{arXiv:v3}{%
Proposition~\ref{prop:existence_continuous_SPDEs}
and, e.g.,
Corollary~\ref{cor:uniqueness_c}.%
}{%
Proposition~3.6
in~\cite{CoxHutzenthalerJentzenVanNeervenWelti2016arXiv}
and, e.g.,
Corollary~3.8
in~\cite{CoxHutzenthalerJentzenVanNeervenWelti2016arXiv}.%
}

\begin{lemma}
\label{lem:a_priori_XN}
Assume the setting in 
Subsection~\ref{ssec:setting},
let 
$ p \in [2,\infty) $,
$ 
  \eta \in [ \max\{ \alpha, 2 \beta \}, 1 )
$,
$ N \in \N_0 $,
and assume that
\begin{equation}
  \E\big[ 
    \| X_0^0 \|_{ H_{ \gamma } }^p
  \big]
  +
  | F |_{
    \calC^1( H_{ \gamma }, \left\| \cdot \right\|_{  H_{ \gamma - \alpha }  } )
  }
  +
  | B |_{
    \calC^1( H_{ \gamma }, \left\| \cdot \right\|_{  \mathrm{HS}( U, H_{ \gamma - \beta } )  } )
  }
  < \infty.
\end{equation}
Then
\begin{align}
&
  \sup_{ t \in [0,T] }
  \left\|
    \max\{ 1, \| X_t^N \|_{ H_{ \gamma } } \}
  \right\|_{
    \mathscr{L}^p( \P;  \left| \cdot \right| )
  }
\leq 
  \sqrt{ 2 }
  \left\|
    \max\{ 1, \| X^0_0 \|_{ H_{ \gamma } } \}
  \right\|_{
    \mathscr{L}^p( \P;  \left| \cdot \right| )
  }
\\ \nonumber & 
  \cdot
  \mathscr{E}_{ ( 1 - \eta ) }\biggl[ 
    \tfrac{ 
      T^{ 1 - \eta }
      \sqrt{ 2 }
    }{
      \sqrt{ 1 - \eta }
    }
    \biggl(
      \sup_{ v \in H_{ \gamma } }
      \tfrac{
        \| F( v ) \|_{ H_{ \gamma - \eta } }
      }{
        \max\{ 1, \| v \|_{ H_{ \gamma } } \}
      }
    \biggr)
    +
    \sqrt{
      T^{ 1 - \eta } p ( p - 1 )
    }
    \biggl(
      \sup_{ v \in H_{ \gamma } }
      \tfrac{
        \| B( v ) \mathscr{P}_N \|_{ \mathrm{HS}( U, H_{ \gamma - \nicefrac{\eta}{2} } ) }
      }{
        \max\{ 1, \| v \|_{ H_{ \gamma } } \}
      }
    \biggr)
  \biggr]
  .
\end{align}
\end{lemma}

\begin{lemma}
\label{lem:spectral_noise}
Assume the setting in 
Subsection~\ref{ssec:setting},
let 
$ p \in [2,\infty) $,
$ 
  \eta \in [ \max\{ \alpha, 2 \beta \}, 1 )
$,
$ N \in \N_0 $,
and assume that
\begin{equation}
  \E\big[ 
    \| X_0^0 \|_{ H_{ \gamma } }^p
  \big]
  +
  | F |_{
    \calC^1( H_{ \gamma }, \left\| \cdot \right\|_{  H_{ \gamma - \alpha }  } )
  }
  +
  | B |_{
    \calC^1( H_{ \gamma }, \left\| \cdot \right\|_{  \mathrm{HS}( U, H_{ \gamma - \beta } )  } )
  }
  < \infty.
\end{equation}
Then
\begin{equation}
\begin{split}
&
  \sup_{ t \in [0,T] }
  \left\|
    X^0_t - X^N_t
  \right\|_{ 
    \mathscr{L}^p( 
      \P ; 
      \left\| \cdot \right\|_{  H_{ \gamma }  }
    ) 
  }
  \leq 
  \Biggl[
    \sqrt{ 2 }
    \sup_{ t \in [0,T] }
    \left\| ( P_0 - P_N ) X_t^0 \right\|_{
      \mathscr{L}^p( \P ; \left\| \cdot \right\|_{  H_{ \gamma }  } )
    }
\\ & 
    +
    \tfrac{
      T^{ \nicefrac{1}{2} - \chi }
      \sqrt{ p \, ( p - 1 ) } 
    }{
      \sqrt{ 1 - 2 \chi }
    }
      \biggl(
        1 +
        \sup_{ t \in [0,T] }
        \| X^N_t \|_{
          \mathscr{L}^p( \P ; \left\| \cdot \right\|_{  H_{ \gamma }  } )
        }
      \biggr)
      \Biggl(
        \sup_{ v \in H_{ \gamma } }
        \frac{
          \| B( v ) ( \mathscr{P}_0 - \mathscr{P}_N ) \|_{
            \mathrm{HS}( U, H_{ \gamma - \chi } )
          }
        }{
          1 + \| v \|_{ H_{ \gamma } }
        }
      \Biggr)
  \Biggr]
\\ & 
  \cdot 
  \mathscr{E}_{ (1 - \eta) } \biggl[ 
  \tfrac{ 
    T^{ 1 - \eta }
    \sqrt{ 2 }
    \,
    \left| F \right|_{
      \calC^1(
        H_{ \gamma }, \left\| \cdot \right\|_{  H_{ \gamma - \eta }  } 
      )
    }
  }{ \sqrt{ 1 - \eta } }
    +
    \sqrt{ 
      T^{ 1 - \eta }
      p ( p - 1 ) 
    } 
  \left| B \right|_{
    \calC^1(
      H_{ \gamma }, \left\| \cdot \right\|_{  \mathrm{HS}( U, H_{ \gamma - \nicefrac{\eta}{2} } )  } 
    )
  }
  \| \mathscr{P}_0 \|_{ L( U ) }
  \biggr]
  < \infty
  .
\end{split}
\end{equation}
\end{lemma}

\begin{proof}[Proof 
of Lemma~\ref{lem:spectral_noise}]
First of all, observe that
Lemma~\ref{lem:a_priori_XN}
ensures that
\begin{equation}
\label{eq:finiteness_XN}
  \sup_{ t \in [0,T] }
  \max \bigl\{
    \| X^0_t \|_{ \mathscr{L}^p( \P; \left\| \cdot \right\|_{ H_{ \gamma } } ) },
    \| X^N_t \|_{ \mathscr{L}^p( \P; \left\| \cdot \right\|_{ H_{ \gamma } } ) }
  \bigr\}
  < \infty
  .
\end{equation}
We can hence apply
\iftoggle{arXiv:v3}{%
Proposition~\ref{prop:strong_perturbation}%
}{%
Proposition~3.5
in~\cite{CoxHutzenthalerJentzenVanNeervenWelti2016arXiv}%
}
to obtain that
\begin{equation}
\begin{split}
&
  \sup_{ t \in [0,T] }
  \left\|
    X^0_t - X^N_t
  \right\|_{ 
    \mathscr{L}^p( 
      \P ; 
      \left\| \cdot \right\|_{  H_{ \gamma }  }
    ) 
  }
\\ & \leq
  \mathscr{E}_{ (1 - \eta) } \biggl[ 
  \tfrac{ 
    T^{ 1 - \eta }
    \sqrt{ 2 }
    \,
    \left| P_N F( \cdot ) \right|_{
      \calC^1(
        H_{ \gamma }, \left\| \cdot \right\|_{  H_{ \gamma - \eta }  } 
      )
    }
  }{ \sqrt{ 1 - \eta } }
    +
    \sqrt{ 
      T^{ 1 - \eta }
      p ( p - 1 ) 
    } 
  \left| P_N B( \cdot ) \mathscr{P}_0 \right|_{
    \calC^1(
      H_{ \gamma }, \left\| \cdot \right\|_{  \mathrm{HS}( U, H_{ \gamma - \nicefrac{\eta}{2} } )  } 
    )
  }
  \biggr]
\\ & \quad 
  \cdot
  \sqrt{2} 
  \,
  \sup_{ t \in [0,T] }
  \bigg\|
  \biggl[
    X_t^0
    -
      \int_0^t
      e^{ ( t - s ) A } P_N F( X^0_s ) 
      \ds
      \biggr]_{ \P, \mathscr{B}( H_\gamma ) }
      -
      \int_0^t
      e^{ ( t - s ) A } P_N B( X^0_s ) \mathscr{P}_0
      \dWs
\\ & \quad
  +
  \biggl[
      \int_0^t
      e^{ ( t - s ) A } P_N F( X^N_s ) 
      \ds
      - X^N_t
      \biggr]_{ \P, \mathscr{B}( H_\gamma ) }
      +
      \int_0^t
      e^{ ( t - s ) A } P_N B( X^N_s ) \mathscr{P}_0
      \dWs
  \bigg\|_{ 
    L^p( 
      \P ; 
      \left\| \cdot \right\|_{  H_{ \gamma }  } 
    ) 
  }
  .
\end{split}
\end{equation}
This shows that
\begin{equation}
\begin{split}
&
  \sup_{ t \in [0,T] }
  \left\|
    X^0_t - X^N_t
  \right\|_{ 
    \mathscr{L}^p( 
      \P ; 
      \left\| \cdot \right\|_{  H_{ \gamma }  } 
    ) 
  }
\\ & \leq
  \mathscr{E}_{ (1 - \eta) } \biggl[ 
  \tfrac{ 
    T^{ 1 - \eta }
    \sqrt{ 2 }
    \,
    \left| F \right|_{
      \calC^1(
        H_{ \gamma }, \left\| \cdot \right\|_{  H_{ \gamma - \eta }  } 
      )
    }
  }{ \sqrt{ 1 - \eta } }
    +
    \sqrt{ 
      T^{ 1 - \eta }
      p ( p - 1 ) 
    } 
  \left| B \right|_{
    \calC^1(
      H_{ \gamma }, \left\| \cdot \right\|_{  \mathrm{HS}( U, H_{ \gamma - \nicefrac{\eta}{2} } )  } 
    )
  }
  \| \mathscr{P}_0 \|_{ L( U ) }
  \biggr]
\\ & \quad 
  \cdot
  \sqrt{2} 
  \,
  \sup_{ t \in [0,T] }
  \left\|
    \bigl[ ( P_0 - P_N ) X^0_t \bigr]_{ \P, \mathscr{B}( H_\gamma ) }
    +
      \int_0^t
      e^{ ( t - s ) A } P_N B( X^N_s ) ( \mathscr{P}_0 - \mathscr{P}_N )
      \dWs
  \right\|_{ 
    L^p( 
      \P ; 
      \left\| \cdot \right\|_{  H_{ \gamma }  } 
    ) 
  }
  .
\end{split}
\end{equation}
The Burkholder-Davis-Gundy-type inequality
in Lemma~7.7 in Da Prato \& Zabczyk~\cite{dz92}
hence implies that
\begin{align}
& \nonumber
  \sup_{ t \in [0,T] }
  \left\|
    X^0_t - X^N_t
  \right\|_{ 
    \mathscr{L}^p( 
      \P ; 
      \left\| \cdot \right\|_{ 
        H_{ \gamma } 
      } 
    ) 
  }
\\ \nonumber & \leq
  \mathscr{E}_{ (1 - \eta) } \biggl[ 
  \tfrac{ 
    T^{ 1 - \eta }
    \sqrt{ 2 }
    \,
    \left| F \right|_{
      \calC^1(
        H_{ \gamma }, \left\| \cdot \right\|_{  H_{ \gamma - \eta }  } 
      )
    }
  }{ \sqrt{ 1 - \eta } }
    +
    \sqrt{ 
      T^{ 1 - \eta }
      p ( p - 1 ) 
    } 
  \left| B \right|_{
    \calC^1(
      H_{ \gamma }, \left\| \cdot \right\|_{  \mathrm{HS}( U, H_{ \gamma - \nicefrac{\eta}{2} } )  } 
    )
  }
  \| \mathscr{P}_0 \|_{ L( U ) }
  \biggr]
\\ & \quad
  \cdot
  \sqrt{2} 
  \,
  \Biggl[
    \sup_{ t \in [0,T] }
    \left\| ( P_0 - P_N ) X_t^0 \right\|_{
      \mathscr{L}^p( \P ; \left\| \cdot \right\|_{  H_{ \gamma }  } )
    }
\\ \nonumber & \quad
    +
      \sup_{ s \in [0,T] }
      \left\|
        B( X^N_s ) [ \mathscr{P}_0 - \mathscr{P}_N ]
      \right\|_{
        \mathscr{L}^p( \P ; \left\| \cdot \right\|_{ \mathrm{HS}( U, H_{ \gamma - \chi } ) } )
      }
    \sqrt{
      \tfrac{ p \, ( p - 1 ) }{ 2 }
      \sup_{ t \in [0,T] }
      \int_0^t
      \left( t - s \right)^{ - 2 \chi }
      \ds
    }
  \Biggr]
  .
\end{align}
This and \eqref{eq:finiteness_XN} complete the proof
of Lemma~\ref{lem:spectral_noise}.
\end{proof}

\begin{corollary}
\label{cor:convGlobLip_0}
Assume the setting in Subsection~\ref{ssec:setting},
let 
$
  \vartheta \in [ 0, \min\{ 1 - \alpha, \nicefrac{ 1 }{ 2 } - \beta \} )
$,
$ p \in [2,\infty) $,
and assume that
$ X^0_0( \Omega ) \subseteq H_{ \gamma + \vartheta } $
and
\begin{align}
\label{eq:ass_FB0}
  \E\big[ 
    \| X_0^0 \|_{ H_{ \gamma + \vartheta } }^p
  \big]
+
  | F |_{
    \calC^1( H_{ \gamma }, \left\| \cdot \right\|_{  H_{ \gamma - \alpha }  } )
  }
+
  | B |_{
    \calC^1( H_{ \gamma }, \left\| \cdot \right\|_{  \mathrm{HS}( U, H_{ \gamma - \beta } )  } )
  }  & < \infty,
\\
\label{eq:ass_B1.general_0}
  \sup_{ N \in \N }
  \sup_{ v \in H_{ \gamma } }
  \left[
  \frac{
  N^{ \iota \vartheta }
  \,
  \|
    B( v ) ( \mathscr{P}_0 - \mathscr{P}_N )
  \|_{ \mathrm{HS}( U, H_{ \gamma - \chi } ) }
  }{
    1 + \| v \|_{ H_{ \gamma } }
  }
  \right]
  & < \infty
  .
\end{align}
Then it holds that
\begin{equation}
  \sup_{ N \in \N_0 }
  \sup_{ t \in [0,T] }
  \bigl(
  \| F( X^N_t ) \|_{
    \mathscr{L}^p( 
      \P ;
      \left\| \cdot \right\|_{  H_{ \gamma - \alpha }  }
    )
  }
  +
  \| B( X^N_t ) \mathscr{P}_N \|_{
    \mathscr{L}^p( 
      \P ;
      \left\| \cdot \right\|_{  \mathrm{HS}( U, H_{ \gamma - \chi } )  }
    )
  } \bigr) < \infty
\end{equation}
and
\begin{equation}
  \sup_{ N \in \N_0 }
  \sup_{ t \in [0,T] }
  \bigl( N^{ \iota \vartheta }
  \,
  \|
    X^0_t - X^N_t
  \|_{
    \mathscr{L}^p( 
      \P ; 
      \left\| \cdot \right\|_{  H_{ \gamma }  } 
    )
  }
  \bigr)
  <
  \infty.
\end{equation}
\end{corollary}

\begin{proof}[Proof of
Corollary~\ref{cor:convGlobLip_0}]
Combining the assumptions that 
$ X^0_0( \Omega ) \subseteq H_{ \gamma + \vartheta } $
and
$
  \E\big[ 
    \| X^0_0 \|_{ H_{ \gamma + \vartheta } }^p
  \big] $\linebreak
  $
  < \infty
$
with, e.g.,
\iftoggle{arXiv:v3}{%
Proposition~\ref{prop:existence_continuous_SPDEs}
and, e.g., Corollary~\ref{cor:uniqueness_c}%
}{%
Proposition~3.6
in~\cite{CoxHutzenthalerJentzenVanNeervenWelti2016arXiv}
and, e.g.,
Corollary~3.8
in~\cite{CoxHutzenthalerJentzenVanNeervenWelti2016arXiv}%
}
ensures that
$
  \forall \, t \in [0,T] \colon
  \P\big( 
    X^0_t \in H_{ \gamma + \vartheta } 
  \big) = 1
$
and 
$
  \sup_{ t \in [0,T] }
  \E\big[ \|  \mathbbm{1}_{ \{ X^0_t \in H_{ \gamma + \vartheta } \} } X_t^0 \|_{ H_{ \gamma + \vartheta } }^p \big] < \infty
$.
This and the assumption that 
$
  \sup_{ N \in \mathbb{N} } 
  N^{ \iota } 
  \sup\!\big( 
    \{
      \nicefrac{ 1 }{ | \lambda_h | }
      \colon
      h \in \mathbb{H} \backslash \mathbb{H}_N 
    \} \cup \{ 0 \}
  \big)
  < \infty
$
imply that
\begin{align}
& \nonumber
  \sup_{ N \in \N }
  \sup_{ t \in [0,T] }
  \left[
    N^{ \iota \vartheta }
    \left\|
      ( P_0 - P_N ) X^0_t
    \right\|_{
      \mathscr{L}^p( \P ; \left\| \cdot \right\|_{  H_{ \gamma }  } )
    }
  \right]
\\ \nonumber & \leq 
  \sup_{ N \in \N }
  \sup_{ t \in [0,T] }
  \left[
    N^{ \iota \vartheta }
    \left\|
      ( - A )^{ - \vartheta }
      ( P_0 |_{ H_\gamma } - P_N |_{ H_\gamma } )
    \right\|_{
      L( H_{ \gamma } )
    }
    \big\|
      \mathbbm{1}_{ \{ X^0_t \in H_{ \gamma + \vartheta } \} } X^0_t
    \big\|_{
      \mathscr{L}^p( \P ; \left\| \cdot \right\|_{ H_{ \gamma + \vartheta } } )
    }
  \right]
\\ \label{eq:N_iota_rate} & \leq 
  \left[
    \sup_{ N \in \N }
    N^{ \iota \vartheta }
    \left\|
      ( - A )^{ - 1 }
      ( \operatorname{Id}_{ H_{ \gamma } } - P_N |_{ H_\gamma } )
    \right\|_{
      L( H_{ \gamma } )
    }^{ \vartheta }
  \right]
  \left[
    \sup_{ t \in [0,T] }
    \big\|
      \mathbbm{1}_{ \{ X^0_t \in H_{ \gamma + \vartheta } \} } X^0_t
    \big\|_{
      \mathscr{L}^p( \P ; \left\| \cdot \right\|_{ H_{ \gamma + \vartheta } } )
    }
  \right]
\\ \nonumber & = 
  \left[
    \sup_{ N \in \N }
    N^{ \iota \vartheta }
    \left[
      \sup\!\big( 
        \{ 
          \nicefrac{ 1 }{ | \lambda_h | }
          \colon
          h \in \mathbb{H} \backslash \mathbb{H}_N
        \} \cup \{ 0 \}
      \big)
    \right]^{ \vartheta }
  \right]
  \left[
    \sup_{ t \in [0,T] }
    \big\|
      \mathbbm{1}_{ \{ X^0_t \in H_{ \gamma + \vartheta } \} } X^0_t
    \big\|_{
      \mathscr{L}^p( \P ; \left\| \cdot \right\|_{ H_{ \gamma + \vartheta } } )
    }
  \right]
\\ \nonumber &   < \infty .
\end{align}
In addition, observe that 
\eqref{eq:ass_FB0},
\eqref{eq:ass_B1.general_0},
and Lemma~\ref{lem:a_priori_XN}
ensure that
\begin{equation}
\label{eq:X_bound}
  \sup_{ N \in \N_0 }
  \sup_{ t \in [0,T] }
  \| X_t^N \|_{
    \mathscr{L}^p( \P ; \left\| \cdot \right\|_{  H_{ \gamma }  } )
  }
  < \infty 
  .
\end{equation}
The triangle inequality and again \eqref{eq:ass_B1.general_0} hence prove that
\begin{equation}
\label{eq:X_B_bound}
\begin{split}
&
  \sup_{ N \in \N_0 }
    \sup_{ t \in [0,T] }
  \| B( X^N_t ) \mathscr{P}_N \|_{
    \mathscr{L}^p( 
      \P ;
      \left\| \cdot \right\|_{  \mathrm{HS}( U, H_{ \gamma - \chi } )  }
    )
  }
\\ & 
\leq 
  \biggl(
    1 +
    \sup_{ N \in \N_0 }
    \sup_{ t \in [0,T] }
    \| X^N_t \|_{
      \mathscr{L}^p( \P ; \left\| \cdot \right\|_{  H_{ \gamma }  } )
    }
  \biggr)
  \biggl(
  \sup_{ N \in \N_0 }
  \sup_{ v \in H_{ \gamma } }
  \tfrac{
    \| 
      B( v ) \mathscr{P}_N
    \|_{ \mathrm{HS}( U, H_{ \gamma - \chi } )
    }
  }{
    1 + \| v \|_{ H_{ \gamma } }
  }
  \biggr)
\\ & \leq
  \biggl(
    1 +
    \sup_{ N \in \N_0 }
    \sup_{ t \in [0,T] }
    \| X^N_t \|_{
      \mathscr{L}^p( \P ; \left\| \cdot \right\|_{  H_{ \gamma }  } )
    }
  \biggr)
  \\ & \quad \cdot
  \biggl(
    \sup_{ v \in H_{ \gamma } }
    \tfrac{
      \| 
        B( v )
      \|_{\mathrm{HS}( U, H_{ \gamma - \chi } )}
      \| \mathscr{P}_0 \|_{ L( U ) }
    }{
      1 + \| v \|_{ H_{ \gamma } }
    }
    +
    \sup_{ N \in \N_0 }
    \sup_{ v \in H_{ \gamma } }
    \tfrac{
      \| 
        B( v ) ( \mathscr{P}_0 - \mathscr{P}_N )
      \|_{\mathrm{HS}( U, H_{ \gamma - \chi } )}
    }{
      1 + \| v \|_{ H_{ \gamma } }
    }
  \biggr)
  < \infty
  .
\end{split}
\end{equation}
In the next step we combine \eqref{eq:X_bound}, 
\eqref{eq:N_iota_rate}, and \eqref{eq:ass_B1.general_0}
with Lemma~\ref{lem:spectral_noise} 
to obtain that
\begin{equation}
\label{eq:X_rate_zero}
  \sup\nolimits_{ N \in \N_0 }
  \sup\nolimits_{ t \in [0,T] }
  \big(
  N^{ \iota \vartheta }
  \,
  \|
    X^0_t - X^N_t
  \|_{
    \mathscr{L}^p( 
      \P ; 
      \left\| \cdot \right\|_{  H_{ \gamma }  } 
    )
  }
  \big)
  <
  \infty
  .
\end{equation}
Furthermore, observe that
\eqref{eq:X_bound}
assures that
$
  \sup_{ N \in \N_0 }
      \sup_{ t \in [0,T] }
  \| F( X^N_t ) \|_{
    \mathscr{L}^p( \P ; \left\| \cdot \right\|_{  H_{ \gamma - \alpha }  } )
  }
  < \infty
$.
This, \eqref{eq:X_B_bound},
and \eqref{eq:X_rate_zero}
complete
the proof of Corollary~\ref{cor:convGlobLip_0}.
\end{proof}
The next result, Corollary~\ref{cor:convGlobLip}, proves
strong convergence rates in
H\"{o}lder norms for spatial spectral Galerkin approximations of
SEEs
with globally Lipschitz continuous nonlinearities.
Note
in the setting of Corollary~\ref{cor:convGlobLip} that,
e.g.,
Becker et al.~\cite[Theorem~1.1 and Lemma~2.6]{BeckerGessJentzenKloeden2018arXiv}
show
in the case $ \iota = 2 $, $ \delta = 0 $
that the convergence rate established in~\eqref{eq:convGlobLip}
is essentially sharp
(cf., e.g., Conus, Jentzen, \& Kurniawan~\cite[Lemma~7.2]{ConusJentzenKurniawan2019}).
\begin{corollary}
\label{cor:convGlobLip}
Assume the setting in Subsection~\ref{ssec:setting},
let 
$
  \vartheta \in ( 0, \min\{ 1 - \alpha, \nicefrac{ 1 }{ 2 } - \beta \} )
$,
$ p \in ( \nicefrac{ 1 }{ \vartheta }, \infty ) $,
and assume that
$ X^0_0( \Omega ) \subseteq H_{ \gamma + \vartheta } $,
$ 
  \E\big[ 
    \| X_0^0 \|_{ H_{ \gamma + \vartheta } }^p
  \big] < \infty $,
$
  | F |_{
    \calC^1( H_{ \gamma }, \left\| \cdot \right\|_{  H_{ \gamma - \alpha }  } )
  } < \infty $,
$
  | B |_{
    \calC^1( H_{ \gamma }, \left\| \cdot \right\|_{  \mathrm{HS}( U, H_{ \gamma - \beta } )  } )
  } < \infty $,
and
\begin{equation}
\label{eq:ass_B1.general_0b}
  \sup_{ N \in \N }
  \sup_{ v \in H_{ \gamma } }
  \left[
  \frac{
  \|
    B( v ) \mathscr{P}_N 
  \|_{ \mathrm{HS}( U, H_{ \gamma - \beta } ) }
  +
  N^{ \iota \vartheta }
  \,
  \|
    B( v ) ( \mathscr{P}_0 - \mathscr{P}_N )
  \|_{ \mathrm{HS}( U, H_{ \gamma - \chi } ) }
  }{
    1 + \| v \|_{ H_{ \gamma } }
  }
  \right]
  < \infty
  .
\end{equation}
Then it holds for all 
$ \delta \in [ 0, \vartheta - \nicefrac{ 1 }{ p } ) $,
$ \varepsilon \in (0,\infty) $
that
\begin{equation}
\label{eq:convGlobLip}
  \sup_{ N \in \N }
  \left[
  \E \Bigl[
      \| X^N \|_{ \calC^{ \delta }( [0,T], \left\| \cdot \right\|_{ H_\gamma } )  }^p
  \Bigr]
  +
  N^{ \iota \, ( \vartheta - \delta - \nicefrac{1}{p} - \eps ) }
  \Bigl( \E \Bigl[
      \| X^0 - X^N \|_{ \calC^{ \delta }( [0,T], \left\| \cdot \right\|_{ H_\gamma } )  }^p
  \Bigr] \Bigr)^{ \nicefrac{1}{p} }
  \right]
  < \infty 
  .
\end{equation}
\end{corollary}
\begin{proof}[Proof of Corollary \ref{cor:convGlobLip}]
Throughout this proof let $ \eta \in \R $ 
be the real number given by
$
  \eta = \max\{ \alpha, 2 \beta \}
$
and let 
$ \theta^N \in \mathscr{P}_T $,
$ N \in \N $,
be a sequence of sets such that
\begin{equation}
\sup_{ N \in \N }
        \biggl[
            \frac{d_{ \max }( \theta^N )}{ N^{ - \iota } }
            + \frac{ N^{ - \iota } }{d_{ \min }( \theta^N )}
        \biggr]
< \infty.
\end{equation}
In particular, this ensures that
$
  \limsup_{ N \to \infty }
  d_{ \max }( \theta^N ) = 0 
$.
In addition, 
Corollary~\ref{cor:convGlobLip_0} 
proves that
\begin{equation}
\label{eq:proof_Lip_to_apply_hoelder}
\begin{split}
&
  \sup_{ N \in \N }
  \left[
  \left| d_{ \max }( \theta^N ) \right|^{ - \vartheta }
  \sup_{ t \in \theta^N }
  \|
    X^0_t - X^N_t
  \|_{
    \mathscr{L}^p( \P; \left\| \cdot \right\|_{  H_{ \gamma }  } ) 
  }
  \right]
\\ & 
\leq
  \biggl[
    \sup_{ N \in \N }
    \tfrac{ 
      \left|
        d_{ \max }( \theta^N )
      \right|^{ - \vartheta }
    }{
      N^{ \iota \vartheta } 
    }
  \biggr]
  \biggl(
  \sup_{ N \in \N }
  \sup_{ t \in \theta^N }
  N^{ \iota \vartheta }
  \,
  \|
    X^0_t - X^N_t
  \|_{
    \mathscr{L}^p( \P; \left\| \cdot \right\|_{  H_{ \gamma }  } ) 
  }
  \biggr)
\\ & 
\leq
  \biggl[
    \sup_{ N \in \N }
    \tfrac{
      N^{ - \iota } 
    }{
      d_{ \min }( \theta^N )
    }
  \biggr]^{ \vartheta }
  \biggl(
  \sup_{ N \in \N }
  \sup_{ t \in [0,T] }
  N^{ \iota \vartheta }
  \,
  \|
    X^0_t - X^N_t
  \|_{
    \mathscr{L}^p( \P; \left\| \cdot \right\|_{  H_{ \gamma }  } ) 
  }
  \biggr)
  < \infty
  .
\end{split}
\end{equation}
Next note that, e.g.,
\iftoggle{arXiv:v3}{%
Corollary~\ref{cor:temporal_regularity_SPDE}%
}{%
Corollary~3.2
in~\cite{CoxHutzenthalerJentzenVanNeervenWelti2016arXiv}%
}
shows
for all $ N \in \N_0 $,
$ 
  \varepsilon \in 
  ( 
    0, \min\{ 1 - \eta , \nicefrac{ 1 }{ 2 } - \beta \}
  )
$
that
\begin{align}
&
  \sup_{ 
    \substack{  
      t_1, t_2 \in [0,T] 
      ,
      \\
      t_1 \neq t_2
    }
  }
  \Biggl(
    \frac{
      \left|
        \min\{ t_1, t_2 \}
      \right|^{ \max\{ \gamma + \varepsilon - ( \gamma + \vartheta ) , 0 \} }
      \| X^N_{ t_1 } - X^N_{ t_2 } \|_{
        \mathscr{L}^p( \P; \left\| \cdot \right\|_{  H_{ \gamma }  } )
      }
    }{
      \left| t_1 - t_2 \right|^{ \varepsilon }
    }
  \Biggr)
\\ \nonumber & 
  \leq 
    \| 
      X_0^N 
    \|_{
      \mathscr{L}^p( 
        \P ; 
        \left\| \cdot \right\|_{  H_{ \min\{ \gamma + \vartheta, \gamma + \varepsilon \} }  }
      ) 
    }
  +
  \biggl[
    \sup_{ s \in [0,T] }
    \left\|
      P_N F( X_s^N ) 
    \right\|_{
      \mathscr{L}^p( \P; \left\| \cdot \right\|_{  H_{ \gamma - \eta }  } )
    }
  \biggr]
    \frac{
      2
      \,
      T^{ 
        ( 1 + \gamma - \eta - \min\{ \gamma + \vartheta, \gamma + \varepsilon \} ) 
      }
    }{
      \left( 1 - \eta - \varepsilon \right) 
    }
\\ \nonumber & \quad
  +
  \biggl[
    \sup_{ s \in [0,T] }
    \left\|
      P_N B( X_s^N ) \mathscr{P}_N
    \right\|_{
      \mathscr{L}^p( \P; \left\| \cdot \right\|_{  \mathrm{HS}( U, H_{ \gamma - \beta } )  } )
    }
  \biggr]
    \frac{ 
      \sqrt{ p \, ( p - 1 ) }
      \,
      T^{ ( \nicefrac{1}{2} + \gamma - \beta - \min\{ \gamma + \vartheta , \gamma + \varepsilon \} ) }
    }{
      \left( 
        1 - 2 \beta - 2 \varepsilon 
      \right)^{ \nicefrac{1}{2} }
    }
  < \infty
  .
\end{align}
This and the fact that
$
    \min\{ 1 - \eta , \nicefrac{ 1 }{ 2 } - \beta \}
  =
    \min\{ 1 - \max\{ \alpha, 2 \beta \} , \nicefrac{ 1 }{ 2 } - \beta \}
  =
    \min\{ 1 - \alpha, \nicefrac{ 1 }{ 2 } - \beta \}
  > \vartheta > 0
$
imply that
\begin{equation}
\begin{split}
&
  \sup_{ N \in \N_0 }
  \sup_{ 
    \substack{  
      t_1, t_2 \in [0,T] 
      ,
      \\
      t_1 \neq t_2
    }
  }
  \Biggl(
    \frac{
      \| X^N_{ t_1 } - X^N_{ t_2 } \|_{
        \mathscr{L}^p( \P; \left\| \cdot \right\|_{  H_{ \gamma }  } )
      }
    }{
      \left| t_1 - t_2 \right|^{ \vartheta }
    }
  \Biggr)
\\ & 
  \leq 
    \sup_{ N \in \N_0 }
    \| 
      X_0^N 
    \|_{
      \mathscr{L}^p( 
        \P ; 
        \left\| \cdot \right\|_{  H_{ \gamma + \vartheta }  } 
      ) 
    }
  +
  \biggl[
    \sup_{ N \in \N_0 }
    \sup_{ s \in [0,T] }
    \|
      F( X_s^N ) 
    \|_{
      \mathscr{L}^p( \P; \left\| \cdot \right\|_{  H_{ \gamma - \eta }  } )
    }
  \biggr]
    \frac{
      2
      \,
      T^{ 
        ( 1 - \eta - \vartheta ) 
      }
    }{
      \left( 1 - \eta - \vartheta \right) 
    }
\\ & \quad
  +
  \biggl[
    \sup_{ N \in \N_0 }
    \sup_{ s \in [0,T] }
    \|
      B( X_s^N ) \mathscr{P}_N
    \|_{
      \mathscr{L}^p( \P; \left\| \cdot \right\|_{  \mathrm{HS}( U, H_{ \gamma - \beta } )  } )
    }
  \biggr]
    \frac{ 
      \sqrt{ p \, ( p - 1 ) }
      \,
      T^{ 
        ( 
          \nicefrac{1}{2} - \beta - \vartheta
        ) 
      }
    }{
      \left( 
        1 - 2 \beta - 2 \vartheta 
      \right)^{ \nicefrac{1}{2} }
    }
  .
\end{split}
\end{equation}
Corollary~\ref{cor:convGlobLip_0} 
and estimate~\eqref{eq:ass_B1.general_0b}
hence prove that
\begin{equation}
\begin{split}
&
  \sup_{ N \in \N_0 }
  |
    X^N
  |_{
    \calC^{ \vartheta }( [0,T] ,
        \left\| \cdot \right\|_{  \mathscr{L}^p( \P ; \left\| \cdot \right\|_{  H_{ \gamma }  } )  } )
  }
  \\ &
\leq
    \| 
      X_0^0 
    \|_{
      \mathscr{L}^p( 
        \P ; 
        \left\| \cdot \right\|_{  H_{ \gamma + \vartheta }  } 
      ) 
    }
  +
  \biggl[
    \sup_{ N \in \N_0 }
    \sup_{ s \in [0,T] }
    \|
      F( X_s^N ) 
    \|_{
      \mathscr{L}^p( \P; \left\| \cdot \right\|_{  H_{ \gamma - \eta }  } )
    }
  \biggr]
    \frac{
      2
      \,
      T^{ 
        ( 1 - \eta - \vartheta ) 
      }
    }{
      \left( 1 - \eta - \vartheta \right) 
    }
\\ & \quad
  +
  \biggl[
    \sup_{ N \in \N_0 }
    \sup_{ s \in [0,T] }
    \|
      B( X_s^N ) \mathscr{P}_N
    \|_{
      \mathscr{L}^p( \P; \left\| \cdot \right\|_{  \mathrm{HS}( U, H_{ \gamma - \beta } )  } )
    }
  \biggr]
    \frac{ 
      \sqrt{ p \, ( p - 1 ) }
      \,
      T^{ 
        ( 
          \nicefrac{1}{2} - \beta - \vartheta
        ) 
      }
    }{
      \left( 
        1 - 2 \beta - 2 \vartheta 
      \right)^{ \nicefrac{1}{2} }
    }
  < \infty
  .
\end{split}
\end{equation}
This, \eqref{eq:proof_Lip_to_apply_hoelder},
and the fact that
$
  \vartheta \in ( \nicefrac{ 1 }{ p } , 1 ]
$
allow us to apply
Corollary~\ref{cor:hoelder3}
to obtain for all 
$ \delta \in [ 0, \vartheta - \nicefrac{ 1 }{ p } ) $,
$ \varepsilon \in (0,\infty) $
that
\begin{equation}
\begin{split}
  \sup_{ N \in \N }
  \bigg[
  &
  \E \Bigl[
      \| X^N \|_{ \calC^{ \delta }( [0,T], \left\| \cdot \right\|_{ H_\gamma } )  }^p
  \Bigr]
  \\ & 
  +
  \bigl|
    d_{ \max }( \theta^N )
  \bigr|^{ - ( \vartheta - \delta - \nicefrac{1}{p} - \eps ) }
  \Bigl( \E \Bigl[
      \| X^0 - X^N \|_{ \calC^{ \delta }( [0,T], \left\| \cdot \right\|_{ H_\gamma } )  }^p
  \Bigr] \Bigr)^{ \nicefrac{1}{p} }
  \bigg]
  < \infty .
\end{split}
\end{equation}
Combining this with the 
fact that
$
\sup_{ N \in \N }
        \bigl[
            \frac{d_{ \max }( \theta^N )}{ N^{ - \iota } }
            \bigr] < \infty
$
completes the proof of
Corollary~\ref{cor:convGlobLip}.
\end{proof}

\subsection{Almost sure convergence in H\"{o}l\-der norms for Galerkin approximations of SEEs with semi-globally Lipschitz continuous nonlinearities}
\label{subsec:convLocLip}

The proof of the following corollary
employs a standard localisation 
argument; see, e.g., \cite{g98,p01}.

\begin{corollary}
\label{cor:convLocLip}
Assume the setting in Subsection~\ref{ssec:setting}, 
let
$ 
  \vartheta \in 
  ( 0, \min\{ 1 - \alpha , \nicefrac{ 1 }{ 2 } - \beta \} ) 
$,
assume that
$ 
  \P\big( X_0^0 \in H_{ \gamma + \vartheta } 
  \big) = 1 
$,
and assume for all non-empty bounded sets $ E \subseteq H_{ \gamma } $ that
\begin{equation}
\label{eq:convLocLip_assumption}
  \sup_{ N \in \N }
  \sup_{ v \in E }
  \left[
  \frac{
  \|
    B( v ) \mathscr{P}_N 
  \|_{ \mathrm{HS}( U, H_{ \gamma - \beta } ) }
  +
  N^{ \iota \vartheta }
  \,
  \|
    B( v ) ( \mathscr{P}_0 - \mathscr{P}_N )
  \|_{ \mathrm{HS}( U, H_{ \gamma - \chi } ) }
  }{
    1 + \| v \|_{ H_{ \gamma } }
  }
  \right]
  < \infty
  .
\end{equation}
Then it holds for all 
$ \delta \in [0,\vartheta) $,
$ \eps \in (0,\infty) $
that
\begin{equation}
\P \biggl(
  \sup_{ N \in \N }
  \big[
    N^{ \iota ( \vartheta - \delta - \eps ) }
    \,
    \|
      X^0 - X^N
    \|_{
      \calC^{ \delta }( [0,T], \left\| \cdot \right\|_{  H_{ \gamma }  } ) 
    }
  \big]
  <
  \infty
\biggr) = 1.
\end{equation}
\end{corollary}

\begin{proof}[Proof of Corollary~\ref{cor:convLocLip}]
Throughout this proof we assume w.l.o.g.\ that $ X^0_0( \Omega ) \subseteq H_{ \gamma + \vartheta } $,
let 
$
  \delta \in [ 0, \vartheta )
$,
let 
$
  \phi_{ r, M } \colon H_r \to H_r
$,
$ r \in \R $,
$ M \in (0,\infty) $,
be the mappings which satisfy
for all $ r \in \R $,
$ M \in (0,\infty) $,
$ v \in H_r $ that
\begin{equation}
  \phi_{ r, M }( v )
  =
  v \cdot
  \min\biggl\{ 
    1, \frac{M + 1}{1 + \| v \|_{ H_r }}
  \biggr\},
\end{equation}
let
$ \xi_M \colon \Omega \to H_{ \gamma } $,
$ M \in \N $,
be the mappings which satisfy
for all $ M \in \N $ that
$
  \xi_M =
  \phi_{ \gamma + \vartheta, M }( X^0_0 )
$,
let 
$ F_M \colon H_{ \gamma } \to H_{ \gamma - \alpha } $,
$ M \in \N $,
and 
$ B_M \colon H_{ \gamma } \to \mathrm{HS}( U, H_{ \gamma - \beta } ) $,
$ M \in \N $,
be the mappings which satisfy
for all $ M \in \N $ that
$
  F_M = F \circ \phi_{ \gamma, M } 
$
and
$
  B_M = B \circ \phi_{ \gamma, M } 
$,
and let 
$
  S_M \subseteq H_{ \gamma }
$,
$ M \in \N $,
be the sets which satisfy for all 
$ M \in \N $ that
$
  S_M = \{ 
    v \in H_{ \gamma } \colon 
    \| v \|_{ H_{ \gamma } } \leq M + 1
  \}
$.
Observe that it holds
for all $ v, w \in H_{ \gamma } $, $ M \in \N $
that
\begin{equation}
\begin{split}
&
  \left\|
    \phi_{ \gamma, M }( v )
    -
    \phi_{ \gamma, M }( w )
  \right\|_{ H_{ \gamma } }
\\ & =
  \left\|
    \frac{
      v 
      \,
      ( 1 + \| w \|_{ H_{ \gamma } } )
      \min\{ 1 + \| v \|_{ H_{ \gamma } } , M + 1 \}
      -
      w 
      \,
      ( 1 + \| v \|_{ H_{ \gamma } } )
      \min\{ 1 + \| w \|_{ H_{ \gamma } } , M + 1 \}
    }{
      ( 1 + \| v \|_{ H_{ \gamma } } )
      \,
      ( 1 + \| w \|_{ H_{ \gamma } } )
    }
  \right\|_{ H_{ \gamma } }
\\ & \leq
  \left\| v - w \right\|_{ H_{ \gamma } }
\\ &
  +
  \left\|
    \frac{
      w 
      \,
      \big[
        ( 1 + \| w \|_{ H_{ \gamma } } )
        \min\{ 1 + \| v \|_{ H_{ \gamma } } , M + 1 \}
        -
        ( 1 + \| v \|_{ H_{ \gamma } } )
        \min\{ 1 + \| w \|_{ H_{ \gamma } } , M + 1 \}
      \big]
    }{
      ( 1 + \| v \|_{ H_{ \gamma } } )
      \,
      ( 1 + \| w \|_{ H_{ \gamma } } )
    }
  \right\|_{ H_{ \gamma } }
\\ & \leq
  \left\| v - w \right\|_{ H_{ \gamma } }
\\ & \quad
  +
    \frac{
      \left|
        ( 1 + \| w \|_{ H_{ \gamma } } )
        \min\{ 1 + \| v \|_{ H_{ \gamma } } , M + 1 \}
        -
        ( 1 + \| v \|_{ H_{ \gamma } } )
        \min\{ 1 + \| w \|_{ H_{ \gamma } } , M + 1 \}
      \right|
    }{
      ( 1 + \| v \|_{ H_{ \gamma } } )
    }
    .
\end{split}
\end{equation}
This ensures for all $ v, w \in H_{ \gamma } $, $ M \in \N $
that
\begin{equation}
\begin{split}
&
  \left\|
    \phi_{ \gamma, M }( v )
    -
    \phi_{ \gamma, M }( w )
  \right\|_{ H_{ \gamma } }
\\ & \leq
  \left\| v - w \right\|_{ H_{ \gamma } }
  +
    \frac{
        \big| 
          \| w \|_{ H_{ \gamma } } 
          -
          \| v \|_{ H_{ \gamma } } 
        \big|
        \min\{ 1 + \| v \|_{ H_{ \gamma } } , M + 1 \}
    }{
      ( 1 + \| v \|_{ H_{ \gamma } } )
    }
\\ & \quad
  +
    \frac{
        ( 1 + \| v \|_{ H_{ \gamma } } )
      \left|
        \min\{ 1 + \| v \|_{ H_{ \gamma } } , M + 1 \}
        -
        \min\{ 1 + \| w \|_{ H_{ \gamma } } , M + 1 \}
      \right|
    }{
      ( 1 + \| v \|_{ H_{ \gamma } } )
    }
\\ & \leq
  \left\| v - w \right\|_{ H_{ \gamma } }
  +
        \left| 
          \| w \|_{ H_{ \gamma } } 
          -
          \| v \|_{ H_{ \gamma } } 
        \right|
\\ & \quad
      +
      \left|
        \min\{ 1 + \| v \|_{ H_{ \gamma } } , M + 1 \}
        -
        \min\{ 1 + \| w \|_{ H_{ \gamma } } , M + 1 \}
      \right|
\\ & \leq
  3
  \left\| v - w \right\|_{ H_{ \gamma } }
  .
\end{split}
\end{equation}
Hence, we obtain for all $ M \in \N $ that
$
  |
    \phi_{ \gamma, M }
  |_{
    \calC^1( H_{ \gamma }, \left\| \cdot \right\|_{  H_{ \gamma }  } )
  }
  \leq 3
$.
This, the fact that
$
  \forall \, M \in \N \colon
  |
    F|_{ S_M }
  |_{
    \calC^1( S_M, \left\| \cdot \right\|_{  H_{ \gamma - \alpha }  } )
  }
  +
  |
    B|_{ S_M }
  |_{
    \calC^1( S_M, \left\| \cdot \right\|_{  \mathrm{HS}( U, H_{ \gamma - \beta } )  } )
  }
  +
  |
    \phi_{ \gamma, M }
  |_{
    \calC^1( H_{ \gamma }, \left\| \cdot \right\|_{  H_{ \gamma }  } )
  }
  < \infty
$,
and the fact that
$ 
  \forall \, M \in \N \colon
  \phi_{ \gamma, M }( H_{ \gamma } )
  \subseteq
  S_M
$
ensure that it holds
for all $ M \in \N $, $ p \in [1,\infty) $
that
\begin{equation}
  |
    F_M
  |_{
    \calC^1( H_{ \gamma }, \left\| \cdot \right\|_{  H_{ \gamma - \alpha }  } )
  }
  +
  |
    B_M
  |_{
    \calC^1( H_{ \gamma }, \left\| \cdot \right\|_{  \mathrm{HS}( U, H_{ \gamma - \beta } )  } )
  }
  +
  \E\big[
    \| \xi_M \|^p_{ H_{ \gamma + \vartheta } }
  \big]
  < \infty.
\end{equation}
E.g.,
\iftoggle{arXiv:v3}{%
Proposition~\ref{prop:existence_continuous_SPDEs}%
}{%
Proposition~3.6
in~\cite{CoxHutzenthalerJentzenVanNeervenWelti2016arXiv}%
}
hence proves that there exist
$
( \mathscr{F}_t )_{ t \in [0,T] }
$/$ \mathscr{B}( H_{ \gamma } ) $-adapted stochastic processes
$
  \mathscr{X}^{ N, M } \colon [0,T] \times \Omega \rightarrow H_{ \gamma }
$,
$ N \in \N_0 $, $ M \in \N $,
with continuous sample paths
such that it holds for all 
$ N \in \N_0 $, $ M \in \N $, $ t \in [0,T] $
that
\begin{equation}
\label{eq:GlobLipGalerkin_solution2}
\begin{split}
 \bigl[ \mathscr{X}_t^{ N, M } \bigr]_{ \P, \mathscr{B}( H_\gamma ) }
 & = 
\biggl[  e^{ t A } P_N \xi_M
  +
  \int_0^t
    e^{ ( t - s ) A }
    P_N F_M( \mathscr{X}^{ N, M }_s )
  \ds
  \biggr]_{ \P, \mathscr{B}( H_\gamma ) }
\\ & \quad
  +
  \int_0^t
    e^{ ( t - s ) A }
    P_N B_M( \mathscr{X}_s^{ N, M } ) \mathscr{P}_N
  \dWs
\end{split}
\end{equation}
(cf., e.g., Theorem~7.1 in
van~Neerven, Veraar, \& Weis~\cite{VanNeervenVeraarWeis2008}).
We now introduce a bit more notation.
Let 
$ 
  \tau_{ N, M } \colon \Omega \rightarrow [0,T]
$,
$ M \in \N $,
$ N \in \N_0 $,
be the mappings which satisfy
for all $ M \in \N $, $ N \in \N_0 $ that
\begin{equation}
\label{eq:def_tau_NM}
  \tau_{ N, M } 
= 
  \min\!\left\{
    T 
    \,
    \mathbbm{1}_{
      \{
        \| X_0^0 \|_{ H_{ \gamma + \vartheta } } \leq M 
      \}
    }
  ,
    \inf\!\left(
      \big\{ 
        t \in [0,T] 
        \colon
        \| \mathscr{X}_t^{ N, M } \|_{ H_{ \gamma } }
        \geq M
      \big\}
      \cup
      \{ T \}
    \right)
  \right\}
  ,
\end{equation}
let $ \varUpsilon \in \calF $
be the set given by 
\begin{align*}
\label{eq:def_varOmega}
&
  \varUpsilon
=  
\\ &
  \left[
    \cap_{ N \in \N_0 }
      \cup_{ M \in \N }
      \cap_{ m \in \{ M, M + 1, \ldots \} }
      \{ \tau_{ N, m } = T \}
  \right]
  \cap 
  \left[ 
    \cap_{ 
      M \in \N ,
      N \in \N_0
    }
    \left\{
      \| \mathscr{X}^{ N, M } \|_{
        \calC^{ \delta }( [0,T], \left\| \cdot \right\|_{  H_{ \gamma }  } )
      }
      < \infty
    \right\}
  \right]
\\ &
  \cap
  \left[
  \cap_{
    M \in \N, N \in \N_0
  }
  \left(
    \left\{
      \| X_0^0 \|_{ H_{ \gamma + \vartheta } } 
      > 
      M
    \right\}
    \cup
    \big\{
    \forall \, t \in [0,T] \colon 
    \mathscr{X}^{ N, M }_{ \min\{ t, \tau_{ N, M } \} }
    =
    X^N_{ \min\{ t, \tau_{ N, M } \} }
    \big\}
  \right)
  \right]
\\ & \yesnumber
  \cap
  \left[
    \cap_{ M, n \in \N } 
      \left\{ 
        \sup\nolimits_{ N \in \N }
        \big(
          N^{ \iota ( \vartheta - \delta - \nicefrac{ 1 }{ n } ) }
          \,
          \|
            \mathscr{X}^{ 0, M }
            - 
            \mathscr{X}^{ N, M }
          \|_{
            \calC^{ \delta }( [0,T], \left\| \cdot \right\|_{  H_{ \gamma }  } ) 
          }
        \big)
        < \infty 
      \right\}
    \right]
  ,
\end{align*}
let 
$ 
  \mathscr{M} \colon 
  \varUpsilon  
  \to \N 
$
be the mapping which satisfies
for all 
$ \omega \in \varUpsilon $ that
\begin{equation}
\label{eq:def_calM}
  \mathscr{M}( \omega ) 
  =
  \min\!\big\{ 
    M \in \N \cap
    \big(
      \| X^0_0( \omega ) \|_{ H_{ \gamma + \vartheta } }
      ,
      \infty
    \big)
    \colon
    \forall \, m \in \{ M, M + 1, \ldots \} \colon
    \tau_{ 0, m }( \omega ) = T
  \big\}
  ,
\end{equation}
and let 
$
  \mathscr{N} \colon 
  \varUpsilon 
  \to \N  
$
be the mapping which satisfies
for all 
\begin{equation*}
  \omega \in
  \varUpsilon 
  \subseteq 
  \left\{
    \mathfrak{w} \in \Omega 
    \colon
    \Bigl[
    \forall \, M \in \N
    \colon
        \limsup\nolimits_{ N \to \infty }
          \|
            \mathscr{X}^{ 0, M }( \mathfrak{w} )
            - 
            \mathscr{X}^{ N, M }( \mathfrak{w} )
          \|_{
            \calC^{ \delta }( [0,T], \left\| \cdot \right\|_{  H_{ \gamma }  } ) 
          }
        = 0
    \Bigr]
  \right\}
\end{equation*}
that
\begin{equation}
\label{eq:def_calN}
  \mathscr{N}( \omega ) 
  =
  \min\!\left\{ 
    N \in \N \colon
    \sup\nolimits_{ n \in \{ N, N + 1, \ldots \} }
    \|
      \mathscr{X}^{ 0, 2 \mathscr{M}( \omega ) }( \omega )
      -
      \mathscr{X}^{ n, 2 \mathscr{M}( \omega ) }( \omega )
    \|_{
      C( [0,T] , \left\| \cdot \right\|_{  H_{ \gamma }  } )
    }
    < 1
  \right\}
  .
\end{equation}
Observe that \eqref{eq:def_varOmega} ensures 
for all 
$ \omega \in \varUpsilon $,
$ N \in \N_0 $,
$ M \in \N $,
$ t \in [ 0, \tau_{ N, M }( \omega) ] $
with
$
  M \geq \| X^0_0( \omega ) \|_{ H_{ \gamma + \vartheta } }
$
that
\begin{equation} 
\label{eq:X_the_same_0}
  \mathscr{X}^{ N, M }_t( \omega )
  =
  X^N_t( \omega )
  .
\end{equation}
This, the fact that 
$
  \forall \, \omega \in \varUpsilon ,
  N \in \N_0 
  \colon
  \exists \, M \in \N \colon
  \forall \, m \in \{ M, M + 1, \ldots \} \colon
  \tau_{ N, m }( \omega ) = T
$,
and the fact that
$
  \forall \, \omega \in \varUpsilon,
  N \in \N_0, m \in \N \colon
  \|
    \mathscr{X}^{ N, m }( \omega )
  \|_{ 
    \calC^{ \delta }( [0,T] , \left\| \cdot \right\|_{  H_{ \gamma }  } )
  }
  < \infty
$
prove that it holds
for all $ \omega \in \varUpsilon $, $ N \in \N_0 $
that
\begin{equation}
\label{eq:Hoelder_reg_X}
  \| X^N( \omega ) \|_{
    \calC^{ \delta }( [0,T] , \left\| \cdot \right\|_{  H_{ \gamma }  } )
  }
  < \infty
  .
\end{equation}
Next note that 
\eqref{eq:def_calM}
ensures
for all $ \omega \in \varUpsilon $,
$ M \in \{ \mathscr{M}( \omega ), \mathscr{M}( \omega ) + 1, \ldots \} $
that
\begin{equation}
\label{eq:property_calM}
  \tau_{ 0, M }( \omega ) = T
\qquad 
  \text{and}
\qquad
  M \geq 
  \mathscr{M}( \omega )
  >
  \| X^0_0( \omega ) \|_{ H_{ \gamma + \vartheta } } 
  .
\end{equation}
This and \eqref{eq:X_the_same_0} show for all 
$ \omega \in \varUpsilon $,
$ M \in \{ \mathscr{M}( \omega ) , \mathscr{M}( \omega ) + 1, \ldots \} $,
$ t \in [ 0, T ] $
that
\begin{equation} 
\label{eq:X_the_same_1}
  \mathscr{X}^{ 0, M }_t( \omega )
=
  X^0_t( \omega )
=
  \mathscr{X}^{ 0, \mathscr{M}( \omega ) }_t( \omega )
  .
\end{equation}
This,
\eqref{eq:property_calM},
and \eqref{eq:def_tau_NM}
prove
for all $ \omega \in \varUpsilon $
that
\begin{equation}
  \sup\nolimits_{ t \in [0,T] }
  \| \mathscr{X}^{ 0, 2 \mathscr{M}( \omega ) }_t( \omega ) \|_{ 
    H_{ \gamma }
  }
  =
  \sup\nolimits_{ t \in [0,T] }
  \| \mathscr{X}^{ 0, \mathscr{M}( \omega ) }_t( \omega ) \|_{ 
    H_{ \gamma }
  }
  \leq \mathscr{M}( \omega ).
\end{equation}
The triangle inequality
and \eqref{eq:def_calN}
hence assure
for all $ \omega \in \varUpsilon $,
$
    N \in 
    \{ 
      \mathscr{N}( \omega ), 
      \mathscr{N}( \omega ) + 1 ,  
      \ldots
    \}
$
that
\begin{equation}
\begin{aligned}
&
  \sup\nolimits_{ t \in [0,T] }
  \| 
    \mathscr{X}^{ N, 2 \mathscr{M}( \omega ) }_t( \omega ) 
  \|_{
    H_{ \gamma }
  }
\\ &
\leq
  \sup\nolimits_{ t \in [0,T] }
  \| 
    \mathscr{X}^{ 0 , 2 \mathscr{M}( \omega ) }_t( \omega ) 
  \|_{
    H_{ \gamma }
  }
  +
  \sup\nolimits_{ t \in [0,T] }
  \|  
    \mathscr{X}^{ 0 , 2 \mathscr{M}( \omega ) }_t( \omega ) 
    - 
    \mathscr{X}^{ N , 2 \mathscr{M}( \omega ) }_t( \omega ) 
  \|_{
    H_{ \gamma } 
  }
\\ & < 
  \sup\nolimits_{ t \in [0,T] }
  \| 
    \mathscr{X}^{ 0 , 2 \mathscr{M}( \omega ) }_t( \omega ) 
  \|_{
    H_{ \gamma } 
  }
  +
  1
\leq
  \mathscr{M}( \omega )
  +
  1
\leq 
  2 \mathscr{M}( \omega )
  .
\end{aligned}
\end{equation}
This and the fact that 
$
  \forall \, \omega \in \varUpsilon \colon
  \| X^0_0( \omega ) \|_{ H_{ \gamma + \vartheta } }
  <
  \mathscr{M}( \omega )
  \leq 
  2 \mathscr{M}( \omega )
$
prove
for all 
$ \omega \in \varUpsilon $,
$
    N \in 
    \{ 
      \mathscr{N}( \omega ), \mathscr{N}( \omega ) + 1 , \ldots
    \}
$
that
$
  \tau_{ N, 2 \mathscr{M}( \omega ) }( \omega ) = T
$.
Again the fact that 
$ 
  \forall \, \omega \in \varUpsilon \colon
  \| X^0_0( \omega ) \|_{ H_{ \gamma + \vartheta } }
  <
  \mathscr{M}( \omega )
  \leq 
  2 \mathscr{M}( \omega )
$
and \eqref{eq:X_the_same_0}
hence show for all
$
  \omega \in \varUpsilon
$,
$
  N \in 
  \{ 
    \mathscr{N}( \omega ), \mathscr{N}( \omega ) + 1 , \ldots 
  \}
$,
$ t \in [0,T] $
that
$
  \mathscr{X}^{ N, 2 \mathscr{M}( \omega ) }_t( \omega )
  =
  X^N_t( \omega )
$.
This and \eqref{eq:X_the_same_1} prove for all
$ \omega \in \varUpsilon $,
$ \varepsilon \in ( 0 , \infty ) $
that
\begin{equation}
\begin{split}
& 
  \sup\nolimits_{ N \in \N }
  N^{ \iota ( \vartheta - \delta - \eps ) }
  \,
  \| 
    X^0( \omega ) - X^N( \omega ) 
  \|_{
    \calC^{ \delta }( [0,T] , \left\| \cdot \right\|_{  H_{ \gamma }  } )
  }
\\ & \leq
  \sup\nolimits_{ 
    N \in \{ 1, 2, \ldots, \mathscr{N}( \omega ) \}
  }
  N^{ \iota \vartheta }
  \,
  \| 
    X^0( \omega ) - X^N( \omega ) 
  \|_{
    \calC^{ \delta }( [0,T] , \left\| \cdot \right\|_{  H_{ \gamma }  } )
  }
\\ & \quad
  +
  \sup\nolimits_{ 
    N \in \{ \mathscr{N}( \omega ) , \mathscr{N}( \omega ) + 1 , \ldots \}
  }
  N^{ \iota ( \vartheta - \delta - \eps ) }
  \,
  \| 
    X^0( \omega ) - X^N( \omega ) 
  \|_{
    \calC^{ \delta }( [0,T] , \left\| \cdot \right\|_{  H_{ \gamma }  } )
  }
\\ & 
\leq
  \left[ 
    \mathscr{N}( \omega ) 
  \right]^{ \iota \vartheta }
  \left[
  \| 
    X^0( \omega )
  \|_{
    \calC^{ \delta }( [0,T], \left\| \cdot \right\|_{  H_{ \gamma }  } )
  }
  +
  \sup\nolimits_{ 
    N \in \{ 1, 2, \ldots, \mathscr{N}( \omega ) \} 
  }
  \| 
    X^N( \omega ) 
  \|_{
    \calC^{ \delta }( [0,T], \left\| \cdot \right\|_{  H_{ \gamma }  } )
  }
  \right]
\\ & \quad
  +
  \sup\nolimits_{ 
    N \in \{ \mathscr{N}( \omega ) , \mathscr{N}( \omega ) + 1 , \ldots \}
  }
  N^{ \iota ( \vartheta - \delta - \eps ) }  
  \,
  \| 
    \mathscr{X}^{ 0, 2 \mathscr{M}( \omega ) }( \omega ) 
    - 
    \mathscr{X}^{ N, 2 \mathscr{M}( \omega ) }( \omega ) 
  \|_{
    \calC^{ \delta }( [0,T], \left\| \cdot \right\|_{  H_{ \gamma }  } ) 
  }
  .
\end{split}
\end{equation}
Combining this with \eqref{eq:Hoelder_reg_X} 
and \eqref{eq:def_varOmega} ensures 
for all
$ \omega \in \varUpsilon $,
$ \varepsilon \in ( 0 , \infty ) $
that
\begin{align*}
& 
  \sup\nolimits_{ N \in \N }
  N^{ \iota ( \vartheta - \delta - \eps ) }
  \,
  \| 
    X^0( \omega ) - X^N( \omega ) 
  \|_{
    \calC^{ \delta }( [0,T] , \left\| \cdot \right\|_{  H_{ \gamma }  } )
  }
\\ & \leq \yesnumber
  \left[ 
    \mathscr{N}( \omega ) 
  \right]^{ \iota \vartheta }
  \textstyle
  \sum_{ N = 0 }^{ \mathscr{N}( \omega ) }
  \displaystyle
  \| 
    X^N( \omega ) 
  \|_{
    \calC^{ \delta }( [0,T], \left\| \cdot \right\|_{  H_{ \gamma }  } )
  }
\\ & \quad
  +
  \sup\nolimits_{ 
    N \in \{ \mathscr{N}( \omega ) , \mathscr{N}( \omega ) + 1 , \ldots \}
  }
  N^{ \iota ( \vartheta - \delta - \eps ) }  
  \,
  \| 
    \mathscr{X}^{ 0, 2 \mathscr{M}( \omega ) }( \omega ) 
    - 
    \mathscr{X}^{ N, 2 \mathscr{M}( \omega ) }( \omega ) 
  \|_{
    \calC^{ \delta }( [0,T], \left\| \cdot \right\|_{  H_{ \gamma }  } ) 
  }
  < \infty 
  .
\end{align*}
It thus remains to prove that
$ \P\big( \varUpsilon ) = 1 $
to complete the proof of Corollary~\ref{cor:convLocLip}.
For this observe that 
the assumption~\eqref{eq:convLocLip_assumption}
shows for all $ M \in \N $
that
\begin{equation}
\begin{split}
&
  \sup_{ N \in \N }
  \sup_{ v \in H_{ \gamma } }
  \left[
  \frac{
  \|
    B_M( v ) \mathscr{P}_N 
  \|_{ \mathrm{HS}( U, H_{ \gamma - \beta } ) }
  +
  N^{ \iota \vartheta }
  \,
  \|
    B_M( v ) ( \mathscr{P}_0 - \mathscr{P}_N )
  \|_{ \mathrm{HS}( U, H_{ \gamma - \chi } ) }
  }{
    1 + \| v \|_{ H_{ \gamma } }
  }
  \right]
\\
&
\leq
  \sup_{ N \in \N }
  \sup_{ v \in S_M }
  \left[
  \frac{
  \|
    B( v ) \mathscr{P}_N 
  \|_{ \mathrm{HS}( U, H_{ \gamma - \beta } ) }
  +
  N^{ \iota \vartheta }
  \,
  \|
    B( v ) ( \mathscr{P}_0 - \mathscr{P}_N )
  \|_{ \mathrm{HS}( U, H_{ \gamma - \chi } ) }
  }{
    1 + \| v \|_{ H_{ \gamma } }
  }
  \right]
  < \infty
  .
\end{split}
\end{equation}
Corollary~\ref{cor:convGlobLip}
hence proves for all 
$ p \in ( \nicefrac{ 1 }{ \vartheta } , \infty) $,
$ r \in [ 0, \vartheta - \nicefrac{ 1 }{ p } ) $,
$ \varepsilon \in (0,\infty) $,
$ M \in \N $
that
\begin{equation}
\label{eq:apply_CorLip}
  \sup_{ N \in \N } 
  \left[
  \E \Bigl[
      \| \mathscr{X}^{ N, M } \|_{ \calC^{ r }( [0,T], \left\| \cdot \right\|_{ H_\gamma } )  }^p
  \Bigr]
  +
  N^{ \iota \, ( \vartheta - r - \eps ) }
  \Bigl( \E \Bigl[
      \| \mathscr{X}^{ 0, M } - \mathscr{X}^{ N, M } \|_{ \calC^{ r }( [0,T], \left\| \cdot \right\|_{ H_\gamma } )  }^p
  \Bigr] \Bigr)^{ \nicefrac{1}{p} }
  \right]
  < \infty 
  .
\end{equation}
A standard Borel-Cantelli-type argument
(see, e.g., Lemma 2.1 in Kloeden \& Neuenkirch \cite{kn07}) 
hence ensures
for all $ \eps \in (0, \infty ) $, $ M \in \N $
that
\begin{equation}
\P \Bigl(
  \sup\nolimits_{ N \in \N }
  \big(
    N^{ \iota ( \vartheta - \delta - \eps ) }
    \|
      \mathscr{X}^{ 0, M } - \mathscr{X}^{ N, M }
    \|_{
      \calC^{ \delta }( [0,T], \left\| \cdot \right\|_{  H_{ \gamma }  } )
    }
  \big)
  < \infty
\Bigr) = 1.
\end{equation}
Hence, we obtain that
\begin{equation}
\label{eq:ASconvergence_GalApprox}
  \P\!\left(
    \forall \, M, n \in \N 
    \colon
  \sup\nolimits_{ N \in \N }
  \big[
    N^{ \iota ( \vartheta - \delta - \nicefrac{ 1 }{ n } ) }
    \,
    \|
      \mathscr{X}^{ 0, M } - \mathscr{X}^{ N, M }
    \|_{
      \calC^{ \delta }( [0,T], \left\| \cdot \right\|_{  H_{ \gamma }  } )
    }
  \big]
  < \infty
  \right) = 1
  .
\end{equation}
In addition,
\eqref{eq:apply_CorLip}
proves for all $ N \in \N_0 $, $ M \in \N $
that
$
  \P\big(
    \mathscr{X}^{ N, M } \in \calC^{ \delta }( [0,T], \left\| \cdot \right\|_{  H_{ \gamma }  } )
  \big) = 1
$.
This, in turn, ensures that 
\begin{equation}
\label{eq:local_Lip_in_Hoelder}
  \P\!\left(
    \forall \, M \in \N, N \in \N_0 \colon
    \mathscr{X}^{ N, M } \in \calC^{ \delta }( [0,T], \left\| \cdot \right\|_{  H_{ \gamma }  } )
  \right) = 1 .
\end{equation}
Next observe 
that it holds
for all $ t \in [0,T] $, $ M \in \N $,
$ N \in \N_0 $ 
that
\begin{align}
\nonumber
&  
  \bigl[
    \mathscr{X}^{ N, M }_t
    -
    e^{ t A }
    P_N
    \mathscr{X}_0^{ 0, M }
  \bigr]_{  \P, \mathscr{B}( H_\gamma )  }
  \mathbbm{1}_{
    \{ t \leq \tau_{ N, M } \}
  }
\\ &  \nonumber
  =
  \biggl( \biggl[
    \int_0^t 
      e^{ ( t - s ) A }
      P_N
      F_M( \mathscr{X}^{ N, M }_s )
    \ds
    \biggr]_{  \P, \mathscr{B}( H_\gamma )  }
    +
    \int_0^t 
      e^{ ( t - s ) A }
      P_N
      B_M( \mathscr{X}^{ N, M }_s )
    \mathscr{P}_N
    \dWs
    \biggr)
  \mathbbm{1}_{
    \{ t \leq \tau_{ N, M } \}
  }
\\ & \nonumber
  =
  \biggl( \biggl[
    \int_0^t 
      \mathbbm{1}_{
        \{ s < \tau_{ N, M } \}
      }
      \,
      e^{ ( t - s ) A }
      P_N
      F_M( \mathscr{X}^{ N, M }_s ) 
    \ds
    \biggr]_{  \P, \mathscr{B}( H_\gamma )  }
\\ & \qquad
    +
    \int_0^t 
      \mathbbm{1}_{
        \{ s < \tau_{ N, M } \}
      }
      e^{ ( t - s ) A }
      P_N 
      B_M( \mathscr{X}^{ N, M }_s )
    \mathscr{P}_N
    \dWs
  \biggr)
  \mathbbm{1}_{
    \{ t \leq \tau_{ N, M } \}
  }
\\ \nonumber &  
  = \biggl( \biggl[
    \int_0^t 
      \mathbbm{1}_{
        \{ s < \tau_{ N, M } \}
      }
      \,
      e^{ ( t - s ) A }
      P_N 
      F( \mathscr{X}^{ N, M }_s )
    \ds
    \biggr]_{  \P, \mathscr{B}( H_\gamma )  }
\\ & \qquad \nonumber
    +    
    \int_0^t 
      \mathbbm{1}_{
        \{ s < \tau_{ N, M } \}
      }
      \,
      e^{ ( t - s ) A }
      P_N 
      B( \mathscr{X}^{ N, M }_s )
      \mathscr{P}_N
    \dWs
  \biggr)
  \mathbbm{1}_{
    \{ t \leq \tau_{ N, M } \}
  }.
\end{align}
E.g.,
\iftoggle{arXiv:v3}{%
Proposition~\ref{prop:uniqueness_c_local}%
}{%
Proposition~3.7
in~\cite{CoxHutzenthalerJentzenVanNeervenWelti2016arXiv}%
}
hence shows for all $ N \in \N_0 $, $ M \in \N $
that
\begin{equation}
  \P\Bigl(
    \forall \, t \in [0,T] \colon 
    \mathbbm{1}_{
      \{ \mathscr{X}^{ N, M }_0 = X^N_0 \}
    }
    \mathscr{X}^{ N, M }_{ \min\{ t, \tau_{ N, M } \} }
    =
    \mathbbm{1}_{
      \{ \mathscr{X}^{ N, M }_0 = X^N_0 \}
    }
    X^N_{ \min\{ t, \tau_{ N, M } \} }
  \Bigr)
  = 1
\end{equation}
(cf., e.g., Lemma~8.2 in van~Neerven, Veraar, \& Weis~\cite{VanNeervenVeraarWeis2008}).
This implies for all 
$ N \in \N_0 $, $ M \in \N $
that
\begin{equation}
  \P\big(
    \{ 
      \| X_0^0 \|_{ H_{ \gamma + \vartheta } }
      > M
    \}
    \cup 
    \big\{
    \forall \, t \in [0,T] \colon 
    \mathscr{X}^{ N, M }_{ \min\{ t, \tau_{ N, M } \} }
    =
    X^N_{ \min\{ t, \tau_{ N, M } \} }
    \big\}
  \big)
  = 1.
\end{equation}
Hence, we obtain that
\begin{equation}
\label{eq:zero_set_1}
  \P\Bigl(
    \cap_{ M \in \N , N \in \N_0 }
    \Bigl[
    \{ 
      \| X_0^0 \|_{ H_{ \gamma + \vartheta } }
      > M
    \}
    \cup 
    \big\{
    \forall \, t \in [0,T] \colon 
    \mathscr{X}^{ N, M }_{ \min\{ t, \tau_{ N, M } \} }
    =
    X^N_{ \min\{ t, \tau_{ N, M } \} }
    \big\}
    \Bigr]
  \Bigr)
  = 1
  .
\end{equation}
In the next step we combine this with 
\eqref{eq:def_tau_NM} to obtain
for all $ M \in \N $, $ N \in \N_0 $
that
\begin{equation}
\label{eq:prop_tau_NM}
\P \Bigl(
  \tau_{ N, M } 
= 
  \min\!\left\{
    T 
    \,
    \mathbbm{1}_{
      \{
        \| X_0^0 \|_{ H_{ \gamma + \vartheta } } \leq M 
      \}
    }
  ,
    \inf\!\left(
      \big\{ 
        t \in [0,T] 
        \colon
        \| X_t^N \|_{ H_{ \gamma } }
        \geq M
      \big\}
      \cup
      \{ T \}
    \right)
  \right\}
\Bigr) = 1.
\end{equation}
This shows for all $ N \in \N_0 $, $ M_1, M_2 \in \N $
with $ M_1 \leq M_2 $
that 
$
  \P\big( \tau_{ N, M_1 } \leq \tau_{ N, M_2 } \big) = 1
$.
This, \eqref{eq:prop_tau_NM},
and the fact that
$
  \forall \, \omega \in \Omega,
  N \in \N_0 
  \colon
  \sup_{ t \in [0,T] }
  \| X_t^N( \omega ) \|_{ H_{ \gamma } }
  < \infty
$
imply that it holds for all $ N \in \N_0 $
that
$
  \P\big(
    \cup_{ M \in \N }
    \cap_{ m \in \{ M, M + 1, \ldots \} }
    \{ \tau_{ N, m } = T \}
  \big)
  = 1
$.
This, in turn, proves that
\begin{equation}
\label{eq:zero_set_2}
  \P\big(
    \cap_{ 
      N \in \N_0 
    }
    \cup_{ 
      M \in \N 
    }
    \cap_{ m \in \{ M, M + 1, \ldots \} }
    \{ \tau_{ N, m } = T \}
  \big)
  = 1
  .
\end{equation}
Combining \eqref{eq:zero_set_2}, 
\eqref{eq:local_Lip_in_Hoelder}, \eqref{eq:zero_set_1}, and \eqref{eq:ASconvergence_GalApprox}
proves that $ \P\big( \varUpsilon ) = 1 $.
The proof of Corollary~\ref{cor:convLocLip}
is thus completed.
\end{proof}

\section{Cubature methods in Banach spaces}\label{sec:cubature}
We first discuss in Subsection~\ref{sec:preliminaries} a number of preliminary definitions related to the Monte Carlo method in Banach spaces.
In Subsection~\ref{sec:montecarlo} we present an elementary error estimate for the Monte Carlo method in Corollary~\ref{cor:BanachMC}.
In Subsection~\ref{sec:mlmc} we then illustrate how expectations of Banach space valued functions of stochastic processes can be approximated.
\subsection{Preliminaries} \label{sec:preliminaries}
As mentioned in the introduction, the rate of convergence of Monte Carlo approximations in a Banach space depends on the so-called \emph{type} of the Banach space; cf., e.g., Section~9.2 in Ledoux \& Talagrand~\cite{lt91}.
In order to define the type of a Banach space, we first reconsider a few concepts from the literature.
\begin{definition}
Let $ ( \Omega, \calF, \P ) $ be a probability space,
let $ J $ be a set, 
and let 
$
  r_j \colon \Omega \to \{ -1, 1 \} 
$,
$ j \in J $,
be 
a family of 
independent random variables 
with
$
  \forall \, j \in J \colon
  \P\big( r_j = 1 \big) = \P\big( r_j = - 1 \big) 
$.
Then we say that 
$ ( r_j )_{ j \in J } $
is a $ \P $-Rademacher family.
\end{definition}

\begin{definition}
Let $ p \in (0,\infty) $
and let $ ( E, \left\| \cdot \right\|_E ) $ be an $ \R $-Banach space.
Then we denote by
$
  \mathscr{T}_p( E ) \in [0,\infty]
$
the extended real number given by
\begin{align}
&
  \mathscr{T}_p( E )
  =
  \sup\!\left(
    \left\{
    r \in [0,\infty) \colon
    \begin{array}{c}
    \exists \, \text{probability space } ( \Omega, \mathscr{F}, \P ) \colon
    \\
    \exists \, \text{$ \P $-Rademacher family } ( r_j )_{ j \in \N } \colon
    \\
    \exists \, k \in \N \colon
    \exists \, x_1, x_2, \dots, x_k \in E \backslash \{ 0 \} \colon
    \\
    r = 
    \frac{
      \left(
        \E \left[
          \|
            \sum_{ j = 1 }^k
            r_j x_j
          \|^p_E
        \right]
      \right)^{ 1 / p }
    }{
      \left(
        \sum_{ j = 1 }^k
        \left\| x_j \right\|_E^p
      \right)^{ 1 / p }
    }
    \end{array}
  \right\} \cup \{ 0 \}
  \right)
\end{align}
and we call $ \mathscr{T}_p( E ) $
the type $ p $-constant of $ E $.
\end{definition}

\begin{definition}
Let $ p \in (0,\infty) $
and let $ ( E, \left\| \cdot \right\|_E ) $ be an $ \R $-Banach space
which satisfies $ \mathscr{T}_p( E ) < \infty $.
Then
we say that $ ( E, \left\| \cdot \right\|_E ) $ has type $ p $
(we say that $ E $ has type $ p $).
\end{definition}

Note that it holds for 
all $ p \in (0,\infty) $,
all $ \R $-Banach spaces
$ ( E, \left\| \cdot \right\|_E ) $
with type $ p $,
all probability spaces
$ ( \Omega, \mathscr{F}, \P ) $,
all $ \P $-Rademacher families
$ ( r_j )_{ j \in \N } $,
and all 
$ k \in \N $,
$ x_1, x_2,\linebreak \dots, x_k \in E $
that
\begin{equation}
      \Biggl\| 
        \sum\limits_{ j = 1 }^k
        r_j x_j
      \Biggr\|_{
        \mathscr{L}^p( \P; \left\| \cdot \right\|_E )
      }
  \leq 
  \mathscr{T}_p( E )
  \Biggl(
    \sum_{ j = 1 }^k
    \left\| x_j \right\|^p_E
  \Biggr)^{ \! \nicefrac{1}{p} }
  .
\end{equation}
In addition, observe that it holds for all 
$ \R $-Banach spaces $ ( E, \left\| \cdot \right\|_E ) $,
all probability spaces $ ( \Omega, \mathscr{F}, \P ) $,
all $ \P $-Rademacher families $ ( r_j )_{ j \in \N } $,
and all $ p \in [2,\infty) $,
$ k \in \N $, $ x \in E \backslash \{ 0 \} $
that
\begin{equation}
\begin{split}
& 
  \mathscr{T}_p( E )
  \geq
  \frac{
      \Bigl\| 
        \sum_{ j = 1 }^k
        r_j x
      \Bigr\|_{
        \mathscr{L}^p( \P; \left\| \cdot \right\|_E )
      }
  }{
    \left[
    \sum_{ j = 1 }^k
    \|
      x
    \|^p_E
    \right]^{ \nicefrac{1}{p} }
  }
\geq
  \frac{
      \left\| x \right\|_E
      \Bigl\| 
        \sum_{ j = 1 }^k
        r_j
      \Bigr\|_{
        \mathscr{L}^2( \P; \left| \cdot \right| )
      }
  }{
    k^{ \nicefrac{ 1 }{ p } }
    \left\| x \right\|_E
  }
  =
  \frac{
    k^{ \nicefrac{ 1 }{ 2 } }
    \left\| x \right\|_E
  }{
    k^{ \nicefrac{ 1 }{ p } }
    \left\| x \right\|_E
  }
  =
  k^{ 
    ( \nicefrac{ 1 }{ 2 } - \nicefrac{ 1 }{ p } )
  }
  .
\end{split}
\end{equation}
In particular, it holds for all 
$ p \in (2,\infty) $
and all
$ \R $-Banach spaces
$ ( E, \left\| \cdot \right\|_E ) $
with $ E \neq \{ 0 \} $
that
$
  \mathscr{T}_p( E ) = \infty
$.
Furthermore, observe that Jensen's inequality together with the 
fact that it holds for all normed $ \R $-vector spaces
$ ( E , \left\| \cdot \right\|_E ) $
and all
$ p \in (0,\infty) $, 
$ q \in [p,\infty) $,
$ k \in \N $,
$ x_1, \dots, x_k \in E $
that
\begin{equation}
\textstyle
  \left(
    \sum_{ j = 1 }^k
    \left\| x_j \right\|^q_E
  \right)^{ \! \nicefrac{1}{q} }
  \leq
  \left(
    \sum_{ j = 1 }^k
    \left\| x_j \right\|^p_E
  \right)^{ \! \nicefrac{1}{p} }
\end{equation}
assures that it holds for all 
$ \R $-Banach spaces $ ( E, \left\| \cdot \right\|_E ) $
and all $ p, q \in (0,\infty) $
with $ p \leq q $ 
that
$ \mathscr{T}_p( E ) \leq \mathscr{T}_q( E ) $.
Hence, it holds for every $ \R $-Banach space 
$ ( E, \left\| \cdot \right\|_E ) $
that the function 
$
  (0,\infty) \ni p \mapsto \mathscr{T}_p( E ) \in [0,\infty]
$
is non-decreasing.
This and the triangle inequality 
ensure for all 
$ p \in (0,1] $
and all
$ \R $-Banach spaces 
$ ( E, \left\| \cdot \right\|_E ) $
with $ E \neq \{ 0 \} $
that
$ \mathscr{T}_p( E ) = 1 $.
In particular, note that it holds for all 
$ \R $-Banach spaces 
$ ( E , \left\| \cdot \right\|_E ) $
that
$ \sup_{ p \in (0,1] } \mathscr{T}_p( E ) \leq 1 < \infty $.
Additionally, observe 
that it holds for all 
$ p \in (0,2] $
and all
$ \R $-Hilbert spaces 
$ ( H, \left< \cdot , \cdot \right>_H , \left\| \cdot \right\|_H ) $
with 
$ H \neq \{ 0 \} $
that
$
  \mathscr{T}_p( H ) = 1
$.
Furthermore, we note that it holds
for 
every probability space 
$ ( \Omega , \mathscr{F}, \P ) $,
every $ p, q \in [1,\infty) $,
and 
every $ \R $-Banach space 
$ ( E, \left\| \cdot \right\|_E ) $
with type $ q $
that
$
  L^p( \P ; \left\| \cdot \right\|_E )
$
has type 
$ \min\{ p, q \} $;
cf., e.g., Proposition~7.1.4 in Hyt\"{o}nen et al.~\cite{HytonenNeervenVeraarWeis2017}, Section~9.2 in Ledoux \& Talagrand~\cite{lt91}, or Theorem~6.2.14 in Albiac \& Kalton~\cite{ak06}.
In particular, it holds 
for every $ p \in [1,\infty) $
and every probability space
$
  ( \Omega, \mathscr{F}, \P )
$
that 
$
  L^p( \P; \left| \cdot \right| )
$
has type $ \min\{ p, 2 \} $.

\begin{definition}
Let $ p, q \in (0, \infty) $.
Then we denote by
$
  \mathscr{K}_{ p, q } \in [0,\infty]
$
the extended real number given by
\begin{align}
\mathscr{K}_{ p, q } =
  \sup
    \left\{
    r \in [0,\infty) \colon
    \begin{array}{c}
      \exists \, \text{$ \R $-Banach space } ( E, \left\| \cdot \right\|_E ) \colon
      \\
      \exists \, \text{probability space } ( \Omega, \mathscr{F}, \P ) \colon \! \! \!
      \\
      \exists \, \text{$ \P $-Rademacher family } ( r_j )_{ j \in \N } \colon
      \exists \, k \in \N \colon
      \\
      \exists \, x_1, x_2, \dots, x_k \in E \backslash \{ 0 \} \colon
      r =  \frac{
            \left(
              \E \left[
                \|
                  \sum_{ j = 1 }^k
                  r_j x_j
                \|^p_E
              \right]
            \right)^{ 1 / p }
          }{
            \left(
              \E \left[
                \|
                  \sum_{ j = 1 }^k
                  r_j x_j
                \|^q_E
              \right]
            \right)^{ 1 / q }
          }
    \end{array}
  \right\}
\end{align}
and we call $ \mathscr{K}_{ p, q } $
the $ (p,q) $-Kahane-Khintchine constant.
\end{definition}

The celebrated {\em Kahane-Khintchine inequality}
asserts that it holds for all $ p , q \in (0,\infty) $
that $ \mathscr{K}_{ p, q } < \infty $;
see, e.g., Theorem~6.2.5 in Albiac \& Kalton~\cite{ak06}.
Observe that Jensen's inequality ensures for all $ p, q \in (0,\infty) $
with $ p \leq q $ that
$
  \mathscr{K}_{ p, q } = 1
$.
The nontrivial assertion of the Kahane-Khintchine inequality 
is the fact that it holds for all $ p, q \in (0,\infty) $
with $ p > q $ that
$
  \mathscr{K}_{ p, q } < \infty
$.
In our analysis below we also use the following two abbreviations.

\begin{definition}
Let $ p, q \in (0,\infty) $
and let $ ( E, \left\| \cdot \right\|_E ) $
be an $ \R $-Banach space.
Then we denote by 
$
  \varTheta_{ p, q }( E ) \in [0,\infty]
$
the extended real number given by
$
  \varTheta_{ p, q }( E )
  =
  2 \mathscr{T}_q( E ) \mathscr{K}_{ p, q } 
$.
\end{definition}

\begin{definition}
Let $ ( \Omega, \mathscr{F}, \P ) $ be a probability space,
let $ p \in (0,\infty) $,
let $ ( E, \left\| \cdot \right\|_E ) $ be an $ \R $-Banach space,
and let $ X \in \mathscr{L}^1( \P; \left\| \cdot \right\|_E ) $.
Then we denote by
$
  \sigma_{p,E}( X ) \in [0,\infty]
$
the extended real number 
given by
$
  \sigma_{p,E}( X ) = 
  \bigl(
    \E\bigl[ 
      \| X - \E[ X ] \|_E^p
    \bigr]
  \bigr)^{ \! \nicefrac{1}{p} }
$.
\end{definition}

\subsection{Monte Carlo methods in Banach spaces} \label{sec:montecarlo}
In this subsection we collect a few elementary results on sums of random variables with values in Banach spaces.
The next result, Lemma~\ref{lem:symm} below, can be found, e.g., in Section~2.2 of Ledoux \& Talagrand~\cite{lt91}.
\begin{lemma}[Symmetrisation lemma] \label{lem:symm}
Consider the notation in Subsection~\ref{notation},
let\linebreak
$ ( E, \norm{\cdot}_E ) $
be an $ \R $-Banach space,
let $ ( \Omega, \calF, \P ) $
be a probability space,
let $ \xi, \tilde{\xi} \in \mathscr{L}^0( \P; \left\| \cdot \right\|_E ) $
be independent mappings which satisfy
$ \E \bigl[ \norm{ \tilde{\xi} }_E \bigr] < \infty $ and $ \E[ \tilde{\xi} ] = 0 $,
and let $ \varphi \colon [0,\infty) \to [0,\infty) $
be a convex and non-decreasing function.
Then
\begin{equation}
\E \bigl[ \varphi( \norm{ \xi }_E ) \bigr]
\leq  \E \bigl[ \varphi( \norm{ \xi - \tilde{\xi} }_E ) \bigr].
\end{equation}
\end{lemma}
\begin{proof}[Proof of Lemma~\ref{lem:symm}]
Jensen's inequality assures that
\begin{equation}
\begin{split}
& \E \bigl[ \varphi( \norm{ \xi }_E ) \bigr]
=  \E \bigl[ \varphi( \norm{ \xi  -  \E[ \tilde{\xi} ] }_E ) \bigr]
= \int_\Omega    \varphi \biggl( \norm[\bigg]{
                                                    \int_\Omega  \xi(\omega)
                                                                          - \tilde{\xi}( \tilde{\omega} )   \,\P( \dd \tilde{\omega} )
                                                  }_E \biggr) \,   \P( \dd \omega ) \\
& \leq \int_\Omega    \varphi \biggl(
                                                    \int_\Omega  \norm{ \xi(\omega)
                                                                          - \tilde{\xi}( \tilde{\omega} ) }_E  \, \P( \dd \tilde{\omega} )
                                                  \biggr) \,   \P( \dd \omega ) \\
& \leq \int_\Omega \int_\Omega   \varphi( \norm{ \xi(\omega)
                                                                   - \tilde{\xi}( \tilde{\omega} )   }_E )
      \,\P( \dd \tilde{\omega} ) \, \P( \dd \omega ) \\
& =
\int_\Omega \int_\Omega
    \mathbbm{1}_{ \overline{ \xi( \Omega ) }^E \times \overline{ \tilde{\xi}( \Omega ) }^E }
    \bigl( \xi(\omega), \tilde{\xi}( \tilde{\omega} ) \bigr) \,
    \varphi( \norm{ \xi(\omega) - \tilde{\xi}( \tilde{\omega} )   }_E )
\,\P( \dd \tilde{\omega} ) \, \P( \dd \omega ) \\
& =
\int_E \int_E
    \mathbbm{1}_{ \overline{ \xi( \Omega ) }^E \times \overline{ \tilde{\xi}( \Omega ) }^E }( x, y ) \,
    \varphi( \norm{ x - y }_E )
\, \bigl( \tilde{\xi} (\P) \bigr)( \dd y ) \bigl( \xi (\P) \bigr)( \dd x ) \\
& =  \int_{ E \times E }
    \mathbbm{1}_{ \overline{ \xi( \Omega ) }^E \times \overline{ \tilde{\xi}( \Omega ) }^E }( x, y ) \,
    \varphi( \norm{ x - y }_E )
\, \bigl( ( \xi, \tilde{\xi} ) (\P) \bigr)( \dd x, \dd y )
= \E \bigl[ \varphi( \norm{ \xi - \tilde{\xi} }_E ) \bigr]. \\
\end{split}
\end{equation}
This completes the proof of Lemma~\ref{lem:symm}.
\end{proof}
\begin{corollary}[Symmetrisation corollary] \label{cor:symm}
Let $ ( E, \norm{\cdot}_E ) $
be an $ \R $-Banach space,
let $ ( \Omega, \calF, \P ) $
be a probability space,
let $ \xi, \tilde{\xi} \in \mathscr{L}^1( \P; \left\| \cdot \right\|_E ) $
be independent and identically distributed mappings which satisfy
$ \E[ \xi ] = 0 $,
and let $ \varphi \colon [0,\infty) \to [0,\infty) $
be a convex and non-decreasing function.
Then
\begin{equation}
\E \bigl[ \varphi( \norm{ \xi }_E ) \bigr]
\leq  \E \bigl[ \varphi( \norm{ \xi - \tilde{\xi} }_E ) \bigr]
\leq  \E \bigl[ \varphi( 2 \, \norm{ \xi }_E ) \bigr].
\end{equation}
\end{corollary}
\begin{proof}[Proof of Corollary~\ref{cor:symm}]
Lemma~\ref{lem:symm} shows that
\begin{equation}
\begin{split}
\E \bigl[ \varphi( \norm{ \xi }_E ) \bigr]
& \leq
\E \bigl[ \varphi( \norm{ \xi - \tilde{\xi} }_E ) \bigr]
\leq
\E \bigl[ \varphi( \norm{ \xi }_E + \norm{ \tilde{\xi} }_E ) \bigr]
=
\E \bigl[ \varphi( \tfrac{1}{2} \, 2 \, \norm{ \xi }_E + \tfrac{1}{2} \, 2 \, \norm{ \tilde{\xi} }_E ) \bigr]
\\ & \leq
\E \bigl[ \tfrac{1}{2} \, \varphi( 2 \, \norm{ \xi }_E ) + \tfrac{1}{2} \, \varphi( 2 \, \norm{ \tilde{\xi} }_E ) \bigr]
=
\tfrac{1}{2} \, \E \bigl[ \varphi( 2 \, \norm{ \xi }_E ) \bigr]
+ \tfrac{1}{2} \,  \E \bigl[ \varphi( 2 \, \norm{ \tilde{\xi} }_E ) \bigr]
\\ & =
\tfrac{1}{2} \, \E \bigl[ \varphi( 2 \, \norm{ \xi }_E ) \bigr]
+ \tfrac{1}{2} \,  \E \bigl[ \varphi( 2 \, \norm{ \xi }_E ) \bigr]
=
\E \bigl[ \varphi( 2 \, \norm{ \xi }_E ) \bigr].
\end{split}
\end{equation}
The proof of Corollary~\ref{cor:symm} is thus completed.
\end{proof}
As a straightforward application we obtain the following randomisation result,  cf., e.g., Lemma~6.3 in Ledoux \& Talagrand~\cite{lt91}.
\begin{lemma}[Randomisation] \label{lem:rand}
Let $ ( E, \norm{\cdot}_E ) $
be an $ \R $-Banach space,
let $ ( \Omega, \calF, \P ) $
be a probability space,
let $ k \in \N $,
let $ \xi_j \in \mathscr{L}^1( \P; \left\| \cdot \right\|_E ) $, $ j \in \{1, \ldots, k \} $,
satisfy
for all  $ j \in \{1, \ldots, k \} $ that $ \E[ \xi_j ] = 0 $,
and let $ r_j \colon \Omega \to \{ -1, 1 \} $, $ j \in \{1, \ldots, k \} $, be a $ \P $-Rademacher family such that
$ \xi_1, \xi_2, \ldots, \xi_k, r_1, r_2, \ldots, r_k $ are independent.
Then it holds for all $p \in [ 1, \infty ) $ that
\begin{equation}
\norm[\Bigg]{ \sum_{j=1}^k \xi_j }_{ \mathscr{L}^p( \P; \left\| \cdot \right\|_E ) }
\leq  2 \, \norm[\Bigg]{ \sum_{j=1}^k r_j \xi_j }_{ \mathscr{L}^p( \P; \left\| \cdot \right\|_E ) }.
\end{equation}
\end{lemma}
\begin{proof}[Proof of Lemma~\ref{lem:rand}]
Throughout this proof let
$ ( \mathbf{\Omega}, \bm{\mathcal{F}}, \mathbf{P} )
=
( \Omega \times \Omega, \calF \otimes \calF, \P \otimes \P ) $,
let $ \mathbf{r}_j \colon \mathbf{\Omega} \to \{ -1, 1 \} $, $ j \in \{1, \ldots, k \} $,
be the mappings which satisfy
for all $ \bm{\omega} = ( \omega, \tilde{\omega} ) \in \mathbf{\Omega} $,
$ j \in \{ 1, \ldots, k \} $
that
$ \mathbf{r}_j( \bm{\omega} ) = r_j( \omega ) $,
and
let $ \bm{\xi}_j \colon \mathbf{\Omega} \to E $, $ j \in \{1, \ldots, k \} $,
and $ \bm{\tilde{\xi}}_j \colon \mathbf{\Omega} \to E $, $ j \in \{1, \ldots, k \} $,
be the mappings which satisfy
for all $ \bm{\omega} = ( \omega, \tilde{\omega} ) \in \mathbf{\Omega} $, $ j \in \{ 1, \ldots, k \} $
that $ \bm{\xi}_j( \bm{\omega} ) = \xi_j( \omega ) $
and $ \bm{\tilde{\xi}}_j( \bm{\omega} ) = \xi_j( \tilde{\omega} ) $.
The fact that
\begin{equation}
\{ 0, 1 \} \times \{1, \ldots, k \}
\ni ( i, j )
\mapsto
\begin{cases}
\bm{\xi}_j -  \bm{\tilde{\xi}}_j & \colon i = 0
\\
\mathbf{r}_j & \colon i = 1
\end{cases}
\end{equation}
is a family of independent mappings
and the fact that
$ \forall \, j \in \{ 1, \ldots, k \} \colon
( \bm{\xi}_j - \bm{\tilde{\xi}}_j ) ( \mathbf{P} )
=
( \bm{\tilde{\xi}}_j - \bm{\xi}_j ) ( \mathbf{P} )
$
prove for all $p \in [ 1, \infty ) $ that
\begin{align}
& \int_{ \mathbf{\Omega} } \,
    \norm[\bigg]{ \sum_{j=1}^{k}
        \mathbf{r}_j( \bm{\omega} )
        \bigl[ \bm{\xi}_j( \bm{\omega} ) - \bm{\tilde{\xi}}_j( \bm{\omega} ) \bigr]
    }_E^p
\, \mathbf{P}( \dd \bm{\omega} ) \nonumber \\
& =
\int_{ \mathbf{\Omega} } \,
    \norm[\bigg]{ \sum_{j=1}^{k}
        \mathbbm{1}_{ \overline{
            \cup_{ i = 0 }^1 [ ( -1 )^i ( \bm{\xi}_j -  \bm{\tilde{\xi}}_j ) ] ( \mathbf{\Omega} )
        }^E }
        \bigl(
            \mathbf{r}_j( \bm{\omega} )
            \bigl[ \bm{\xi}_j( \bm{\omega} ) - \bm{\tilde{\xi}}_j( \bm{\omega} ) \bigr]
        \bigr) \,
        \mathbf{r}_j( \bm{\omega} )
        \bigl[ \bm{\xi}_j( \bm{\omega} ) - \bm{\tilde{\xi}}_j( \bm{\omega} ) \bigr]
    }_E^p
\, \mathbf{P}( \dd \bm{\omega} ) \nonumber \\
& =
\int_{ (\{ -1, 1 \} \times E )^k } \,
    \norm[\bigg]{ \sum_{j=1}^{k}
        \mathbbm{1}_{ \overline{
            \cup_{ i = 0 }^1 [ ( -1 )^i ( \bm{\xi}_j -  \bm{\tilde{\xi}}_j ) ] ( \mathbf{\Omega} )
        }^E }( z_j x_j ) \,
        z_j x_j
    }_E^p
\nonumber \\
& \quad
\bigl( ( \mathbf{r}_1, \bm{\xi}_1 -  \bm{\tilde{\xi}}_1, \ldots, \mathbf{r}_k, \bm{\xi}_k -  \bm{\tilde{\xi}}_k )( \mathbf{P} ) \bigr)
( \dd z_1, \dd x_1, \ldots, \dd z_k, \dd x_k ) \nonumber \\
& =
\int_{\{ -1, 1 \}}\int_E \cdots \int_{\{ -1, 1 \}}\int_E \,
    \norm[\bigg]{ \sum_{j=1}^{k}
        \mathbbm{1}_{ \overline{
            \cup_{ i = 0 }^1 [ ( -1 )^i ( \bm{\xi}_j -  \bm{\tilde{\xi}}_j ) ] ( \mathbf{\Omega} )
        }^E }( z_j x_j ) \,
        z_j x_j
    }_E^p
\label{eq:rand,1} \\
& \quad
\bigl( ( \bm{\xi}_k -  \bm{\tilde{\xi}}_k ) ( \mathbf{P} ) \bigr) ( \dd x_k )
\bigl( ( \mathbf{r}_k ) ( \mathbf{P} ) \bigr) ( \dd z_k )
\, \ldots \,
\bigl( ( \bm{\xi}_1 -  \bm{\tilde{\xi}}_1 ) ( \mathbf{P} ) \bigr) ( \dd x_1 )
\bigl( ( \mathbf{r}_1 ) ( \mathbf{P} ) \bigr) ( \dd z_1 )
\nonumber \\
& =
\int_E \cdots \int_E \,
    \norm[\bigg]{ \sum_{j=1}^{k}
        \mathbbm{1}_{ \overline{
            \cup_{ i = 0 }^1 [ ( -1 )^i ( \bm{\xi}_j -  \bm{\tilde{\xi}}_j ) ] ( \mathbf{\Omega} )
        }^E }( x_j ) \,
        x_j
    }_E^p
\nonumber \\
& \quad
\bigl( ( \bm{\xi}_k -  \bm{\tilde{\xi}}_k ) ( \mathbf{P} ) \bigr) ( \dd x_k )
\, \ldots \,
\bigl( ( \bm{\xi}_1 -  \bm{\tilde{\xi}}_1 ) ( \mathbf{P} ) \bigr) ( \dd x_1 )
\nonumber \\
& =
\int_{ E^k } \,
    \norm[\bigg]{ \sum_{j=1}^{k}
        \mathbbm{1}_{ \overline{
            \cup_{ i = 0 }^1 [ ( -1 )^i ( \bm{\xi}_j -  \bm{\tilde{\xi}}_j ) ] ( \mathbf{\Omega} )
        }^E }( x_j ) \,
        x_j
    }_E^p
\bigl( ( \bm{\xi}_1 -  \bm{\tilde{\xi}}_1, \ldots, \bm{\xi}_k -  \bm{\tilde{\xi}}_k ) ( \mathbf{P} ) \bigr) ( \dd x_1, \ldots, \dd x_k )
\nonumber  \\
& =
\int_{ \mathbf{\Omega} } \,
    \norm[\bigg]{ \sum_{j=1}^{k}
        \bigl[ \bm{\xi}_j( \bm{\omega} ) - \bm{\tilde{\xi}}_j( \bm{\omega} ) \bigr]
    }_E^p
\, \mathbf{P}( \dd \bm{\omega} ).  \nonumber
\end{align}
Furthermore, the fact that $ \sum_{j=1}^{k} \bm{\xi}_j \colon \mathbf{\Omega} \to E $
and $ \sum_{j=1}^{k} \bm{\tilde{\xi}}_j \colon \mathbf{\Omega} \to E  $
are independent,
the facts that
$ \int_{ \mathbf{\Omega} } \,
    \norm[\big]{ \sum_{j=1}^{k} \bm{\tilde{\xi}}_j( \bm{\omega} ) }_E
\, \mathbf{P}( \dd \bm{\omega} )
< \infty $
and
$ \int_{ \mathbf{\Omega} } \,
    \sum_{j=1}^{k} \bm{\tilde{\xi}}_j( \bm{\omega} )
\, \mathbf{P}( \dd \bm{\omega} )
= 0 $,
Lemma~\ref{lem:symm}, and \eqref{eq:rand,1} imply that it holds for all $p \in [ 1, \infty ) $ that
\begin{equation}
\begin{split}
& \norm[\Bigg]{ \sum_{j=1}^k \xi_j }_{ \mathscr{L}^p( \P; \left\| \cdot \right\|_E ) }
= \norm[\Bigg]{ \sum_{j=1}^k \bm{\xi}_j }_{ \mathscr{L}^p( \mathbf{P}; \left\| \cdot \right\|_E ) }
\leq  \norm[\Bigg]{ \sum_{j=1}^k  ( \bm{\xi}_j - \bm{\tilde{\xi}}_j ) }_{ \mathscr{L}^p( \mathbf{P}; \left\| \cdot \right\|_E ) } \\
& = \norm[\Bigg]{ \sum_{j=1}^k  \mathbf{r}_j ( \bm{\xi}_j - \bm{\tilde{\xi}}_j )  }_{ \mathscr{L}^p( \mathbf{P}; \left\| \cdot \right\|_E ) }
\leq  \norm[\Bigg]{ \sum_{j=1}^k  \mathbf{r}_j \bm{\xi}_j  }_{ \mathscr{L}^p( \mathbf{P}; \left\| \cdot \right\|_E ) }
        + \norm[\Bigg]{ \sum_{j=1}^k  \mathbf{r}_j \bm{\tilde{\xi}}_j  }_{ \mathscr{L}^p( \mathbf{P}; \left\| \cdot \right\|_E ) } \\
& = 2 \, \norm[\Bigg]{ \sum_{j=1}^k r_j \xi_j }_{ \mathscr{L}^p( \P; \left\| \cdot \right\|_E ) }.
\end{split}
\end{equation}
The proof of Lemma~\ref{lem:rand} is thus completed.
\end{proof}
The next result, Proposition~\ref{prop:type} below, is the key to estimating the statistical error term in the Banach space valued Monte Carlo method in the next subsection. Proposition~\ref{prop:type} is similar to, e.g., Proposition~9.11 in Ledoux \& Talagrand~\cite{lt91}.
\begin{proposition}[Sums of independent, centred, Banach space valued random variables] \label{prop:type}
Let $ k \in \mathbb{N} $, $ q \in [1,2] $,
let $ ( E, \norm{\cdot}_E ) $
be an $ \R $-Banach space with type $ q $,
let $ ( \Omega, \calF, \P ) $
be a probability space,
and let $ \xi_j \in \mathscr{L}^1( \P; \left\| \cdot \right\|_E ) $, $ j \in \{1, \ldots, k \} $,
be independent mappings which satisfy
for all  $ j \in \{1, \ldots, k \} $ that $ \E[ \xi_j ] = 0 $.
Then it holds for all $ p \in [q,\infty) $ that
\begin{equation} 
\norm[\Bigg]{ \sum_{j=1}^k \xi_j }_{ \mathscr{L}^p( \P; \left\| \cdot \right\|_E ) }
\leq   \varTheta_{ p, q }( E )  \Biggl(  \sum_{j=1}^k \norm{ \xi_j }_{ \mathscr{L}^p( \P; \left\| \cdot \right\|_E ) }^q \Biggr)^{\nicefrac{1}{q}}.
\end{equation}
\end{proposition}
\begin{proof}[Proof of Proposition~\ref{prop:type}]
Throughout this proof let $ ( \tilde{\Omega}, \tilde{\calF}, \tilde{\P} ) $
be a probability space,
let $ r_j \colon \tilde{\Omega} \to \{ -1, 1 \} $, $ j \in \{1, \ldots, k \} $,
be a $ \tilde{\P} $-Rademacher family,
and let $ \bm{\xi}_j \colon \Omega \times \tilde{\Omega} \to E $, $ j \in \{1, \ldots, k \} $,
and $ \mathbf{r}_j \colon \Omega \times \tilde{\Omega} \to \{ -1, 1 \} $, $ j \in \{1, \ldots, k \} $,
be the mappings which satisfy
for all $ \bm{\omega} = ( \omega, \tilde{\omega} ) \in \Omega \times \tilde{\Omega} $, $ j \in \{ 1, \ldots, k \} $
that $ \bm{\xi}_j( \bm{\omega} ) = \xi_j( \omega ) $
and $ \mathbf{r}_j( \bm{\omega} ) = r_j( \tilde{\omega} ) $.
Lemma~\ref{lem:rand} and the triangle inequality show for all $ p \in [ q, \infty ) $ that
\begin{equation}
\begin{split}
& \norm[\Bigg]{ \sum_{j=1}^k \xi_j }_{ \mathscr{L}^p( \P; \left\| \cdot \right\|_E ) }
= \norm[\Bigg]{ \sum_{j=1}^k \bm{\xi}_j }_{ \mathscr{L}^p( \P \otimes \tilde{\P}; \left\| \cdot \right\|_E ) }
\leq  2 \, \norm[\Bigg]{ \sum_{j=1}^k \mathbf{r}_j \bm{\xi}_j }_{ \mathscr{L}^p( \P \otimes \tilde{\P}; \left\| \cdot \right\|_E ) }
\\ & =
2 \,
\Biggl( \int_{ \Omega } \,
    \norm[\bigg]{ \sum_{j=1}^k r_j( \cdot ) \xi_j( \omega ) }_{  \mathscr{L}^p( \tilde{\P}; \left\| \cdot \right\|_E )  }^p
\P( \dd \omega ) \Biggr)^{ \nicefrac{1}{p} }
\\ & \leq
2 \mathscr{K}_{ p, q }
\Biggl( \int_{ \Omega } \,
    \norm[\bigg]{ \sum_{j=1}^k r_j( \cdot ) \xi_j( \omega ) }_{  \mathscr{L}^q( \tilde{\P}; \left\| \cdot \right\|_E )  }^p
\P( \dd \omega ) \Biggr)^{ \nicefrac{1}{p} }
\\
& \leq  2 \mathscr{K}_{ p, q }  \mathscr{T}_q( E )
           \norm[\Bigg]{  \Biggl( \sum_{j=1}^k \norm{ \xi_j }_E^q \Biggr)^{\nicefrac{1}{q}} }_{ \mathscr{L}^p( \P; \left| \cdot \right| ) }
=  2 \mathscr{K}_{ p, q }  \mathscr{T}_q( E )
               \norm[\Bigg]{  \sum_{j=1}^k \norm{ \xi_j }_E^q }_{ \mathscr{L}^{\nicefrac{p}{q}}( \P; \left| \cdot \right| ) }^{ \nicefrac{1}{q} } \\
& \leq  2 \mathscr{K}_{ p, q }  \mathscr{T}_q( E )
           \Biggl(  \sum_{j=1}^k \norm{ \xi_j }_{ \mathscr{L}^p( \P; \left\| \cdot \right\|_E ) }^q \Biggr)^{\nicefrac{1}{q}}.\end{split}
\end{equation}
This finishes the proof of Proposition~\ref{prop:type}.
\end{proof}
The result in Corollary~\ref{cor:sums} below is a direct consequence of Proposition~\ref{prop:type}.
\begin{corollary}[Sums of independent Banach space valued random variables] \label{cor:sums}
Let $ M \in \mathbb{N} $, $ q \in [1,2] $,
let $ ( E, \norm{\cdot}_E ) $
be an $ \R $-Banach space with type $ q $,
let $ ( \Omega, \calF, \P ) $
be a probability space,
and let $ \xi_j \in \mathscr{L}^1( \P; \left\| \cdot \right\|_E ) $, $ j \in \{1, \ldots, M \} $,
be independent.
Then it holds for all $ p \in [q,\infty) $ that
\begin{equation} 
\sigma_{p,E} \Biggl(  \sum_{j=1}^M \xi_j  \Biggr)
\leq   \varTheta_{ p, q }( E )  \Biggl(  \sum_{j=1}^M \abs{ \sigma_{p,E}( \xi_j ) }^q  \Biggr)^{\nicefrac{1}{q}}.
\end{equation}
\end{corollary}
\begin{corollary}[Monte Carlo methods in Banach spaces] \label{cor:BanachMC}
\label{cor:monteCarloBanach}
Let $ M \in \mathbb{N} $, $ q \in [1,2] $,
let $ ( E, \norm{\cdot}_E ) $
be an $ \R $-Banach space with type $ q $,
let $ ( \Omega, \calF, \P ) $
be a probability space,
and let $ \xi_j \in \mathscr{L}^1( \P; \left\| \cdot \right\|_E ) $, $ j \in \{1, \ldots, M \} $,
be independent and identically distributed.
Then it holds for all $ p \in [q,\infty) $ that
\begin{equation}
\norm[\Bigg]{  \E[ \xi_1 ] - \frac{1}{M} \sum_{j=1}^{M} \xi_j  }_{ \mathscr{L}^p( \P; \left\| \cdot \right\|_E ) }
=  \frac{ \sigma_{p,E} \bigl(  \sum_{j=1}^{M} \xi_j  \bigr) }{ M }
\leq  \frac{  \varTheta_{ p, q }( E )  \, \sigma_{p,E} ( \xi_1 )  }{ M^{ 1 - \nicefrac{1}{q} } }.
\end{equation}
\end{corollary}
Results on lower and upper error bounds related to Corollary~\ref{cor:BanachMC} can be found, e.g., in Theorem~1 in Daun \& Heinrich~\cite{DaunHeinrich2013} and in Corollary 2 in Heinrich \& Hinrichs~\cite{HeinrichHinrichs2014}.
Note that Corollary~\ref{cor:BanachMC} does not imply convergence if the underlying Banach space $ ( E, \norm{\cdot}_E ) $ has only type $1$, in the sense that it holds for all $ q \in ( 1, \infty ) $ that $ \mathscr{T}_q( E ) = \infty $.
\subsection{Multilevel Monte Carlo methods in Banach spaces} \label{sec:mlmc}
In many situations the work required to obtain a certain accuracy of an
approximation using the Monte Carlo method can
be reduced by using a multilevel Monte Carlo method. Heinrich
\cite{h98,h01}
was first to observe this and
established
multilevel Monte Carlo methods concerning
convergence in a Banach (function) space.
However, these methods do not apply to
SDEs.
Then
Giles~\cite{g08b} derived the complexity reduction
of multilevel Monte Carlo methods for
SDEs.
The minor contribution
of Proposition~\ref{prop:multilevelMC} below
to the literature on multilevel Monte Carlo methods
is to combine
the approaches of Heinrich~\cite{h98} and of Giles~\cite{g08b}
into a single result on multilevel Monte Carlo methods in Banach spaces.
The useful observation of Proposition~\ref{prop:multilevelMC} 
generalises the discussion in Section~4 of Heinrich~\cite{h01}.

\begin{proposition} [Abstract multilevel Monte Carlo methods in Banach spaces] \label{prop:multilevelMC}
Let 
$ q \in [1,2] $,
let
$
  ( 
    \Omega, \mathscr{F}, \P
  )
$
be a probability space,
let 
$ 
  ( 
    V_1, \left\| \cdot \right\|_{V_1}
  ) 
$
be an $ \R $-Banach space with 
type $ q $,
let 
$ 
  ( 
    V_2, \left\| \cdot \right\|_{V_2}
  ) 
$
be an $ \R $-Banach space with 
$ V_1 \subseteq V_2 $ 
continuously,
let $ v \in V_2 $, 
$ L \in \N $,
$ M_1, \dots, M_L \in \N $,
and for every
$ \ell \in \{ 1, \dots, L \} $
let
$ 
  D_{ \ell, k } \in \mathscr{L}^1( \P; \left\| \cdot \right\|_{ V_1 } )
$,
$ k \in \{ 1, \ldots, M_\ell \} $,
be independent and identically distributed.
Then it holds for all $ p \in [q,\infty) $ that
\begin{equation}  
\begin{split}
&  \left\|
    v -
    \sum_{ \ell = 1 }^{ L }
    \frac{1}{ M_\ell  }
      \sum_{ k = 1 }^{ M_\ell } 
      D_{ \ell, k }
  \right\|_{
    \mathscr{L}^p(\P; \left\| \cdot \right\|_{V_2} )
  } \\
& \leq
  \left\|
    v 
    -
    \sum_{ \ell = 1 }^L
    \E\big[
      D_{ \ell, 1 }
    \big]
  \right\|_{ V_2 }
+  
  \left\| 
    \operatorname{Id}_{V_1}
  \right\|_{
    L( V_1, V_2 )
  }
  \varTheta_{ p, q }( V_1 )
  \sum_{ \ell = 1 }^L
    \frac{ 
      \sigma_{ p, V_1 }( D_{ \ell, 1 } )
    }{
      (
        M_\ell
      )^{ 1 - \nicefrac{ 1 }{ q } }
    }.
\end{split}     
\end{equation}
\end{proposition}
\begin{proof}[Proof of Proposition~\ref{prop:multilevelMC}]
The triangle inequality
and Corollary~\ref{cor:BanachMC}
imply for all $ p \in [q,\infty) $ that
\begin{equation}
\begin{split}
&
  \left\|
    v -
    \sum_{ \ell = 1 }^{ L }
    \frac{1}{ M_\ell }
      \sum_{ k = 1 }^{ M_\ell } 
      D_{ \ell, k }
  \right\|_{
    \mathscr{L}^p(\P; \left\| \cdot \right\|_{V_2} )
  }
\\ & \leq
  \left\|
    v 
    -
    \sum_{ \ell = 1 }^L
    \E\!\left[
      D_{ \ell, 1 }
    \right]
  \right\|_{ V_2 }
  +
  \left\|
    \sum_{ \ell = 1 }^{ L }
    \E\!\left[
      D_{ \ell, 1 }
    \right]
    -
    \sum_{ \ell = 1 }^{ L }
    \frac{1}{M_\ell}
      \sum_{ k = 1 }^{ M_\ell } 
      D_{ \ell, k }
  \right\|_{
    \mathscr{L}^p(\P; \left\| \cdot \right\|_{V_2} )
  }
\\ & \leq
  \left\|
    v 
    -
    \sum_{ \ell = 1 }^L
    \E\!\left[
      D_{ \ell, 1 }
    \right]
  \right\|_{ V_2 } 
  +
  \left\|
    \operatorname{Id}_{ V_1 }
  \right\|_{
    L( V_1, V_2 )
  }
  \sum_{ \ell = 1 }^{ L }
    \left\|
      \E\!\left[
        D_{ \ell, 1 }
      \right]
      -
      \frac{1}{M_\ell}
        \sum_{ k = 1 }^{ M_\ell } 
        D_{ \ell, k }
    \right\|_{
      \mathscr{L}^p(\P; \left\| \cdot \right\|_{V_1} )
    }
\\ & \leq
  \left\|
    v 
    -
    \sum_{ \ell = 1 }^L
    \E\!\left[
      D_{ \ell, 1 }
    \right]
  \right\|_{V_2}
  +
  \left\|
    \operatorname{Id}_{ V_1 }
  \right\|_{
    L( V_1, V_2 )
  }
  \varTheta_{p,q}( V_1 )
  \sum_{ \ell = 1 }^L
    \frac{ 
      \sigma_{p,V_1}( D_{ \ell, 1 } )
    }{
      (M_\ell)^{1-\nicefrac{1}{q}}
    }.
\end{split}
\end{equation}
This completes the proof
of Proposition~\ref{prop:multilevelMC}.
\end{proof}

\begin{corollary}[Multilevel Monte Carlo methods in Banach spaces] \label{lem:montecarlo}
Consider the notation in Subsection~\ref{notation},
let $ q \in [ 1, 2 ] $,
$ L \in \N_0 $,
$ M_0, M_1, \dots, M_{L+1}, N_0, N_1, \dots, N_L \in \N $,
let
$
  \left( 
    \Omega, \mathscr{F}, \P
  \right)
$
be a probability space,
let 
$ 
  ( V_i, \left\| \cdot \right\|_{V_i} ) 
$,
$ i \in \{ 1, 2 \} $,
be separable $ \R $-Banach spaces
such that
$
  ( 
    V_1, \left\| \cdot \right\|_{ V_1 }
  )
$
has type $ q $ and such that
$
  V_1 \subseteq V_2
$
continuously,
let
$ ( V_3, \left\| \cdot \right\|_{V_3} ) $
be an $ \R $-Banach space,
let
$ 
  f \in \mathscr{ M }( \mathscr{B}( V_3 ), \mathscr{B}( V_2 ) )
$,
$ 
  g \in \mathscr{ M }( \mathscr{B}( V_3 ), \mathscr{B}( V_1 ) )
$,
$ X \in \mathscr{M}( \mathscr{ F }, \mathscr{ B }( V_3 ) ) $
satisfy
$
  \E \bigl[
    \norm{ f( X ) }_{ V_2 }
  \bigr]
  < \infty
$,
for every
$ n \in \N $
let
$
  Y^{ n, l, k } \in \mathscr{ M }( \mathscr{ F }, \mathscr{ B }( V_3 ) )
$,
$ k \in \N $,
$ l \in \N_0 $,
satisfy
$ 
  \E \bigl[
    \norm{ g( Y^{n,0,1} ) }_{ V_1 }
  \bigr]
  < \infty
$,
assume that
$ Y^{ N_0, 0, k } $, $ k \in \N $,
are independent and identically distributed,
and assume
for every $ \ell \in \N \cap [ 0, L ] $
that
$ ( Y^{ N_{ (\ell - 1) }, \mathfrak{l}, k }, Y^{ N_\ell, \mathfrak{l}, k } ) $,
$ k \in \N $,
$ \mathfrak{l} \in \N_0 $,
are independent and identically distributed.
Then it holds for all $ p \in [q, \infty ) $ that
\begin{align}
\label{eq:multimontecarlo}
&
  \left\|
    \E\!\left[ f(X) \right]
    -
    \tfrac{
      1
    }{ 
      M_0 
    }
      \sum_{ k = 1 }^{ M_0 }
      g( Y^{ N_0, 0, k } )
    -
    \sum_{ \ell = 1 }^{ L }
    \tfrac{
      1
    }{ M_\ell }
      \sum_{ k = 1 }^{ M_\ell } 
      \left[
        g( Y^{ N_\ell, \ell, k } ) -
        g( Y^{ N_{ (\ell - 1) }, \ell, k } ) 
      \right]
  \right\|_{
    \mathscr{L}^p(\P; \left\| \cdot \right\|_{V_2} )
  }
\nonumber \\
& \leq
  \left\|
    \E\!\left[ f(X) \right]
    -
    \E\!\left[ 
      g( Y^{N_L,0,1} ) 
    \right]
  \right\|_{ V_2 }
  \\&\quad\nonumber
  +
  \left\|
    \operatorname{Id}_{ V_1 }
  \right\|_{
    L( V_1, V_2 )
  }
    \varTheta_{ p, q }( V_1 )
  \biggl(
  \tfrac{ 
    \sigma_{p,V_1}\!\left(
      g( Y^{ N_0, 0, 1 } )
    \right)
  }{ (M_0)^{1-\nicefrac{1}{q}} }
  +
    \sum_{ \ell = 1 }^{ L }
    \tfrac{
      \sigma_{p,V_1}\!\left(
        g( Y^{ N_\ell, 0, 1 } ) -
        g( Y^{ N_{ (\ell - 1) }, 0, 1 } )
      \right)
    }{
      (M_\ell)^{1-\nicefrac{1}{q}}
    }
  \biggr)
\\ & \leq
\nonumber
  \left\|
    \E\!\left[ f(X) \right]
    -
    \E\!\left[ g( Y^{N_L,0,1} ) 
    \right]
  \right\|_{ V_2 }
  \\&
  \quad\nonumber
  +
  \left\|
    \operatorname{Id}_{ V_1 }
  \right\|_{
    L( V_1, V_2 )
  }
    \varTheta_{ p, q }( V_1 )
  \biggl(
  \tfrac{ 
    2 \, \left\|
      g( X )
    \right\|_{
      \mathscr{L}^p(\P; \left\| \cdot \right\|_{V_1} )
    }
  }{  (M_0)^{1-\nicefrac{1}{q}} }
  +
    \sum_{ \ell = 0 }^{ L }
    \tfrac{
      4 \, \left\| 
        g( Y^{ N_\ell, 0, 1 } ) - 
        g( X )
      \right\|_{
        \mathscr{L}^p(\P; \left\| \cdot \right\|_{V_1} )
      }
    }{
      (\min\{ M_\ell, M_{\ell+1} \})^{1-\nicefrac{1}{q}}
    }
  \biggr).
\end{align}
\end{corollary}
\begin{proof}[Proof
of Corollary~\ref{lem:montecarlo}]
Proposition~\ref{prop:multilevelMC}
and the identity
\begin{equation}
  \E\!\left[ g( Y^{N_L,0,1} ) 
  \right]
=
  \E\!\left[ 
    g( Y^{N_0,0,1} )
  \right]
  +
  \sum_{ \ell = 1 }^{ L } 
    \E\!\left[ 
      g( Y^{N_\ell,0,1} ) -
      g( Y^{N_{(\ell-1)},0,1} ) 
    \right]
\end{equation}
imply the first inequality
in \eqref{eq:multimontecarlo}.
Next note that
the triangle inequality implies
for all $ \xi \in \mathscr{L}^1( \P; \left\| \cdot \right\|_{ V_1 } ) $,
$ p \in [q, \infty ) $ that
$ \sigma_{p,V_1}( \xi ) \leq 2 \, \| \xi \|_{\mathscr{L}^p(\P;\left\| \cdot \right\|_{V_1})} $.
This and again the triangle inequality show for all $ p \in [q, \infty ) $ that
\begin{equation}
\begin{split}
&
\tfrac{ 
\sigma_{p,V_1}\!\left(
  g( Y^{ N_0, 0, 1 } )
\right)
}{ (M_0)^{1-\nicefrac{1}{q}} }
+
\sum_{ \ell = 1 }^{ L }
\tfrac{
  \sigma_{p,V_1}\!\left(
    g( Y^{ N_\ell, 0, 1 } ) -
    g( Y^{ N_{ (\ell - 1) }, 0, 1 } )
  \right)
}{
  (M_\ell)^{1-\nicefrac{1}{q}}
}
\\
& \leq
\tfrac{ 
2 \, \left\|
  g( Y^{ N_0, 0, 1 } )
\right\|_{\mathscr{L}^p(\P;\left\| \cdot \right\|_{V_1})}
}{ (M_0)^{1-\nicefrac{1}{q}} }
+
\sum_{ \ell = 1 }^{ L }
\tfrac{
2 \, \left\|
    g( Y^{ N_\ell, 0, 1 } ) -
    g( Y^{ N_{ (\ell - 1) }, 0, 1 } )
\right\|_{\mathscr{L}^p(\P;\left\| \cdot \right\|_{V_1})}
}{
  (M_\ell)^{1-\nicefrac{1}{q}}
}
\\
&\leq
\tfrac{
2 \, \left\|
  g( X )
\right\|_{\mathscr{L}^p(\P;\left\| \cdot \right\|_{V_1})}
\mathop{+}
2 \, \left\|
  g( Y^{ N_0, 0, 1 } ) - g(X)
\right\|_{\mathscr{L}^p(\P;\left\| \cdot \right\|_{V_1})}
}{ (M_0)^{1-\nicefrac{1}{q}} }
\\ & \quad
+
\sum_{ \ell = 1 }^{ L }
\tfrac{
2 \, \left\|
    g( Y^{ N_\ell, 0, 1 } ) -
    g( X )
\right\|_{\mathscr{L}^p(\P;\left\| \cdot \right\|_{V_1})}
  \mathop{+}
2 \, \left\|
    g( Y^{ N_{ (\ell - 1) }, 0, 1 } )
    -
    g( X )
\right\|_{\mathscr{L}^p(\P;\left\| \cdot \right\|_{V_1})}
}{
  (M_\ell)^{1-\nicefrac{1}{q}}
}
\\ &
\leq
\tfrac{ 
2 \, \left\|
  g( X )
\right\|_{\mathscr{L}^p(\P;\left\| \cdot \right\|_{V_1})}
}{ (M_0)^{1-\nicefrac{1}{q}} }
+
\sum_{ \ell = 0 }^{ L }
\tfrac{
 4 \, \left\|
    g( Y^{ N_\ell, 0, 1 } ) -
    g( X )
  \right\|_{\mathscr{L}^p(\P;\left\| \cdot \right\|_{V_1})}
}{
  (\min\{M_{\ell},M_{\ell+1}\})^{1-\nicefrac{1}{q}}
}.
\end{split}
\end{equation}
This
implies the second inequality in~\eqref{eq:multimontecarlo}.
The proof of Corollary~\ref{lem:montecarlo} is thus completed.
\end{proof}

\begin{corollary}[Convergence of
multilevel Monte Carlo approximations]
\label{c:mlmc.conv}
Consider the notation in Subsection~\ref{notation},
let $ T \in (0,\infty) $,
$ \beta \in (0,1] $,
$ 
  \alpha \in (0,\beta) 
$,
$
  c,r \in [0,\infty)
$,
let 
$ ( \Omega, \mathscr{F}, \P ) $
be a probability space,
let 
$ ( E, \norm{\cdot}_{E} ) $
be a separable $ \R $-Banach space with type $ 2 $,
let
$ 
  X \colon [0,T] \times \Omega
  \rightarrow E
$ 
be a stochastic
process with 
continuous sample paths
which satisfies for all $p\in[1,\infty)$, $\gamma\in[0,\beta)$ that
$
   X \in \calC^{\gamma}\bigl([0,T], \left\| \cdot \right\|_{  \mathscr{L}^p(\P; \left\| \cdot \right\|_E )  } \bigr)
$,
for every $ N \in \N $, $ \ell \in \N_0 $, $ k \in \N $
let 
$ Y^{ N, \ell, k } \colon
  [0,T] \times \Omega
  \rightarrow E $
be a stochastic process which satisfies
for all
$ 
  n \in \{ 0, 1, \dots, N-1 \} 
$,
$ 
  t \in 
  \big[ 
    \frac{ n T }{ N }, 
    \frac{ (n+1) T }{ N }
  \big] 
$,
$ p \in [1,\infty) $, $ \rho \in [0,\beta) $
that
\begin{align}
&  Y_{t}^{ N, \ell, k }
  = 
  \bigl(
    n + 1 - \tfrac{ t N }{ T }
  \bigr) \cdot
  Y_{ \frac{ n T }{ N } }^{ N, \ell, k } 
  +
  \bigl(
    \tfrac{ t N }{ T } - n
  \bigr) \cdot
  Y_{ \frac{ (n+1) T }{ N } }^{ N, \ell, k }, \\
& \sup_{ M \in \N }
  \sup_{ m\in\{0,1,\ldots,M\} }
  \Bigl(
    M^{ \rho } \,
      \norm[\big]{ 
        X_{ \frac{ m T }{ M } } - 
        Y_{ \frac{ m T }{ M } }^{ M, 0, 1 }
      }_{ 
        \mathscr{L}^p(\P; \left\| \cdot \right\|_E ) 
      }
  \Bigr)
  < \infty, \label{eq:strong.conv}
\end{align}
assume
for every $ N_1, N_2 \in \N $
that
$ ( Y^{ N_1, \ell, k }, Y^{ N_2, \ell, k } ) $, $ k \in \N $, $ \ell \in \N_0 $,
are independent and identically distributed,
and let
$
  f \colon 
  C( [0,T], E)
  \rightarrow
  C( [0,T], E)
$
be a
$ \mathscr{B}\bigl( C( [0,T], E) \bigr) / $\linebreak$ \mathscr{B}\bigl( C( [0,T], E) \bigr) $-measurable
function which satisfies
for all $ v, w \in \calC^{ \alpha }( [0,T], \left\| \cdot \right\|_E ) $
that
$ f(v), f(w) \in \calC^{\alpha}( [ 0, T ], \left\| \cdot \right\|_E ) $ and
\begin{equation}
\label{eq:f.local.lipschitz}
  \norm{ f(v) - f(w) }_{
    \calC^{ \alpha }( [0,T], \left\| \cdot \right\|_E )
  }
\leq
  c \,
  \Bigl( 
    1 
    + 
    \norm{v}_{ 
      \calC^{ \alpha }( [0,T], \left\| \cdot \right\|_E ) 
    }^r
    +
    \norm{w}_{ 
      \calC^{ \alpha }( [0,T], \left\| \cdot \right\|_E ) 
    }^r
  \Bigr)
  \norm{ v - w }_{ 
    \calC^{ \alpha }( [0,T], \left\| \cdot \right\|_E ) 
  }.
\end{equation}
Then
it holds that
\begin{equation}
\E \Bigl[ \norm{ f(X) }_{ \calC^{ \alpha }( [0,T], \left\| \cdot \right\|_E ) } \Bigr]
<\infty,
\end{equation}
it holds for all
$ p \in [1,\infty) $, $ \rho \in [0,\beta-\alpha) $
that
\begin{equation}
\sup_{ N \in \N } 
\bigg[
N^{ \rho }
\Bigl( \E \Bigl[
\norm{
  f(X) - f(Y^{N,0,1})
}_{ 
    \calC^{ \alpha }( [0,T], \left\| \cdot \right\|_E )
}^p
\Bigr] \Bigr)^{\nicefrac{1}{p}}
\bigg]
 < \infty,
\end{equation}
and it holds
for all $p\in[1,\infty)$, $\gamma\in[0,\alpha) $, $\rho \in [0, \beta-\alpha)$
that
\begin{equation}
\label{eq:mlmc.conv}
\begin{multlined}[c][0.903\textwidth]
\sup_{L\in\N}
\left[
2^{L \cdot \min\{ \rho, \nicefrac{1}{2} \} }
L^{ -\mathbbm{1}_{ \{ \nicefrac{1}{2} \} }( \rho ) }
\Biggl| \E \Biggl[
\biggl\|
\E[ f( X ) ]
-
 {\textstyle \sum\limits_{ k = 1 }^{ 2^L } }
  \tfrac{
f( Y^{ 1 , 0, k } )
  }{ 2^L }
\vphantom{\Biggr|^{\nicefrac{1}{p}}}
\right.
\\
\left.
-
{\textstyle \sum\limits_{ \ell = 1 }^{L}   }
    {\textstyle \sum\limits_{ k = 1 }^{ 2^{L-\ell} } }
\tfrac{
 f( Y^{ 2^\ell , \ell, k } ) 
 -
 f( Y^{ 2^{(\ell-1)} , \ell, k } ) 
}{
2^{L-\ell}
}
\biggr\|_{ \calC^{ \gamma }( [0,T], \left\| \cdot \right\|_E )  }^p
\Biggr] \Biggr|^{\nicefrac{1}{p}}
\right]
<\infty.
\hspace{-1em}
\end{multlined}
\end{equation}
\end{corollary}
\begin{proof}[Proof of Corollary~\ref{c:mlmc.conv}]
Throughout this proof
let
$ \gamma \in [ 0, \alpha ) $,
$ \delta \in ( \gamma, \frac{3\gamma + \alpha}{4} ) $,
let $ C^1( [ 0, T ], E ) $
be the $ \R $-vector space of continuously Fr\'echet differentiable functions
from $ [ 0, T ] $ to $ E $,
let $ \norm{\cdot}_{ C^1( [ 0, T ], E ) } \colon C^1( [ 0, T ], E ) \to [ 0, \infty ) $
be the function which satisfies
for all $ v \in C^1( [ 0, T ], E ) $ that
$ \norm{v}_{ C^1( [ 0, T ], E ) } = \norm{v}_{ C( [ 0, T ], \left\| \cdot \right\|_E ) } + \norm{v'}_{ C( [ 0, T ], \left\| \cdot \right\|_E ) } $,
let
$ \mathcal{W}^{ \nicefrac{(\alpha + \gamma)}{2}, \nicefrac{4}{(\alpha - \gamma)} }( [ 0, T ], E ) $
be the Sobolev space
with regularity parameter $ \nicefrac{(\alpha + \gamma)}{2} \in ( 0, 1 ) $
and integrability parameter $ \nicefrac{4}{(\alpha - \gamma)} \in ( 4, \infty ) $
of continuous functions from $ [ 0, T ] $ to $ E $,
let
\begin{equation*}
\norm{\cdot}_{ \mathcal{W}^{ \nicefrac{(\alpha + \gamma)}{2}, \nicefrac{4}{(\alpha - \gamma)} }( [ 0, T ], E ) } \colon
\mathcal{W}^{ \nicefrac{(\alpha + \gamma)}{2}, \nicefrac{4}{(\alpha - \gamma)} }( [ 0, T ], E ) \to [ 0, \infty )
\end{equation*}
be the function which satisfies
for all $ v \in \mathcal{W}^{ \nicefrac{(\alpha + \gamma)}{2}, \nicefrac{4}{(\alpha - \gamma)} }( [ 0, T ], E ) $ that
\begin{equation}
\norm{v}_{ \mathcal{W}^{ \nicefrac{(\alpha + \gamma)}{2}, \nicefrac{4}{(\alpha - \gamma)} }( [ 0, T ], E ) }
= \!
\Biggl[
\int_{0}^{T}
    \norm{ v( t ) }_E^{ \frac{4}{\alpha - \gamma} }
\ud t
\mathop{+} \!
\int_{0}^{T}
\int_{0}^{T}
    \frac{ \norm{ v( t ) - v( s ) }_E^{ \frac{4}{\alpha - \gamma} } }{ \abs{ t - s }^{ \frac{ 3 \alpha + \gamma }{ \alpha - \gamma } } }
\ud t
\ds
\Biggr]^{\frac{\alpha - \gamma}{4}},
\end{equation}
let $ V_1, V_2 \subseteq \calC^\gamma( [ 0, T ], \left\| \cdot \right\|_E ) $ be the sets given by
$ V_1 = \mathcal{W}^{ \nicefrac{(\alpha + \gamma)}{2}, \nicefrac{4}{(\alpha - \gamma)} }( [ 0, T ], E ) $
and
\begin{equation}
V_2
=
\Biggl\{
    v \in \calC^\gamma( [ 0, T ], \left\| \cdot \right\|_E ) \colon
    \limsup_{ n \to \infty } \sup_{ \substack{ s, t \in [ 0, T ], \: 0 < \abs{ s - t } < \nicefrac{1}{n} } }
    \frac{ \norm{ v( s ) - v( t ) }_E }{ \abs{ s - t }^\gamma } = 0
\Biggr\}
\end{equation}
(cf., e.g., Lunardi~\cite[Section~0.2]{Lunardi1995reprint}),
let
$ \norm{\cdot}_{V_1} \colon V_1 \to [ 0, \infty ) $
be the function given by
$ \norm{\cdot}_{V_1} =  \norm{\cdot}_{ \mathcal{W}^{ \nicefrac{(\alpha + \gamma)}{2}, \nicefrac{4}{(\alpha - \gamma)} }( [ 0, T ], E ) } $,
let
$ \norm{\cdot}_{V_2} \colon V_2 \to [ 0, \infty ) $
be the function which satisfies
for all $ v \in V_2 $ that
$ \norm{v}_{V_2}
=
\norm{v}_{ \calC^\gamma( [ 0, T ], \left\| \cdot \right\|_E ) } $,
let $ ( V_3,  \norm{\cdot}_{V_3} ) $
be the $ \R $-Banach space
given by
\begin{align}
( V_3, \norm{\cdot}_{V_3} ) & = \Bigl( C( [ 0, T ], E ), \norm{\cdot}_{ C( [ 0, T ], \left\| \cdot \right\|_E ) }
                    \text{\raisebox{-0.1\baselineskip}{$\bigr|_{ C( [ 0, T ], E ) }$}}  \Bigr),
\end{align}
and
let $ \mathfrak{f} \colon V_3 \to V_2 $
and
$ g \colon V_3 \to V_1 $
be the functions which satisfy
for all $ v \in V_3 $ that
$ \mathfrak{f}(v) = g(v) = \mathbbm{1}_{ \calC^\alpha( [ 0, T ], \left\| \cdot \right\|_E ) }(v) f(v) $.
Observe that
the Kolmogorov-Chentsov continuity theorem (see Theorem~\ref{thm:Kolmogorov})
together with the assumptions that
$
   X \in \cap_{p\in[1,\infty)}\cap_{\eta\in[0,\beta)}
   \calC^{\eta}\bigl([0,T], \left\| \cdot \right\|_{  \mathscr{L}^p(\P; \left\| \cdot \right\|_E )  } \bigr)
$
and that $ X $ has continuous sample paths
implies
for all $ p \in [ 1, \infty ) $
that
$
\E \bigl[
    \norm{ X }_{  \calC^{\alpha}([0,T], \left\| \cdot \right\|_E) }^p
\bigr]
< \infty
$.
This,
assumption~\eqref{eq:f.local.lipschitz},
H\"older's inequality,
and
Corollary~\ref{cor:hoelder3}
show
for all
$ p \in [1,\infty) $,
$ \rho \in [0,\beta-\alpha) $
that
\begin{align}
& \sup_{ N \in \N } 
  \Bigl(
    N^{ \rho } \,
      \E \Bigl[
    \norm{
      \mathfrak{f}(X) - g(Y^{N,0,1})
    }_{ 
        \calC^{ \alpha }( [0,T], \left\| \cdot \right\|_E )
    }
   \Bigr]
  \Bigr) \nonumber \\
& \leq \sup_{ N \in \N } 
  \bigg[
    N^{ \rho }
      \Bigl( \E \Bigl[
    \norm{
      f(X) - f(Y^{N,0,1})
    }_{ 
        \calC^{ \alpha }( [0,T], \left\| \cdot \right\|_E )
    }^p
   \Bigr]
  \Bigr)^{\nicefrac{1}{p}} \bigg] \nonumber
  \\&
  \nonumber
  \leq
  \sup_{ N \in \N } 
  \bigg[
    N^{ \rho }
    \Bigl( \E \Bigl[ \Bigl(
  c \,
  \bigl( 
    1 
    + 
    \norm{ 
      X
    }_{ 
      \calC^{ \alpha }( [0,T], \left\| \cdot \right\|_E ) 
    }^r
    +
    \norm{ 
      Y^{N,0,1} 
    }_{ 
      \calC^{ \alpha }( [0,T], \left\| \cdot \right\|_E ) 
    }^r
  \bigr)
\\ & \label{eq:Calpha.fX-fY}  \qquad\qquad\qquad\quad\, \cdot
  \norm{ 
    X - Y^{N,0,1} 
  }_{ 
    \calC^{ \alpha }( [0,T], \left\| \cdot \right\|_E ) 
  }
   \Bigr)^p \, \Bigr] \Bigr)^{ \nicefrac{1}{p} }
  \bigg]
  \\& \nonumber
  \leq
  c \,
  \biggl[
  1+
  \Bigl( \E \Bigl[
  \norm{
    X
  }_{ \calC^{\alpha}([0,T],\left\| \cdot \right\|_E)  }^{2pr}
  \Bigr] \Bigr)^{ \nicefrac{1}{(2p)} }
  +
  \sup_{ N \in \N }
  \Bigl( \E \Bigl[
  \norm{
    Y^{N,0,1}
  }_{ \calC^{\alpha}([0,T],\left\| \cdot \right\|_E) }^{2pr}
  \Bigr] \Bigr)^{ \nicefrac{1}{(2p)} }
  \biggr]
  \\&
  \qquad
  \cdot
  \sup_{ N \in \N } 
  \bigg[
  N^{ \rho }
  \Bigl( \E \Bigl[
  \norm{
    X - Y^{N,0,1} 
  }_{ \calC^{ \alpha }( [0,T], \left\| \cdot \right\|_E ) }^{2p}
  \Bigr] \Bigr)^{ \nicefrac{1}{(2p)} }
  \bigg]
  < \infty.
  \nonumber
\end{align}
Assumption~\eqref{eq:f.local.lipschitz} also ensures for all $ p \in [ 1, \infty ) $ that
\begin{align}
& \E \Bigl[ \norm{ f(X) }_{ \calC^{ \alpha }( [0,T], \left\| \cdot \right\|_E ) } \Bigr]
\leq
\Bigl( \E \Bigl[
    \norm{ f( X ) }_{ \calC^{ \alpha }( [0,T], \left\| \cdot \right\|_E ) }^p
\Bigr] \Bigr)^{ \nicefrac{1}{p} }
\label{eq:Calpha.fX} \\
& \leq
\norm{ f(0) }_{ \calC^{ \alpha }( [0,T], \left\| \cdot \right\|_E ) }
+
c \,
\biggl[
\Bigl( \E \Bigl[
    \norm{ X }_{ \calC^{ \alpha }( [0,T], \left\| \cdot \right\|_E ) }^p
\Bigr] \Bigr)^{ \nicefrac{1}{p} }
+
\Bigl( \E \Bigl[
    \norm{ X }_{ \calC^{ \alpha }( [0,T], \left\| \cdot \right\|_E ) }^{(r+1)p}
\Bigr] \Bigr)^{ \nicefrac{1}{p} }
\biggr]
< \infty. \nonumber
\end{align}
Next note that $ ( V_1, \norm{\cdot}_{V_1} ) $
is a separable $ \R $-Banach space with type 2.
In addition,
the fact that
$ ( C^1( [ 0, T ], E ), \norm{\cdot}_{C^1( [ 0, T ], E )} ) $
is a separable $ \R $-Banach space,
the fact that
$ C^1( [ 0, T ], E ) \subseteq \calC^\gamma( [ 0, T ], \left\| \cdot \right\|_E ) $
continuously,
and the fact that
\begin{equation}
 \overline{C^1( [ 0, T ], E )}^{ \calC^\gamma( [ 0, T ], \left\| \cdot \right\|_E ) } = V_2 
\end{equation}
(cf., e.g., Lunardi~\cite[Proposition~0.2.1]{Lunardi1995reprint})
prove that
$ ( V_2, \norm{\cdot}_{V_2} ) $
is a separable $ \R $-Banach space.
Moreover,
the Sobolev embedding theorem proves that
$ V_1 \subseteq \calC^\delta( [ 0, T ], \left\| \cdot \right\|_E ) $ continuously.
This and the fact that
$ \calC^\delta( [ 0, T ], \left\| \cdot \right\|_E ) \subseteq V_2 $ continuously
establish that
$ V_1 \subseteq V_2 $ continuously.
Combining~\eqref{eq:Calpha.fX} with~\eqref{eq:Calpha.fX-fY} and the fact that
$ \calC^{\alpha}([ 0, T ], \left\| \cdot \right\|_E ) \subseteq V_1 $ continuously
hence implies for all $ p \in [1,\infty) $, $ \rho \in [0,\beta-\alpha) $ that
$ \E \bigl[ \norm{ \mathfrak{f}(X) }_{V_2} \bigr]
+ \sup_{ N \in \N } \E \bigl[ \norm{ g(Y^{N,0,1}) }_{V_1} \bigr] < \infty $,
$ \norm{ g( X ) }_{ \mathscr{L}^p(\P; \left\| \cdot \right\|_{V_1} ) } < \infty $, and
\begin{equation} \label{eq:convergence}
\sup_{ N \in \N } \Bigl( N^{ \rho } \,
                            \E \Bigl[ \norm{ \mathfrak{f}( X ) - g( Y^{N,0,1} ) }_{V_2} \Bigr] \Bigr)
+ \sup_{ N \in \N } \Bigl( N^{ \rho } \,
                            \norm{ g( X ) - g( Y^{N,0,1} ) }_{ \mathscr{L}^p(\P; \left\| \cdot \right\|_{V_1} ) } \Bigr)
                            < \infty.
\end{equation}
Furthermore, observe that it holds
for all $ L \in \N $, $ \rho \in [ 0, \beta-\alpha ) \setminus \{ \tfrac{1}{2} \} $ that
\begin{equation} \label{eq:sum}
\begin{split}
\sum_{\ell=1}^{L} (2^{\ell})^{-\rho} \, 2^{-\frac12(L-\ell)}
& = 2^{-\frac{L}{2}}
      \sum_{\ell=1}^{L} 2^{(\frac{1}{2}-\rho)\ell}
= 2^{-\frac{L}{2}}
   \tfrac{ 1-2^{(\frac12-\rho)L} }{ 2^{\rho - \frac{1}{2}} - 1 }
= 2^{- L \cdot \min \{ \rho, \frac{1}{2} \} }
          \tfrac{ 1-2^{-|\frac12-\rho| L} }{ {\displaystyle |} 1-2^{\rho-\frac{1}{2}} {\displaystyle |} }
\leq \tfrac{ 2^{- L \cdot \min \{ \rho, \frac{1}{2} \} } }{ {\displaystyle |} 1-2^{\rho-\frac{1}{2}} {\displaystyle |} }
\end{split}
\end{equation}
and
\begin{equation} \label{eq:sum2}
\begin{split}
\sum_{\ell=1}^{L} (2^{\ell})^{-\frac12} \, 2^{-\frac12(L-\ell)}
& =
2^{-\frac{L}{2}} L
.
\end{split}
\end{equation}
Combining Corollary~\ref{lem:montecarlo} with \eqref{eq:convergence}, \eqref{eq:sum}, and \eqref{eq:sum2} implies \eqref{eq:mlmc.conv}.
This finishes the proof of Corollary~\ref{c:mlmc.conv}.
\end{proof}
Corollary~\ref{c:mlmc.conv} can be applied to many SDEs.
Under general conditions on the coefficient functions of the SDEs
(see, e.g., Theorem~1.3 and Subsection~3.1 in~\cite{HutzenthalerJentzen2014}),
suitable stopped-tamed Euler approximations (cf.\ (6) in~\cite{HutzenthalerJentzenWang2013} or (10) in~\cite{HutzenthalerJentzenKloeden2012})
converge in the strong sense with convergence rate $\nicefrac{1}{2}$.
We note that the classical Euler-Maruyama approximations do not satisfy
condition~\eqref{eq:strong.conv} for most SDEs with superlinearly growing coefficients;
see Theorem 2.1 in~\cite{hjk11} and Theorem 2.1 in~\cite{HutzenthalerJentzenKloeden2013}.
Moreover, under general conditions on the coefficients
it holds that the solution process is strongly $\nicefrac{1}{2}$-H\"older continuous
in time.
In conclusion, provided that a suitable numerical scheme is employed,
Corollary~\ref{c:mlmc.conv} can be applied to many SDEs
with $\beta=\nicefrac{1}{2}$.

\subsection*{Acknowledgements}
This project has been partially supported
by the research project
``Numerical approximations of stochastic differential equations with non-globally Lipschitz continuous coefficients'' funded by the German Research Foundation,
by the ETH Research Grant \mbox{ETH-47 15-2}
``Mild stochastic calculus and numerical approximations for nonlinear sto\-chastic evolution equations with L\'evy noise'',
and by the project
``Construction of New Smoothness Spaces on Domains''
(project number I 3403)
funded by the Austrian Science Fund (FWF).

\printbibliography

\end{document}